\documentstyle[txmac,a4,%leqno,%
amssymb,%
bibenc,%
case,%
twoside,%
nocaphead2,%
epsf,%
leqno,%
myrot,%
mypic,%
times,mathptm%
]{article}

\advance\oddsidemargin by -1.9cm
\advance\evensidemargin by -1.9cm
\advance\textwidth by 3.8cm
\def\mynewtheo#1#2{%
\newtheorem{@#1}{#2}[section]%
\newenvironment{#1}{\begin{@#1}\rm}{\end{@#1}}}

\mynewtheo{lemma}{Lemma}
\mynewtheo{exer}{Exercise}
\mynewtheo{theo}{Theorem}
\mynewtheo{rem}{Remark}
\mynewtheo{defi}{Definition}
\mynewtheo{conj}{Conjecture}
\mynewtheo{corr}{Corollary}
\mynewtheo{prop}{Proposition}
\mynewtheo{question}{Question}
\mynewtheo{exam}{Example}

\makeatletter

\newenvironment{theorem}{\begin{theo}}{\end{theo}}
\newenvironment{conjecture}{\begin{conj}}{\end{conj}}

\begin{document}

% typing `equation' breaks my fingers!
\newenvironment{eqn}{\begin{equation}}{\end{equation}\@ignoretrue}

% eqlabel=(#1)
\newenvironment{myeqn*}[1]{\begingroup\def\@eqnnum{\reset@font\rm#1}%
\xdef\@tempk{\arabic{equation}}\begin{equation}\edef\@currentlabel{#1}}
{\end{equation}\endgroup\setcounter{equation}{\@tempk}\ignorespaces}

% eqlabel=#1
\newenvironment{myeqn}[1]{\begingroup\let\eq@num\@eqnnum
\def\@eqnnum{\bgroup\let\r@fn\normalcolor % an extremely UGLY hack !!!
\def\normalcolor####1(####2){\r@fn####1#1}%
%\show\reset@font
\eq@num\egroup}%
\xdef\@tempk{\arabic{equation}}\begin{equation}\edef\@currentlabel{#1}}
{\end{equation}\endgroup\setcounter{equation}{\@tempk}\ignorespaces}

% eqlabel=(eqnnr) \qed
\newenvironment{myeqn**}{\begin{myeqn}{(%\arabic{equation}
%\show\theequation
\theequation)\es\es\mbox{\qed}}\edef\@currentlabel{\theequation}}
{\end{myeqn}\stepcounter{equation}}

\renewenvironment{eqn}{\begin{equation}}{\end{equation}\ignorespaces}

\def\rottab#1#2#3{
\expandafter\advance\csname c@#3\endcsname by -1\relax
\centerline{%
%\hfill
%\fboxsep0pt\relax
%\fbox{%
\rbox{\centerline{\vbox{\setbox1=\hbox{#1}%
\hbox to \wd1{\hfill\vbox{{%
%\tracingmacros1
\caption{#2}}}\hfill}%
\vskip6mm
\box1}}%
}
%}
}%
}

\newcommand{\mybin}[2]{\text{$\Bigl(\begin{array}{@{}c@{}}#1\\#2%
\end{array}\Bigr)$}}
\newcommand{\mybinn}[2]{\text{$\biggl(\begin{array}{@{}c@{}}%
#1\\#2\end{array}\biggr)$}}

\def\overtwo#1{\mbox{\small$\mybin{#1}{2}$}}
\newcommand{\mybr}[2]{\text{$\Bigl\lfloor\mbox{%
\small$\displaystyle\frac{#1}{#2}$}\Bigr\rfloor$}}
\def\mybrtwo#1{\mbox{\mybr{#1}{2}}}

\def\myfrac#1#2{\raisebox{0.2em}{\small$#1$}\!/%
\!\raisebox{-0.2em}{\small$#2$}}
\def\ffrac#1#2{\mbox{\small$\ds\frac{#1}{#2}$}}

\def\noloop{{\diag{0.5cm}{0.5}{1}{\picline{0.25 0}{0.25 1}}}}

\def\vrt#1{{\picfillgraycol{0}\picfilledcircle{#1}{0.09}{}}}

\def\ReidI#1#2{
  \diag{0.5cm}{0.9}{1}{
    \pictranslate{0.4 0.5}{
      \picscale{#11 #21}{
        \picmultigraphics[S]{2}{1 -1}{
	  \picmulticurve{-6 1 -1 0}{0.5 -0.5}{0.5 0}{0.1 0.3}{-0.2 0.3}
	} %picmultigraphics
	\piccirclearc{-0.2 0}{0.3}{90 270}
      }%picscale
    }%pictranslate
  }%diag
}
\def\Pos#1#2{{\diag{#1}{1}{1}{#2
\picmultiline{-5 1 -1 0}{0 1}{1 0}
\picmultiline{-5 1 -1 0}{0 0}{1 1}
}}}
\def\Neg#1#2{{\diag{#1}{1}{1}{#2
\picmultiline{-5 1 -1 0}{0 0}{1 1}
\picmultiline{-5 1 -1 0}{0 1}{1 0}
}}}
\def\Nul#1#2{{\diag{#1}{1}{1}{#2
\piccirclearc{0.5 1.4}{0.7}{-135 -45}
\piccirclearc{0.5 -0.4}{0.7}{45 135}
}}}
\def\Inf#1#2{{\diag{#1}{1}{1}{#2
\piccirclearc{0.5 1.4 x}{0.7}{135 -135}
\piccirclearc{0.5 -0.4 x}{0.7}{-45 45}
}}}
\def\pos{\Pos{0.5em}{\piclinewidth{10}}}
\def\neg{\Neg{0.5em}{\piclinewidth{10}}}
\def\nul{\Nul{0.5em}{\piclinewidth{10}}}
\def\inf{\Inf{0.5em}{\piclinewidth{10}}}

\def\noloop{{\diag{0.5cm}{0.5}{1}{\picline{0.25 0}{0.25 1}}}}

\def\ReidI#1#2{
  \diag{0.5cm}{0.9}{1}{
    \pictranslate{0.4 0.5}{
      \picscale{#11 #21}{
        \picmultigraphics[S]{2}{1 -1}{
	  \picmulticurve{-6 1 -1 0}{0.5 -0.5}{0.5 0}{0.1 0.3}{-0.2 0.3}
	} %picmultigraphics
	\piccirclearc{-0.2 0}{0.3}{90 270}
      }%picscale
    }%pictranslate
  }%diag
}

\def\GD{\szCD{6mm}}
\def\szCD#1#2{{\let\@nomath\@gobble\small\diag{#1}{2.4}{2.4}{
  \picveclength{0.27}\picvecwidth{0.1}
  \pictranslate{1.2 1.2}{
    \piccircle{0 0}{1}{}
    #2
}}}}
\def\CD{\szCD{4mm}}

\def\point#1{{\picfillgraycol{0}\picfilledcircle{#1}{0.08}{}}}
\def\labpt#1#2#3{\pictranslate{#1}{\point{0 0}\picputtext{#2}{$#3$}}}
\def\vrt#1{{\picfillgraycol{0}\picfilledcircle{#1}{0.09}{}}}

\def\chrd#1#2{\picline{1 #1 polar}{1 #2 polar}}
\def\arrow#1#2{\picvecline{1 #1 polar}{1 #2 polar}}

\def\labch#1#2#3{\chrd{#1}{#2}\picputtext{1.3 #2 polar}{$#3$}}
\def\labar#1#2#3{\arrow{#1}{#2}\picputtext{1.3 #2 polar}{$#3$}}
\def\labbr#1#2#3{\arrow{#1}{#2}\picputtext{1.3 #1 polar}{$#3$}}

\def\@dcont{}
\def\svCD#1{\ea\glet\csname #1\endcsname\@dcont}
\def\rsCD#1{\ea\glet\ea\@dcont\csname #1\endcsname\ea\glet
\csname #1\endcsname\relax}

\def\addCD#1{\ea\gdef\ea\@dcont\ea{\@dcont #1}}
\def\drawCD#1{\szCD{#1}{\@dcont}}

\def\epsfs#1#2{{\catcode`\_=11\relax\ifautoepsf\unitxsize#1\relax\else
\epsfxsize#1\relax\fi\epsffile{#2.eps}}}
\def\epsfsv#1#2{{\vcbox{\epsfs{#1}{#2}}}}
\def\vcbox#1{\setbox\@tempboxa=\hbox{#1}\parbox{\wd\@tempboxa}{\box
  \@tempboxa}}
\newbox\@tempboxb
\def\vtbox#1{\setbox\@tempboxa=\hbox{\shortstack[t]{#1}}%
  \setbox\@tempboxb=\hbox{2rut}\@tempdima\ht\@tempboxa\relax
  \advance\@tempdima-\ht\@tempboxb\relax\raise-\@tempdima\box \@tempboxa
  %\dp\@tempboxa=\@tempdima
  %\box \@tempboxa
  }
\def\p{\epsfsv{2cm}}

\def\@test#1#2#3#4{%
  \let\@tempa\go@
  \@tempdima#1\relax\@tempdimb#3\@tempdima\relax\@tempdima#4\unitxsize\relax
  \ifdim \@tempdimb>\z@\relax
    \ifdim \@tempdimb<#2%
      \def\@tempa{\@test{#1}{#2}}%
    \fi
  \fi
  \@tempa
}

\def\go@#1\@end{}
\newdimen\unitxsize
\newif\ifautoepsf\autoepsftrue

\unitxsize4cm\relax
\def\epsfsize#1#2{\epsfxsize\relax\ifautoepsf
  {\@test{#1}{#2}{0.1 }{4   }
		{0.2 }{3   }
		{0.3 }{2   }
		{0.4 }{1.7 }
		{0.5 }{1.5 }
		{0.6 }{1.4 }
		{0.7 }{1.3 }
		{0.8 }{1.2 }
		{0.9 }{1.1 }
		{1.1 }{1.  }
		{1.2 }{0.9 }
		{1.4 }{0.8 }
		{1.6 }{0.75}
		{2.  }{0.7 }
		{2.25}{0.6 }
		{3   }{0.55}
		{5   }{0.5 }
		{10  }{0.33}
		{-1  }{0.25}\@end
		\ea}\ea\epsfxsize\the\@tempdima\relax
		\fi
		}

\def\vrt#1{{\picfillgraycol{0}\picfilledcircle{#1}{0.09}{}}}
\def\fdot#1{{\picfillgraycol{0}\picfilledcircle{#1}{0.02}{}}}
\def\cycl#1#2#3#4{\vrt{#1}\vrt{#2}\vrt{#3}\vrt{#4}
\picline{#1}{#2}\picline{#2}{#3}\picline{#3}{#4}\picline{#4}{#1}}
\def\xcycl#1#2#3#4#5#6{\vrt{#1}\vrt{#2}\vrt{#3}\vrt{#4}\vrt{#5}
\vrt{#6}\picline{#1}{#2}\picline{#2}{#3}\picline{#3}{#4}
\picline{#4}{#5}\picline{#5}{#6}\picline{#6}{#1}}

%usage ODD # of  coords, then {}
\def\@curvepath#1#2#3{%
  \@ifempty{#2}{\piccurveto{#1 }{@stc}{@std}#3}%
    {\piccurveto{#1 }{#2 }{#2  #3  0.5 conv}
    \@curvepath{#3}}%
}
\def\curvepath#1#2#3{%
  \piccurve{#1 }{#2 }{#2 }{#2  #3  0.5 conv}%
  \picPSgraphics{/@stc [ #1  #2  -1 conv ] $ D /@std [ #1  ] $ D }%
  \@curvepath{#3}%
}

%usage EVEN # of  coords, then {}
\def\@opencurvepath#1#2#3{%
  \@ifempty{#3}{\piccurveto{#1 }{#1 }{#2 }}%
    {\piccurveto{#1 }{#2 }{#2  #3  0.5 conv}\@opencurvepath{#3}}%
}
\def\opencurvepath#1#2#3{%
  \piccurve{#1 }{#2 }{#2 }{#2  #3  0.5 conv}%
  \@opencurvepath{#3}%
}

\author{A. Stoimenow\footnotemark[1]\\[2mm]
\small Research Institute for Mathematical Sciences, \\
\small Kyoto University, Kyoto 606-8502, Japan\\
\small e-mail: {\tt stoimeno@kurims.kyoto-u.ac.jp}\\
\small WWW: {\hbox{\web|http://www.kurims.kyoto-u.ac.jp/~stoimeno/|}}
}

{\def\thefootnote{\fnsymbol{footnote}}
\footnotetext[1]{Supported by 21st Century COE Program.}
}

\title{\large\bf \uppercase{Properties of closed 3-braids}\\
[1mm]\uppercase{and other link braid representations}\\
[4mm]
\phantom{\small\it This is a preprint. I would be grateful
for any comments and corrections %but I prefer you not circulate it
!}}

\date{\large Current version: \curv\ \ \ First version:
\makedate{22}{9}{2004}}

\maketitle

\makeatletter

\let\vn\varnothing
\let\point\pt
\let\ay\asymp
\let\pa\partial
\let\ap\alpha
\let\bt\beta
\let\be\beta
\let\Dl\Delta
\let\Gm\Gamma
\let\gm\gamma
\let\de\delta
\let\dl\delta
\let\eps\epsilon
\let\lm\lambda
\let\Lm\Lambda
\let\sg\sigma
\let\vp\varphi
\let\zt\zeta
\let\om\omega
\let\diagram\diag
\let\nb\nabla
\let\prt\partial
\let\wh\widehat
\let\wt\widetilde
\def\tW{\wt W}
\let\sm\setminus
\let\tl\tilde
\def\inx{\mathop {\operator@font ind}}
\def\spn{\mathop {\operator@font span}\nolimits}
\def\dig{\mathop {\operator@font diag}}
\def\Mc{\max\cf}
\def\Md{\max\deg}
\def\md{\min\deg}
\def\mc{\min\cf}
\def\vol{\text{\rm vol}\,}
\def\Ra{\Rightarrow}
\def\Lra{\Longrightarrow}
\def\lra{\longrightarrow}
\def\so{\Rightarrow}
\def\So{\Longrightarrow}
\def\nin{\not\in}
\let\ds\displaystyle
\let\llra\longleftrightarrow
\let\reference\ref
\let\ol\overline
\let\ul\underline
\let\u\underline
\let\h\hat
\let\es\enspace
\def\lfra{\leftrightarrow}
\def\cf{\text{\rm cf}\,}

\def\TM{$^\text{\raisebox{-0.2em}{${}^\text{TM}$}}$}
\def\ssim{\stackrel{\ds \sim}{\vbox{\vskip-0.2em\hbox{$\scriptstyle
*$}}}}

\long\def\@makecaption#1#2{%
   % \tm
   \vskip \abovecaptionskip 
   {\let\label\@gobble
   \let\ignorespaces\@empty
   \xdef\@tempt{#2}%
   }%
   \ea\@ifempty\ea{\@tempt}{%
   \sbox\@tempboxa{%
      \fignr#1#2}%
      }{%
   \sbox\@tempboxa{%
      {\fignr#1:}\capt\ #2}%
      }%
   \ifdim \wd\@tempboxa >\captionwidth {%
      %\centerline{\parbox{\captionwidth}{\unhbox \@tempboxa}}%
      \rightskip=\@captionmargin\leftskip=\@captionmargin
      \unhbox\@tempboxa\par
     }%
   \else
      \centerline{\box \@tempboxa}%
      % \hbox to\captionwidth{\hfil\box\@tempboxa\hfil}%
   \fi
   \vskip \belowcaptionskip
   }%
\def\fignr{\small\sffamily\bfseries}%
\def\capt{\small\sffamily}%

% \long\def\@makecaption#1#2{%
%    % \tm
%    \vskip 10pt
%    {\let\label\@gobble
%    \let\ignorespaces\@empty
%    \xdef\@tempt{#2}%
%    %\typeout{`#2'}%
%    }%
%    \ea\@ifempty\ea{\@tempt}{%
%    \setbox\@tempboxa\hbox{%
%       \fignr#1#2}%
%       }{%
%    \setbox\@tempboxa\hbox{%
%       {\fignr#1:}\capt\ #2}%
%       }%
%    \ifdim \wd\@tempboxa >\captionwidth {%
%       \rightskip=\@captionmargin\leftskip=\@captionmargin
%       \unhbox\@tempboxa\par}%
%    \else
%       \hbox to\captionwidth{\hfil\box\@tempboxa\hfil}%
%    \fi}%
% %
% \def\fignr{\small\sffamily\bfseries}%
% \def\capt{\small\sffamily}%

\newdimen\@captionmargin\@captionmargin2cm\relax
\newdimen\captionwidth\captionwidth0.8\hsize\relax

\def\eqref#1{(\protect\ref{#1})}

\def\proof{\@ifnextchar[{\@proof}{\@proof[\unskip]}}
\def\@proof[#1]{\noindent{\bf Proof #1.}\enspace}

\def\hint{\noindent Hint: }
\def\problem{\noindent{\bf Problem.} }

\def\@mt#1{\ifmmode#1\else$#1$\fi}
\def\qed{\hfill\@mt{\Box}}
\def\qqed{\hfill\@mt{\Box\enspace\Box}}

\def\cU{{\cal U}}
\def\cC{{\cal C}}
\def\cP{{\cal P}}
\def\fg{{\frak g}}
\def\tr{\text{tr\,}}
\def\cZ{{\cal Z}}
\def\cD{{\cal D}}
\def\bR{{\Bbb R}}
\def\cE{{\cal E}}
\def\bZ{{\Bbb Z}}
\def\bN{{\Bbb N}}

\def\bysame{\same[\kern2cm]\,}

\def\br#1{\left\lfloor#1\right\rfloor}
\def\BR#1{\left\lceil#1\right\rceil}

\def\abstractname{}

\@addtoreset {footnote}{page}

\renewcommand{\section}{%
   \@startsection
         {section}{1}{\z@}{-2.5ex \@plus -1ex \@minus -.2ex}%
               {3ex \@plus.2ex}{\Large\bf}%
}

\renewcommand{\subsubsection}{%
   \@startsection
         {subsubsection}{1}{\z@}{-1.5ex \@plus -1ex \@minus -.2ex}%
               {1ex \@plus.2ex}{\large\bf}%
}
\renewcommand{\@seccntformat}[1]{\csname the#1\endcsname .
\quad}

\def\bC{{\Bbb C}}
\def\bP{{\Bbb P}}

\makeatletter

\let\old@tl\~\def\~{\raisebox{-0.8ex}{\tt\old@tl{}}}
\let\lra\longrightarrow
\let\sm\setminus
\let\eps\varepsilon
\let\ex\exists
\let\fa\forall
\let\ps\supset

\def\rs#1{\raisebox{-0.4em}{$\big|_{#1}$}}

{\let\@noitemerr\relax
\vskip-2.7em\kern0pt\begin{abstract}
\noindent{\bf Abstract.}\enspace
We show that 3-braid links with given (non-zero) Alexander
or Jones polynomial are finitely many, and can be effectively
determined. We classify among closed 3-braids strongly
quasi-positive and fibered ones, and show that 3-braid links have
a unique incompressible Seifert surface. We also classify the
positive braid words with Morton-Williams-Franks
bound 3 and show that closed positive braids of braid index 3 are
closed positive 3-braids. For closed braids on more
strings, we study the alternating links occurring. In particular
we classify those of braid index 4, and show that their Morton-%
Williams-Franks inequality is exact. Finally, we use the Burau
representation to obtain new braid index criteria, including an
efficient 4-braid test.\\[2mm]
\\[2mm]
(secondary)
\end{abstract}
}\vspace{7mm}

{\parskip0.2mm\tableofcontents}
\vspace{7mm}

\section{Introduction}

Originating from the pioneering work of Alexander
\cite{Alexander2} and Artin \cite{Artin,Artin2}, braid
theory has become intrinsically interwoven with knot theory,
and over the years, braid representations of different types
have been studied, many of them with motivation coming
from fields outside of knot theory, for example dynamical
systems \cite{Williams} or 4-dimensional QFTs \cite{Kreimer}.

The systematic study of closed 3-braids was begun by Murasugi
\cite{Murasugi2}. Later 3-braid links have been classified
in \cite{BirMen}, in the sense that 3-braid representations
of the same link are exactly described. This is, up to
a few exceptions, mainly just conjugacy. The conjugacy
problem of 3-braids has a series of solutions, starting with
Schreier's algorithm \cite{Sc}, going over Garside
\cite{Garside} (for arbitrary braid groups), Xu \cite{Xu}
(and later more generally Birman-Ko-Lee \cite{BKL}), until,
for example, a recent algorithm of linear complexity due to
Fiedler and Kurlin \cite{FK}. So Birman-Menasco's work allows
to decide if two 3-braids have the same closure link.

However, many properties of links (except achirality and
invertibility) are not evident from (3-strand) braid
representations, and thus to classify 3-braid links with
special properties remains a non-trivial task. A first result
was given in \cite{Murasugi} for rational links, and then
improved in \cite{3br}, where we completed this project for
alternating links.

{}From the opposite point of view, a natural question concerning
braid representations of links is: 

\begin{question}\label{q1}
If a braid representation of a particular type exists, does
also one exist with the minimal number of strands (among all
braid representations of the link)?
\end{question}

The minimal number of
strands for a braid representation of a link $L$ is called the
\em{braid index} of the link, and will be denoted by $b(L)$.

\begin{figure}[h]
\[
% \begin{array}{c}
\epsfs{5cm}{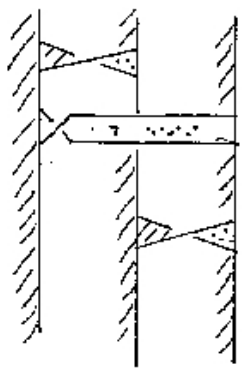}
% \end{array}
\]\\[-14mm]
\caption{A braided Seifert surface obtained from a band representation
\label{figbd}}
\end{figure}

In \cite{Bennequin}, Bennequin studied in relation to contact
structures braid representations by bands, independently considered
by Rudolph \cite{Rudolph} in the context of algebraic curves, and
more recently in \cite{BKL} from a group theoretic point of view.
A band representation naturally spans a Seifert surface of the
link. See figure \reference{figbd}.
 
% In figure \reference{fig2}, for the move (a) $\to$ (b), we slide
% $B$ and $C$, respectively, along $D$ and $A$, then delete $D$
% (by deplumbing a Hopf band), and then slide $C$ back along
% $A$. For the move (b) $\to$ (c), we slide $C$ along $B$,
% then slide $B$ along $C$. The $k$ bands $A$ can be subsequently
% removed by Murasugi desumming a $(2,k)$-torus link fiber surface.
% Thus the surface $F$, spanned by $\hat\bt$, turns after
% (de)summing Hopf bands into a surface $F'$ consisting of $k$
% copies of $F'_0$ and one negative band $N$. See figure \ref{fig3} (a).
%  
% \begin{figure}[h]
% \[
% \begin{array}{cc}
%  (a) & (b) \\
% \end{array}
% \]
% \caption{\label{fig3}}
% \end{figure}
%  
% In figure \reference{fig3}, we perform the move (a) $\to$ (b), sliding
% $B$ along $N$. The last surface $F'$, in figure \reference{fig3} (b),
% is Murasugi sum of a fibre surface spanning the $(2,2,-2)$-pretzel
% link and $\tl F$, where $\tl F$ consists 
% of $k - 1$ copies of $F'_0$ and the band $N$. By induction on
% $k$, we see that $F'$ is a fibre surface, and hence $\hat\bt$
% is a fibred link. \qed
%  
% \end{appendix}
Bennequin proved in \cite{Bennequin} that a minimal genus band
representation always exists on 3-braids. Such surface was called
by Birman-Menasco \cite{BirMen2} a Bennequin surface. Later, by
examples found by Ko and Lee, and Hirasawa and myself, the
existence of Bennequin surfaces on the minimal number of strands
was found not to extend to 4-braids. (We will see in this paper,
though, among many other things, that it does extend for
alternating links.) Since by Rudolph's work \cite{Rudolph},
a Bennequin surface exists for any link, on a possibly large
number of strands, the answer to the question \reference{q1}
is negative for minimal genus band representations. 

A special case of band representations are the positive ones,
called strongly quasi-positive in \cite{Rudolph}. For such
representations the above examples do not apply. The existence
of minimal such representations is a question of Rudolph, whose
answer is not known. Our first result is a positive answer to
Rudolph's question for 3-braids. 

\begin{theorem}\label{thqp}
If $L$ is a strongly quasi-positive link, and has braid index 3,
then $L$ is has a strongly quasi-positive 3-braid representation.
\end{theorem}

This result is an application of the work in \cite{3br} on
Xu's normal form \cite{Xu}. It legitimizes the term ``strongly
quasi-positive 3-braid link'', which otherwise would be
linguistically sloppy, for it may mean a 3-braid link that is strongly
quasi-positive, or the link (which is the closure) of a strongly
quasi-positive 3-braid.

As an (\em{a-priori}-so-seeming) improvement of Bennequin's
theorem, Birman-Menasco proved \cite{BirMen2} that
for 3-braid links not only some but in fact \em{any} minimal genus
surface is a Bennequin surface. Xu used this result and the
classification in \cite{BirMen} to conclude that most 3-braid links
have a unique minimal genus surface. Unfortunately, she failed to
deal (in the oriented sense) with some exceptions, which we complete
in section \reference{SA} (see theorem \reference{wu}).

The following work involves a further study of Xu's normal form.
% A proof of this result is given by the authors in an appendix.
%
A natural relation of this form to the skein polynomial was
exhibited in \cite{3br}. Here we extend this study to the Alexander
polynomial $\Dl$. As a result, we can describe 3-braid links with
given $\Dl$, and in particular show no non-trivial 3-braid
knot has trivial polynomial. In \cite{Birman}, Birman constructed
some pairs of different 3-braid links with the same polynomials
and proposed the problem to understand the values of the various
link polynomials on 3-braid links. Our work solves (to a large
extent at least) Birman's problem for the Alexander polynomial.
We also have a solution for the Jones and Brandt-Lickorish-%
Millett-Ho $Q$ polynomial. 

Parallel to the Alexander polynomial description, and using
a recent result of Hirasawa-Murasugi, we classify among 3-braid
links the fibered ones. (This result was proved independently
by Yi Ni in \cite{Ni}, apparently unaware of this work, which
I did to make widely public immediately.) We will see then
that a 3-braid knot is fibered if (and only if) its Alexander
polynomial is monic, but that this is \em{not} true for 2- or
3- component links.

The fact that most 3-braid links were found fibered,
and this could be proved independently from Bennequin and
Birman-Menasco, opened the hope for a more natural approach to
some of their results. Finally we were indeed able to complete this
project, with the assistance of M. Hirasawa and M. Ishiwata, obtaining

\begin{theorem}\label{thuq}
Any 3-braid link has a unique (non-closed) incompressible surface.
\end{theorem}

This result subsumes all previous uniqueness theorems. Besides,
its proof is entirely different, in that it fully avoids the contact
geometry approach of Bennequin and the considerable complexity
in Birman-Menasco's subsequent braid foliation work. Instead,
we use the sutured manifold theory of Gabai \cite{Gabai,Gabai2}
(to deal with the fibered links) and Kobayashi-Kakimizu \cite{%
Kobayashi,Kakimizu} (for the remaining cases), as well as some
of our preceding work on the Alexander polynomial (to rule out
disconnected Seifert surfaces). Since Kobayashi-Kakimizu use
a slightly stronger notion of uniqueness (see the beginning of
\S\reference{Uq}), theorem \reference{thuq} is an improvement even
for minimal genus surfaces. Moreover, its proof underscores the
geometric meaning of Xu's form, that remained unclear in \cite{BirMen}.
While we still cannot provide an ``easy'' proof of Birman-Menasco's
classification theorem (see remark on p.~34 of \cite{BirMen}), our
work will likely lead to simplifications in their very lengthy
treatment, which makes heavy use of incompressible surfaces.

The study and applications of Xu's form occupy the first main
part of the paper. Then some related results are included,
whose treatment requires different methods. These
results are dealt with in the later sections of the paper.

As a first such topic, more substantial effort will be made to
prove an analogous result to theorem \reference{thqp} for
positive (in the usual Artin generators) braid representations.

\begin{theorem}\label{thg}
If $L$ is a positive braid link, and has braid index 3, then 
$L$ is a positive 3-braid link.
\end{theorem}

So our work here can be understood to answer question
\reference{q1} for braid index 3 and strongly quasi-positive
or positive braid links. (The state of quasi-positive links,
and of quasi-positive braid representations with regard to 
question \reference{q1}, in contrast remains open. There is,
however, an algorithm, found by orevkov \cite{Orevkov}
to decide if a 3-braid is quasi-positive.)

The proof of theorem \reference{thg} requires a study of positive
braids where the braid index bound in the Morton-Williams-Franks
inequality \cite{Morton,WilFr} is 3. This study builds on and extends
(but also considerably simplifies) the work of Nakamura
\cite{Nakamura} for Morton-Williams-Franks bound 2. Then some
detailed case distinction and calculation are necessary.
In \cite{posex_bcr} it was shown, again, that theorem
\reference{thg} is not true for 4-braids: there are two 16 crossing
knots that have positive 16 crossing 5-braid representations,
but braid index 4 (and consequently only non-positive 4-braid
representations). It is in fact for these examples that I was
led to investigate about the theorem.

In section \ref{Spos},
%we come back to theorem \reference{thg}
% and show how to obtain its alternative proof
%. As a partial step, which is equally
% noteworthy,
using the arguments in the proof of theorem \reference{thg}
and a criterion of Yokota \cite{Yokota}, we consider 3-braid
representations of links that are (diagrammatically) positive.
We obtain strong restrictions on such representations,
and in particular determine which of the links are not fibered.
(See theorem \ref{thposq}.)

% We will show two completely different approaches to
% theorem \reference{thg}. The first uses the
% Morton-Williams-Franks inequality and tries to identify
% the positive braids with Morton-Williams-Franks bound 3,
% extending the work of Nakamura \cite{Nakamura} for
% bound 2. The proof is slightly technical, and requires
% the verification of certain special cases.
% Another proof, following a longer, but much
% more insightful road, is then given by continuing
% the arguments that establish theorem \reference{thqp}.

Another problem proposed by Birman, in \cite{MortonPb},
was to understand the relation of alternation of links
(as a pre-eminent diagram defined property) and braid
representations. Early substantial results on the braid
index of alternating links were due to Murasugi \cite{Murasugi},
who determined the braid index for rational and fibered
alternating links. In \cite{3br} we used Xu's form to classify
the braid index 3 alternating links. Here, in \S\reference{Sfg},
we will easily recover this result and push it forward to braid
index 4. For this we use an argument based on the Jones polynomial,
and connect the celebrated Kauffman-Murasugi-Thistlethwaite work
\cite{Kauffman2,Murasugi3,Thistle2} to braid representations.
The existence of the Bennequin surface of alternating links on a
4 string braid (corollary \reference{cfT}) is among several easy
consequences we obtain.

The last section \S\reference{SBurau} discusses some applications
of the Squier unitarization \cite{Squier} of the Burau
representation, that concern also braids on more strings. A
first series of conditions are estimates on the norm of special
values of the Alexander and Jones polynomial. They relate to and
partly extend estimates of Jones in \cite{Jones}. In particular, we
will see that the Alexander polynomial provides conditions for every
given braid index. Some properties of the skein polynomial, and a
relation to Mahler measure are also discussed.

In the case of a 4-braid, one can go further and (almost)
identify the Eigenvalues of the Burau matrix
from the Alexander and Jones polynomial (of its closure).
Using this, a criterion for braid index 4 is derived.
We show examples exhibiting the efficiency of this test, including
such where not only the Morton-Williams-Franks inequality itself,
but also its 2-cable version fails (and so our test seems the only
practicable option).

% and is included here to complete the picture.

In the appendix we collect some work of Hirasawa, Ishiwata and
Murasugi, which completes the proof of several of our results.

\section{Preliminaries, basic definitions and conventions}
% \noindent{\bf Convention.} 

Basic concepts that appear throughout the paper are summarized.
`W.l.o.g.' will abbreviate `without loss of generality';
`r.h.s.' will stand for `right hand-side'.

\subsection{Links and link diagrams\label{Sdg}}

Links are represented by diagrams; we assume diagrams are
oriented (though sometimes orientation is not relevant).

A crossing $p$ in a link diagram $D$ is called \em{reducible} (or
nugatory) if it looks like on the left of figure \reference{figred}.
$D$ is called reducible if it has a reducible crossing, else it is
called \em{reduced}. The reducing of the reducible crossing $p$
is the move depicted on figure \reference{figred}. Each diagram $D$
can be (made) reduced by a finite number of these moves.

\begin{figure}[htb]
\begin{eqn}\label{eqred}
\diag{6mm}{4}{2}{
  \picrotate{-90}{\rbraid{-1 2}{1 1.4}}
  % \piccirclearc{2 0.5}{1.3}{45 135}
  % \piccirclearc{2 1.5}{1.3}{-135 -45}
  \picputtext[u]{2 0.7}{$p$}
  \picscale{1 -1}{
    \picfilledcircle{0.7 -1}{0.8}{$P$}
  }
  \picfilledcircle{3.3 1}{0.8}{$Q$}
} \qquad\lra\qquad
\diag{6mm}{4}{2}{
  \piccirclearc{2 0.5}{1.3}{45 135}
  \piccirclearc{2 1.5}{1.3}{-135 -45}
  \picfilledcircle{1 1}{0.8}{$P$}
  \picfilledcircle{3 1}{0.8}{$Q$}
} 
\end{eqn}
\caption{\label{figred}}
\end{figure}

We assume in the following all diagrams reduced, unless otherwise
stated.

The diagram on the right of figure \reference{figtan}
is called \em{connected sum} $A\# B$ of the diagrams $A$ and $B$.
If a diagram $D$ can be represented as the connected sum of 
diagrams $A$ and $B$, such that both $A$ and $B$ have at least one
crossing, then $D$ is called \em{disconnected} (or composite), else
it is called \em{connected} (or prime). $K$ is \em{prime} if whenever
$D=A\# B$ is a composite diagram of $K$, one of $A$ and $B$
represent an unknotted arc (but not both; the unknot is
not considered to be prime per convention).

By $c(D)$ we denote the
\em{number of crossings} of $D$, $n(D)$ the number of components
of $D$ (or $K$, $1$ if $K$ is a knot), and $s(D)$ the \em{number
of Seifert circles} of $D$. The \em{crossing number} $c(K)$ of a
knot or link $K$ is the minimal crossing number of all diagrams
$D$ of $K$.
$!D$ is the \em{mirror image} of $D$, and $!K$ is the mirror image
of $K$. Clearly $g(!D)=g(D)$ and $g(!K)=g(K)$.

\begin{figure}[htb]
\[
\diag{6mm}{3}{2}{
  \piccirclearc{1.8 1}{0.5}{-120 120}
  \picfilledcircle{1 1}{0.8}{$A$}
}\, \#\,
\diag{6mm}{3}{2}{
  \piccirclearc{1.2 1}{0.5}{60 300}
  \picfilledcircle{2 1}{0.8}{$B$}
}\quad =\quad
\diag{6mm}{4}{2}{
  \piccirclearc{2 0.5}{1.3}{45 135}
  \piccirclearc{2 1.5}{1.3}{-135 -45}
  \picfilledcircle{1 1}{0.8}{$A$}
  \picfilledcircle{3 1}{0.8}{$B$}
} 
\]
\caption{\label{figtan}}
\end{figure}

The \em{Seifert graph} $\Gm(D)$ of a diagram $D$ is defined to
be the graph whose vertices are the Seifert circles of $D$ and
whose edges are the crossings. The \em{reduced Seifert graph}
$\Gm'(D)$ is defined by removing multiple copies of an edge
between two vertices in $\Gm(D)$, so that a simple edge
remains. (See \cite{MurPrz,MurPrz2} for example.)

The \em{(Seifert) genus} $g(K)$ resp.\ \em{Euler characteristic}
$\chi(K)$ of a knot or link $K$ is said to be the minimal genus
resp.\ maximal Euler characteristic of Seifert surface of $K$.
For a diagram $D$ of $K$, $g(D)$ is defined to be the genus
of the Seifert surface obtained by Seifert's algorithm on $D$,
and $\chi(D)$ its Euler characteristic. We have $\chi(D)=s(D)-c(D)$
and $2g(D)=2-n(D)-\chi(D)$.

The numbering of knots we use is as in the tables of \cite
[appendix]{Rolfsen} for prime knots of crossing number $\le 10$,
and as in \cite{KnotScape} for those of crossing number $11$
to $16$. KnotScape's numbering is reorganized so that for
given crossing number non-alternating knots are appended after
alternating ones, instead of using `a' and `n' superscripts.

\subsection{Polynomial link invariants\label{brD}}

Let $X\in\bZ[t,t^{-1}]$. The \em{minimal} or \em{maximal degree} 
$\md X$ or $\Md X$ is the minimal resp.\ maximal exponent of $t$ 
with non-zero coefficient in $X$. Let $\spn_tX=\Md_tX-\md_tX$.
The coefficient in degree $d$ of $t$ in $X$ is denoted $[X]_{t^d}$
or $[X]_{d}$. The \em{leading coefficient} $\Mc\,X$ of $X$ is its
coefficient in degree $\Md X$. If $X\in\bZ[x_1^{\pm 1},x_2^{\pm 1}]$,
then $\Md_{x_1}X$ denotes the maximal degree in $x_1$. Minimal
degree and coefficients are defined similarly, and of course
$[X]_{x_1^k}$ is regarded as a polynomial in $x_2^{\pm 1}$.

Let $P(v,z)$ be the \em{skein polynomial} \cite{HOMFLY,LickMil}.
It is a Laurent polynomial in two variables of oriented knots and
links. We use here the convention of \cite{Morton}, i.e. with the
polynomial taking the value $1$ on the unknot, having the variables
$v$ and $z$ and satisfying the skein relation
\begin{eqn}\label{1}
v^{-1}\,P\bigl(
\diag{5mm}{1}{1}{
\picmultivecline{-5 1 -1.0 0}{1 0}{0 1}
\picmultivecline{-5 1 -1.0 0}{0 0}{1 1}
}
\bigr)\,-\,
v \,P\bigl(
\diag{5mm}{1}{1}{
\picmultivecline{-5 1 -1 0}{0 0}{1 1}
\picmultivecline{-5 1 -1 0}{1 0}{0 1}
}
\bigr)\,=\,
z\,P\bigl(
\diag{5mm}{1}{1}{
\piccirclevecarc{1.35 0.5}{0.7}{-230 -130}
\piccirclevecarc{-0.35 0.5}{0.7}{310 50}
}
\bigr)\,.
\end{eqn}

We will denote in each triple as in \eqref{1} 
the diagrams (from left to right) by $D_+$, $D_-$ and $D_0$.
For a diagram $D$ of a link $L$, we will use all of the notations
$P(D)=P_D=P_D(l,m)=P(L)$ etc.\ for its skein polynomial, with
the self-suggestive meaning of indices and arguments. So we
can rewrite \eqref{1} as
\begin{eqn}\label{srel}
v^{-1}P_+(v,z)-vP_-(v,z) = zP_0(v,z)\,.
\end{eqn}

The \em{writhe} is a number ($\pm1$), assigned to any crossing in a
link diagram. A crossing as on the left in \eqref{1} has writhe 1 and
is called \em{positive}. A crossing as in the middle of \eqref{1} has
writhe $-1$ and is called \em{negative}. The writhe of a link diagram
is the sum of writhes of all its crossings.

Let $c_{\pm}(D)$ be the number of positive, respectively negative
crossings of a diagram $D$, so that $c(D)=c_+(D)+c_-(D)$ and
$w(D)=c_+(D)-c_-(D)$. 

The \em{Jones polynomial} \cite{Jones2} $V$, and (one variable)
\em{Alexander polynomial} \cite{Alexander} $\Dl$ are obtained
from $P$ by variable substitutions
\begin{eqn}\label{DP}
\Dl(t)=P(1,t^{1/2}-t^{-1/2})\,,
\end{eqn}
and
\begin{eqn}\label{VPsub}
V(t)=P(t,t^{1/2}-t^{-1/2})\,.
\end{eqn}

Hence these polynomials also satisfy corresponding skein relations.
(In algebraic topology, the Alexander polynomial is usually
defined only up to units in $\bZ[t,t^{-1}]$; the present
normalization is so that $\Dl(t)=\Dl(1/t)$ and $\Dl(1)=1$.)

In very contrast to its relatives, the \em{range} of the 
Alexander polynomial (i.e., set of values it takes) is known.
Let us call a polynomial $\Dl\in\bZ[t^{1/2},t^{-1/2}]$
\em{admissible} if it satisfies for some natural number
$n\ge 1$ the three properties\\
\mbox{}\ \ (i)\quad $t^{(n-1)/2}\Dl\in \bZ[t^{\pm 1}]$,\\
\mbox{}\ (ii)\quad  $\Dl(t)=(-1)^{n-1}\Dl(1/t)$ and \\
\mbox{}(iii)\quad   $(t^{1/2}-t^{-1/2})^{n-1}\mid
                    \Dl$ for $n>1$, or $\Dl(1)=1$ for $n=1$. \\
It is well-known that these are exactly the polynomials that occur
as (1-variable) Alexander polynomials of some $n$-component link. 

The \em{Kauffman polynomial} \cite{Kauffman} $F$ is usually defined
via a regular isotopy invariant $\Lm(a,z)$ of unoriented links.
We use here a slightly different convention for the variables
in $F$, differing from \cite{Kauffman,Thistle2} by the interchange
of $a$ and $a^{-1}$. Thus in particular we have for a link diagram $D$
the relation $F(D)(a,z)=a^{w(D)}\Lm(D)(a,z)$, where $\Lm(D)$ is the
writhe-unnormalized
version of the polynomial, given in our convention by the properties
\begin{eqnarray}
& \Lambda\bigl(\Pos{0.5cm}{}\bigr)\ +\ \Lambda\bigl(\Neg{0.5cm}{}\bigr)\ =\ z\ \bigl(\ 
\Lambda\bigl(\Nul{0.5cm}{}\bigr)\ +\ \Lambda\bigl(\Inf{0.5cm}{}\bigr)\ \bigr)\,,\label{wseven}\\[2mm]
& \label{wseven.5} \Lambda\bigl(\ \ReidI{-}{-}\bigr) = a^{-1}\ \Lambda\bigl(\noloop\bigr);\quad
\Lambda\bigl(\ \ReidI{-}{ }\bigr) = a\ \Lambda\bigl(\noloop\bigr)\,,\\[2mm]
& \label{wseven.7} \Lambda\bigl(\,\mbox{\Large $\bigcirc$}\,\bigr) = 1\,.
\end{eqnarray}

The \em{Brandt-Lickorish-Millett-Ho polynomial} \cite{BLM} $Q$
is given by $Q(z)=F(1,z)$.

\subsection{Semiadequacy\label{Ssaq}}

An alternative description of $V$ is given by the Kauffman bracket in
\cite{Kauffman2}. We do not need this description directly, but
for self-containedness is it useful to recall the related concept of
semiadequacy that was popularized in \cite{LickThis}.

\begin{figure}[htb]
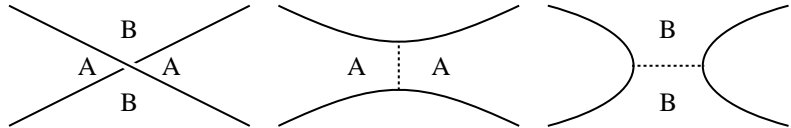

\[
\diag{8mm}{4}{2}{
   \picline{0 0}{4 2}
   \picmultiline{-5.0 1 -1.0 0}{0 2}{4 0}
   \picputtext{2.7 1}{A}
   \picputtext{1.3 1}{A}
   \picputtext{2 1.6}{B}
   \picputtext{2 0.4}{B}
} \quad
\diag{8mm}{4}{2}{
   \pictranslate{2 1}{
       \picmultigraphics[S]{2}{1 -1}{
           \piccurve{-2 1}{-0.3 0.2}{0.3 0.2}{2 1}
       }
       {\piclinedash{0.05 0.05}{0.01}
        \picline{0 -0.4}{0 0.4}
       }
   }
   \picputtext{2.7 1}{A}
   \picputtext{1.3 1}{A}
} \quad
\diag{8mm}{4}{2}{
   \pictranslate{2 1}{
       \picmultigraphics[S]{2}{-1 1}{
           \piccurve{2 -1}{0.1 -0.5}{0.1 0.5}{2 1}
       }
       {\piclinedash{0.05 0.05}{0.01}
        \picline{0 -0.6 x}{0 0.6 x}
       }
   }
   \picputtext{2 1.6}{B}
   \picputtext{2 0.4}{B}
}
\]
\caption{\label{figsplit}The A- and B-corners of a
crossing, and its both splittings. The corner A (resp. B)
is the one passed by the overcrossing strand when rotated 
counterclockwise (resp. clockwise) towards the undercrossing 
strand. A type A (resp.\ B) splitting is obtained by connecting 
the A (resp.\ B) corners of the crossing. It is useful to
put a ``trace'' of each splitted crossing as an arc connecting
the loops at the splitted spot.}
\end{figure}

Let $D$ be an unoriented link diagram.
A \em{state} is a choice of splittings of type $A$ or 
$B$ for any single crossing (see figure \ref{figsplit}), 
Let the \em{$A$-state} of $D$ be the state where all
crossings are $A$-spliced; similarly define the $B$-state.

We call a diagram \em{$A$-(semi)adequate} if in the $A$-state
no crossing trace (one of the dotted lines in figure
\reference{figsplit}) connects a loop with itself.
Similarly we define $B$-(semi)adequate. A diagram is
\em{semiadequate} if it is $A$- or $B$-semiadequate, and
\em{adequate} if it is simultaneously $A$- and $B$-semiadequate.
A link is adequate/semiadequate if it has an adequate/semiadequate
diagram. (Here $A$-adequate and $B$-adequate is what is called
$+$adequate resp.\ $-$adequate in \cite{Thistle}.)

We also need a part of the semiadequacy formulas for the Jones
polynomial. As in \cite{ntriv}, for an non-negative integer $i$, we
write $V_i$ for $[V]_{\md V+i}$ and $\bar V_i=[V]_{\Md V-i}$. (These 
are the $i+1$-st and $i+1$-st last coefficient of $V$ in its
coefficient list.) For an A- (resp. B-) semiadequate diagram, 
in \cite{ntriv}, and independently in \cite{DasLin}, formulas were 
obtained for $V_{1,2}$ (resp. $\bar V_{1,2}$).

\parbox[b]{12.5cm}{
\parskip 5pt plus 3pt minus 2pt\relax
{}From the $A$-state $A(D)$ of a diagram $D$ we define two
graphs.  The first graph, we call it the \em{$A$-graph}
$G(A)=G(A(D))$ has vertices for each loop in $A(D)$, and
an edge between each pair of loops that are connected
by a trace of at least one crossing in $D$. If there are
at least two such traces, we call the edge \em{multiple}.  

Let $\bigtriangleup(D)=\bigtriangleup(A(D))$ be the number of
cycles of length 3 (\em{triangles}) in $G(A(D))$.  

The second graph, we call it the \em{intertwining graph}
$IG(A)=IG(A(D))$, has vertices for each multiple edge
in $A(D)$, and an edge in $IG(A)$ is drawn between each
pair of vertices in $IG(A)$, whose corresponding multiple
traces in $A(D)$ contain traces of the form shown on the right.
}
\enspace\parbox[t]{4.0cm}{
\vbox to 4.0cm{\vss
\[
\diag{10mm}{4}{4}{
  \picline{2 0}{2 4}
  {\piclinedash{0.1 0.1}{0.01}
   \picline{1 0.8}{2 0.8}
   \picline{3 1.6}{2 1.6}
   \picline{1 2.4}{2 2.4}
   \picline{3 3.2}{2 3.2}
  }
  \picfilledellipse{1.6 1 x}{0.3 1.2}{}
  \picfilledellipse{2.4 3 x}{0.3 1.2}{}
}
\]\vss}
}

\begin{theorem}(\cite{ntriv,DasLin})\label{T02}
Assume that $D$ is a connected $A$-adequate diagram. Then
$V_0V_1=\chi(G(A))-1$ and
\[
V_0V_2=\mybin{2-\chi(G(A))}{2}+\chi(IG(A))-\bigtriangleup(A(D))\,.
\]
\end{theorem}

Here $\chi$ is the Euler characteristic (number of vertices minus
number of edges); $V_0=\pm 1$ by \cite{LickThis}.
Analogous formulas hold for $B$-adequate diagrams.

\subsection{Braids and braid words\label{brs}}

The $n$-string \em{braid group} $B_n$ is considered generated by
the Artin \em{standard generators} $\sg_i$ for $i=1,\dots,n-1$. These
are subject to relations of the type $[\sg_i,\sg_j]=1$ for $|i-j|>1$,
which we call \em{commutativity relations} (the bracket denotes the
commutator) and $\sg_{i+1}\sg_i\sg_{i+1}= \sg_i\sg_{i+1}\sg_i$,
which we call \em{Yang-Baxter} (or shortly YB) \rm{relations}.

We will make one noteworthy modification of this notation.
In the following $\sg_3$ will stand for the usual Artin generator
for braids of $4$ or more strands (such braids are considered
explicitly only in sections \reference{S2.2} and \reference{SBurau}),
while on $3$-braids (considered in all other sections) it will denote
the ``band'' generator $\sg_2\sg_1\sg_2^{-1}=\sg_1^{-1}\sg_2\sg_1$.

It will be often convenient in braid words to write $\pm i$ for
$\sg_i^{\pm 1}$. For example, $[21(33-4)^232]$ is an alternative
writing for $\sg_2\sg_1(\sg_3^2\sg_4^{-1})^2\sg_3\sg_2$.

The following definition summarizes basic terminology of
braid words used throughout the paper.

\begin{defi}\label{dfext}
Choose a word 
\begin{eqn}\label{wd}
\bt=\prod_{i=1}^n\sg_{k_i}^{l_i}
\end{eqn}
with $k_i\ne k_{i+1}$ and $l_i\ne 0$. We understand such words
\em{in cyclic order}.

Call $\sg_{k_i}^{l_i}$ the \em{syllables}
of $\bt$. For such a syllable, let $k_i$ be called the \em{index}
of the syllable, and $l_i$ its \em{exponent} or \em{length}.
We call a syllable $\sg_{k_i}^{l_i}$ \em{non-trivial} if $|l_i|>1$
and \em{trivial} if $|l_i|=1$.

We say $n$ is the \em{syllable length} of $\bt$ in \eqref{wd}.
The \em{(word) length} $c(\bt)$ of $\bt$ is $\sum_{i=1}^n |l_i|$, and
$c_{\pm}(\bt)=\sum_{\pm l_i>0}\pm l_i$ are the \em{positive/negative
length} of $\bt$. A word is \em{positive} if $c_-(\bt)=0$, or
equivalently, if all $l_i>0$. The \em{exponent sum} of $\bt$ is
defined to be $[\bt]=\sum_{i=1}^n l_i$. The \em{index sum} of
$\bt$ is $\sum_{i=1}^n k_i\cdot |l_i|$ (i.e.
each letter, not just syllable in $\bt$, contributes to that sum).

For $3$-braid words, the $k_i$ interchange between $1$ and $2$.
Thus the vector of the $l_i$, considered up to cyclic permutations,
determines the conjugacy class unambiguously. We call it
the \em{Schreier vector}.

Let $\bt$ be a positive word.
We call $\tl\bt$ an \em{extension} of $\bt$ if $\tl\bt$ is obtained by
replacing some (possibly no) trivial syllables in $\bt$ by non-trivial
ones of the same index. Contrarily, we call $\bt$ a \em{syllable
reduction} of $\tl\bt$. We call $\bt$ \em{non-singular} if $[\bt]_k
>1$ for all $k=1,\dots,n-1$, where $[\bt]_k=\sum_{k_i=k}l_i$ is the
exponent sum of $\sg_k$ in $\bt$. If $[\bt]_k=1$, we say that the
syllable of index $k$ in $\bt$ is \em{isolated} or \em{reducible}.
\end{defi}

To avoid confusion, it seems useful to clarify \em{a priori}
the following use of symbols (even though we recall it at
appropriate places later).

A comma separated list of integers
will stand for a \em{sequence of syllable indices} of a braid word. 
The `at' sign `@', written after such an index means that the
corresponding syllable is trivial, while by an exclamation mark
`$!$' we indicate that the syllable is non-trivial. (As the
exponent for non-trivial syllables will be immaterial, it is
enough to distinguish only whether the syllable is trivial or not.)
If none of $!$ and @ is specified, we do not exclude explicitly
any of either types. 

A bracketed but non-comma separated list of integers will stand
for a braid word. An asterisk `$*$' put after a letter (number)
in such a word means that this letter may be repeated
(it need not be repeated, but it must not be omitted).
So a (possibly trivial) index-2 syllable can be written
as $22*$. The expression `$[23]$' should mean a letter which
is either `2' or `3', and `$[23]+$' means a possibly empty
sequence of letters `$2$' and `$3$'.

\subsection{Braid representations of links\label{brt}}

By a theorem of Alexander \cite{Alexander}, any link is the \em{closure}
$\hat\bt$ of a braid $\bt$. The \em{braid index} $b(L)$ of a link $L$
is the smallest number of strands of a braid $\be$ whose closure
$\hat\be$ is $L$. See \cite{Morton,WilFr,Murasugi}. Such $\bt$ are
also called \em{braid representations} of $L$. The closure operation
gives for a particular braid word $\bt$ a link diagram $D=\hat\bt$.
Then we have for example $w(\hat\bt)=[\bt]$, $c(\bt)=c(\hat\bt)$,
$c_\pm(\bt)=c_\pm(\hat\bt)$, and $s(\hat\bt)$
is the number of strings of $\bt$ (i.e. $n$ for $\bt\in B_n$).

Many properties of braid words we will deal with
relate to the corresponding properties of their link diagrams.
For example, a braid word is called \em{positive} if it contains no
$\sg_i^{-1}$, or in other words, its closure diagram is positive.
A braid is positive if it has a positive braid word.
In a similar fashion, we say that braid is ($A$/$B$-)adequate if
it has an ($A$/$B$-)adequate word representation, and a word is
adequate if the link diagram obtained by its closure is adequate. 

If a braid word $\bt$ is written as $\sg_1^{\pm 1}\ap\sg_1^p\ap'$,
where $p\in\bZ$ and none of $\ap$ and $\ap'$ contains a syllable of
index $1$, then the diagram admits a \em{flype}, which exchanges
the syllables of index $1$ in $\bt$, so that we obtain $\sg_1^p
\ap\sg_1^{\pm 1}\ap'$. This operation preserves the isotopy
type of the closure link $\hat\bt$, but in general changes
the braid conjugacy class. The phenomenon is explained in
\cite{BirMen}. In the context of general link diagrams,
the flype has been studied also extensively, most prominently
in \cite{MenThis}.

Alternatively to the standard Artin generators, one
considers also a representation of the braid groups by
means of an extended set of generators (and their inverses)
\[
\sg_{i,j}^{\pm 1}=\sg_i\dots\sg_{j-2}\sg_{j-1}^{\pm 1}
\sg_{j-2}^{-1}\dots\sg_i^{-1}
\]
for $1\le i<j\le n$. Note that 
\begin{eqn}\label{op}
\sg_i=\sg_{i,i+1}\,.
\end{eqn}
A representation of a braid $\bt$, and its closure link $L=\hat\bt$,
as word in $\sg_{i,j}^{\pm 1}$ is called a \em{band representation}
\cite{BKL}. A band representation of $\bt$ spans naturally a
Seifert surface of the link $L$ as in figure \ref{figbd}: one
glues disks into the strands, and connects them by half-twisted
bands along the $\sg_{i,j}$. The resulting surface is called
\em{braided Seifert surface} of $L$. A minimal genus Seifert
surface of $L$, which is a braided Seifert surface, is also
called a \em{Bennequin surface} \cite{BirMen2}.

In this paper we will deal exclusively with band representations
in $B_3$. Then we have three band generators $\sg_{i,i\bmod 3+1}$
(where $i=1,2,3$, and `$\bmod$' is taken with values between $0$ and
$2$). With \eqref{op}, we have $\sg_1=\sg_{1,2}$ and $\sg_2=\sg_{2,3}$,
and with the special meaning of $\sg_3\in B_3$ introduced above,
$\sg_3=\bar\sg_{1,3}$, where bar denotes the mirror image.
(This mirroring is used here for technical reasons related to
Xu's normal form, as explained below.)

If a band representation contains only positively half-twisted bands
(i.e. no $\sg_{i,j}^{-1}$ occur), it is called \em{band-positive}
or \em{strongly quasi-positive}. A link with a strongly quasi-positive
band (braid) representation is called strongly quasi-positive.
Such links have an importance in connection to algebraic curves;
see \cite{Rudolph2}.

Using the skein polynomial, define a quantity by
\begin{eqn}\label{mwfb}
MWF(L)\,=\,\frac 12\,\big(\Md_vP-\md_vP\big)+1\,.
\end{eqn}
The \em{Morton-Williams-Franks braid index inequality}
\cite{Morton,WilFr} (abbreviated as MWF) states that
\begin{eqn}
b(L)\ge MWF(L)
\end{eqn}
for every link $L$. This inequality is often exact (i.e. an
equality). The study of links where it is exact or not
has occupied a significant part of previous literature.
Most noteworthy is the work of Murasugi \cite{Murasugi}
and Murasugi-Przytycki \cite{MurPrz2}.

The Morton-Williams-Franks 
results from two other inequalities, due to Morton,
namely that for a diagram $D$, we have
\begin{eqn}\label{mq}
1-s(D)+w(D)\,\le\,\md_lP(D)\,\le\,\Md_lP(D)\,\le\,s(D)-1+w(D)\,.
\end{eqn}
Williams-Franks showed these inequalities for the case of braid
representations (i.e. when $D=\hat\bt$ for some braid $\bt$).
Later it was observed from the
algorithm of Yamada \cite{Yamada} and Vogel \cite{Vogel}
that the braid version is actually equivalent to (and not just
a special case of) the diagram version.

These inequalities were later improved in \cite{MurPrz2} in a way 
that allows to settle the braid index problem for many links
(see \S\ref{Sfg} or also \cite{Ohyama}). 

In the special case of $3$-braids (and with the special meaning
of $\sg_3$ as described above), Xu \cite{Xu} gives a normal form
of a conjugacy class in $\sg_{1,2,3}$. By Xu's algorithm, each
$\bt\in B_3$ can be written in one of the two forms
{\def\theenumi{\Alph{enumi}}
\def\labelenumi{(\theenumi)}
\begin{enumerate}
\item\label{itA} $[21]^kR$ or $L^{-1}[21]^{-k}$ ($k\ge 0$), or
\item\label{itB} $L^{-1}R$,
\end{enumerate}
}where $L$ and $R$ are positive words in $\sg_{1,2,3}$ with (cyclically)
non-decreasing indices (i.e. each $\sg_i$ is followed by $\sg_i$
or $\sg_{i\bmod 3+1}$, with `$\bmod$' taken with values between $0$
and $2$). Since the form \reference{itB} must be cyclically
reduced, we may assume that $L$ and $R$ do not start or end with
the same letter. This form is the shortest word in $\sg_{1,2,3}$
of a conjugacy class. By Bennequin's aforementioned result, the
braided surface is then a minimal genus (or Bennequin) surface.

\subsection{Gau\ss{} sum invariants}

We recall briefly the definition of Gau\ss{} sum invariants.
They were introduced first in \cite{Fied} for braids, and later
\cite{Fied2,VirPol} for knots. It is known that all
they give formulas for Vassiliev invariants.

\begin{defi}(\protect\cite{Fied2})\ 
A Gau\ss{} diagram of a knot diagram is an oriented circle with
arrows connecting points on it mapped to a crossing and
oriented from the preimage of the undercrossing (underpass)
to the preimage of the overcrossing (overpass).
\end{defi}

We will call a pair of crossings whose arrows intersect in
the Gau\ss{} diagram a \em{linked} pair.

\begin{figure}[htb]
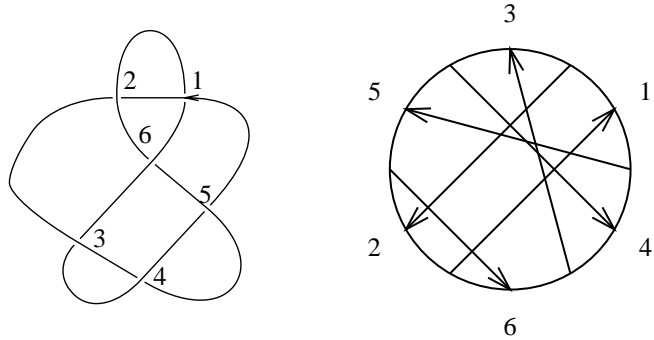

{%\small
\[
\begin{array}{c@{\qquad}c}
     \diag{2cm}{2.0}{2.0}{
  \picputtext[dl]{0 0}{\autoepsffalse\epsfs{4cm}{6_2}}
  \picputtext[u]{1.25 1.65}{$1$} 
  \picputtext[u]{0.8 1.65}{$2$} 
  \picputtext[u]{0.6 0.6}{$3$} 
  \picputtext[u]{1.0 0.35}{$4$} 
  \picputtext[u]{1.3 0.9}{$5$} 
  \picputtext[u]{0.9 1.3}{$6$} 
	}
 & 
\diag{1.6cm}{2.4}{2.4}{\pictranslate{1.2 1.2}{
\piclinewidth{60}
\piccircle{0 0}{1}{}
\labar{240}{ 30}{1}
\labar{ 60}{210}{2}
\labar{300}{ 90}{3}
\labar{120}{330}{4}
\labar{  0}{150}{5}
\labar{180}{270}{6}
} % \pictranslate
} % \diag
\end{array}
\]}
\caption{The standard diagram of the knot $6_2$ and its Gau\ss{}
diagram.\label{fig6_2}}
\end{figure}

\begin{exam}
As an example, figure \ref{fig6_2} shows the knot $6_2$ in its
commonly known projection and the corresponding Gau\ss{} diagram.
\end{exam}

The simplest (non-trivial) Vassiliev knot invariant is the
\em{Casson invariant} $v_2$, with $v_2=\Dl''(1)/2=-V''(1)/6$,
for which Polyak-Viro
\cite{VirPol,VirPol2} gave the simple Gau\ss{} sum formula
\begin{eqn}\label{pv2_}
v_2\,=\,\GD{\arrow{225}{45}
\arrow{315}{135}
\picfillgraycol{0}\picfilledcircle{1 90 polar}{0.09}{}
}\,.
\end{eqn}
Here the point on the circle corresponds to a point on the
knot diagram, to be placed arbitrarily except on a crossing.
(The expression does not alter with the position of the
basepoint; we will hence have, and need, the freedom to place
it conveniently.)

We will use the symmetrized version of \eqref{pv2_}
w.r.t. taking the mirror image of the knot diagram:
\begin{eqn}\label{pv2}
v_2\,=\,\frac{1}{2}\es\left(
\GD{\arrow{225}{45}
\arrow{315}{135}
\picfillgraycol{0}\picfilledcircle{1 90 polar}{0.09}{}
}+
\GD{\arrow{45}{225}
\arrow{135}{315}
\picfillgraycol{0}\picfilledcircle{1 90 polar}{0.09}{}
}\right)\,.
\end{eqn}

\section{Xu's form and Seifert surfaces\label{SA}}

We studied the relation of Xu's algorithm and the skein
polynomial in \cite{3br}, and here we will go further to
connect fiberedness, the and Alexander and Jones polynomial
to Xu's form.

\subsection{Strongly quasi-positive links among links of braid index 3}

Theorem \reference{thqp} follows relatively easily from
the work in \cite{3br}, but it is a good starting point for
the later more substantial arguments basing on Xu's form.

\proof[of theorem \reference{thqp}]
The inequalities of \eqref{mq} for the $v$-degree
of $P$ applied on a positive band representation show that,
with $P=P(L)$ and $\chi=\chi(L)$,
\begin{eqn}\label{x}
\md_vP\ge 1-\chi\,.
\end{eqn}
Because of \cite[proposition 21]{LickMil}, which now says
\begin{eqn}\label{p21}
P(v,v^{-1}-v)=1\,,
\end{eqn}
we have $\md_lP\le
\Md_zP$, and by \cite{3br} $\Md_zP=1-\chi$. So from \eqref{x}
we obtain 
\begin{eqn}\label{y}
\md_lP=1-\chi\,. 
\end{eqn}

It is known that the minimal degree term in $z$ of the skein polynomial
of a $n$-component link is divisible by $(v-v^{-1})^{n-1}$. (In
\cite{Kneissler} in fact all occurring terms are classified.) So
if $MWF (K)=1$, then $K$ is a knot. Now, the identity \eqref{p21}
implies that if $MWF=1$ for some knot $K$, then $P(K)=1$. For
any non-trivial knot $1-\chi=2g>0$, so for any non-trivial strongly
quasi-positive knot from \eqref{x} we have $\md_vP>0$, and so $P\ne 1$.

So a strongly quasi-positive link cannot have $P=1$, and always
$MWF\ge 2$.  % (The split links are easy to deal with.)
Therefore, when the braid index is 3, we have
$MWF\in \{2, 3\}$. Then \eqref{y} and the inequalities \eqref{mq}
for the $v$-degree of $P$ show that a 3-braid representation $\bt$
has exponent sum $[\bt]=3-\chi$, unless $MWF=2$ and $[\bt]=1-\chi$.

In former case, we can find a minimal genus band representation
from $\bt$ by Xu's algorithm \cite{Xu}, and this representation
must be positive. In latter case, we will have one negative band,
and have Xu's form $L^{-1}R$. Here $L$ and $R$ are positive words in
the letters $\sg_1,\sg_2,\sg_3=\sg_2\sg_1\sg_2^{-1}$ with $\sg_i$
followed by $\sg_i$ or $\sg_{i\bmod 3+1}$ (as described
in \S\reference{brt}), and $L$ is a single letter.

W.l.o.g. assume $L=\sg_3$. If the first or last letter in $R$
is a $\sg_3$ then we can cancel two bands, and have a positive
representation. If some $\sg_3$ occurs in $R$, then we can write
$L^{-1}R=\sg_3^{-1}\ap\sg_3\ap'$ where $\ap,\ap'$ are positive
band words. Such a representation is quasi-positive with a Seifert
ribbon \cite{Rudolph} of smaller genus. However, it is known that
for strongly quasi-positive links the genus and 4-genus coincide
(see for example \cite{posqp}). If $R$ has no $\sg_3$, then
we have up to cyclic letter permutation a braid of the form
$\sg_1^k\sg_2^l\sg_1^{-1}\sg_2^{-1}$ with $k\ge 0$, $l>0$. This
braid is easily seen to reduce to (a positive one on) two strands.
\qed

\begin{rem}\label{rqw}
It is not true that all 3-braids that have the skein polynomial of
$(2,n)$-torus links have positive 3-band representations. Birman's
\cite{Birman} construction (see definition \reference{DEFB} below)
yields examples like $12_{2037}$ in \cite{KnotScape} with one
negative band. This construction gives in fact all non-obvious
3-braid $MWF=2$ examples, a circumstance shown in \cite{ntriv} by
applying theorem \ref{thq} (and reportedly in previous unpublished
work of El-Rifai). The result is explained and used below in the
proof of theorem \reference{thJC}.
\end{rem}

\subsection{Uniqueness of minimal genus Seifert surfaces\label{Uq}}

The discussion here came about from the desire to complete Xu's
uniqueness theorem for Seifert surfaces of 3-braid links. While
our result will be further improved later using work to follow,
we need to introduce some notation and basic tools. We will use
the work of Kobayashi \cite{Kobayashi}, implying that the property
a surface to be a unique minimal genus surface is invariant under
Hopf (de)plumbing.

Subsequently Mikami Hirasawa advised me about a subtlety
concerning the notion of `uniqueness' which must be explained.
Xu's work considers Seifert surfaces unique in the sense isotopic
to each other, if we \em{may move} the link by the isotopy.
In particular, such isotopy may interchange link components.
However, we demand the isotopy to preserve component orientation.
(For unoriented isotopy Xu's result is complete.) Contrarily,
Kobayashi's setting assumes uniqueness in the sense that
Seifert surfaces are isotopic to each other \em{fixing} the link.
Hirasawa explained that the definitions are not equivalent, and
that a unique (if we may move the link) minimal genus surface
may get not unique when one (de)plumbs a single Hopf band. To
account for this discrepancy, we should establish proper language.

\begin{defi}
Let us use the term \em{unique} for
Kobayashi's notion of uniqueness (up to isotopy fixing the
link), and let us call Xu's notion of uniqueness (up to
isotopy which may move the link or permute components,
but preserves orientation) \em{weakly unique}.
\end{defi}

Here we state the following extension of Xu's uniqueness theorem.

\begin{theorem}\label{wu}
Every 3-braid link has a weakly unique minimal genus Seifert surface.
\end{theorem}

\proof Birman-Menasco showed that any minimal genus Seifert surface
of a 3-braid link is isotopic to a Bennequin surface. Then, Xu
showed that every conjugacy class of 3-braids carries a canonical
(up to oriented isotopy) Bennequin surface. The conjugacy classes
of 3-braids with given closure link were classified in \cite{BirMen}.
Most links admit a single conjugacy class, and we are done, as in
\cite{Xu}. The exceptional cases are easy to deal with, except the
``flype admitting'' braids 
\begin{eqn}\label{fl3}
\sg_1^{\pm p}\sg_2^{\pm q}\sg_1^{\pm r}\sg_2^{\pm 1}\,.
\end{eqn}
The flype interchanges $\pm q$ and $\pm 1$, and (in general) gives
a different conjugacy class, which differs from the original one
by orientation. So we have (as in \cite{Xu}), two Bennequin surfaces
isotopic only up to orientation. We will settle this case now also
for oriented isotopy, that is, show that these surfaces are
isotopic to themselves with the opposite orientation.

Recall that for a fibered link, a minimal genus surface is the
same as a fiber surface, and such a surface is unique (Neuwirth-%
Stallings theorem). Moreover, by work of Gabai \cite{Gabai,Gabai2,%
Gabai3} and Kobayashi \cite{Kobayashi} the properties of a surface
to ba a minimal genus surface, a unique minimal genus surface or a
fiber surface are invariant under Murasugi (de)sum with a fiber
surface. In particular, this invariance holds for (de)plumbing
a Hopf band (which is understood to be an unknotted annulus with 
one full, positive or negative,
twist). Now we note that in the cases in \eqref{fl3} where $p=r$,
the flype is trivial (i.e. realized by a conjugacy), so that the
surface is weakly unique. However, since Kobayashi's theorem may
fail for weakly unique surfaces, we cannot reduce our surfaces to
this case. We will use Hopf (de)plumbings to recur all cases to fiber
surfaces or a 2-full twisted annulus. Then we understand that our
surfaces are unique, and in particular weakly unique. The type of
Hopf (de)plumbings we will apply is to interconvert all powers
of a given band generator of given sign (for example $\sg_3^k$ is
equivalent to $\sg_3^{-1}$ for each $k<0$).

W.l.o.g. assume in \eqref{fl3} that we have $-$ in $\pm 1$ and that none
of $p,q,r$ is $0$ (the other cases are easy). We assume $p,q,r>0$
and vary the signs before $p,q,r$ properly. Also, since the flype
interchanges $\pm p$ and $\pm r$, we may assume $\pm p\ge \pm r$.
So we exclude the sign choice $(-p,+r)$.

\begin{caselist}
\case $p,q,r>0$. We can write \eqref{fl3} up to cyclic permutation
as $(\sg_2^{-1}\sg_1\sg_2)^p\sg_2^{q-1}\sg_1^r$, and conjugating with
$\sg_2\sg_1\sg_2$ we have $\sg_3^p\sg_1^{q-1}\sg_2^r$. Since $p,r>0$,
this is Hopf plumbing equivalent to $\sg_3\sg_1\sg_2$ or $\sg_3\sg_2$.
These cases are a disk and a 2-full twisted annulus for
the reverse $(2,4)$-torus link, and we are done.\label{cs1}

\case $p,r>0$, $-q<0$. We have an alternating braid, which is
fibered.

\case $p,q,r<0$. We have a negative braid, which is fibered.

\case $-p,-r<0$, $q>0$. We have a word of the form 
$\sg_2^{q}\sg_1^{-r}\sg_2^{-1}\sg_1^{-p}$. If one of
$p$ and $r$ is $1$, then we have a $(2,\,.\,)$-torus link.
If $q=1$, then we can go over to the mirror image, and land in
case \reference{cs1}. Now, when $p,q,r\ge 2$, the
minimal word in Xu's form is $2^{q-2}(-3)(-2)^{r-1}(-1)^{p-1}$.
Since $p-1$ and $r-1$ are non-zero, its surface is plumbing
equivalent to the one of $[2-3-2-1]$ (if $q>2$) or $[-3-2-1]$
(if $q=2$). These are the annulus for the reverse $(2,-4)$-torus link,
and the fiber of the $(-2,-2,2)$-pretzel link, and we are done. 

\case $p>0$, $-r<0$, $q>0$.
We have 
\[
\sg_1^p\sg_2^q\sg_1^{-r}\sg_2^{-1}=\sg_1^p\sg_2^{q-1}
(\sg_2\sg_1^{-1}\sg_2^{-1})^r=\sg_1^p\sg_2^{q-1}\sg_3^{-r}\,.
\]
Such a surface is plumbing equivalent to the one for $\sg_1\sg_3^{-1}$
(if $q=1$) resp. $\sg_1\sg_2\sg_3^{-1}$ (if $q>1$), which are the
fibers of the unknot and Hopf link resp.

\case $p>0$, $-r<0$, $-q<0$. We have 
\[
\sg_1^{-r}\sg_2^{-1}\sg_1^p\sg_2^{-q}=\sg_1^{1-r}(\sg_1^{-1}\sg_2^
{-1}\sg_1)\sg_1^{p-1}\sg_2^{-q}=\sg_1^{1-r}\sg_3^{-1}\sg_1^{p-1}
\sg_2^{-q}\,.
\]
If $-r=-1$ then we have a $(2,p-q-1)$-torus link. If $-r<-1$, then
the band surface from the right word above is plumbing equivalent to the
one for $\sg_1^{-1}\sg_3^{-1}\sg_2^{-1}$ (if $p=1$) or 
$\sg_1^{-1}\sg_3^{-1}\sg_1\sg_2^{-1}$ (if $p>1$), which are the
annulus for the reverse $(2,4)$-torus link, and the fiber of the
$(2,2,-2)$-pretzel link. \qed
\end{caselist}

\subsection{Fiberedness\label{ppp}}

For the rest of the paper we normalize
$\Dl$ so that $\Dl(1)=1$ and $\Dl(t)=\Dl(1/t)$.

\begin{theorem}\label{Thq}
Let $L$ be a strongly quasi-positive 3-braid link. Then the following
are equivalent:
\def\labelenumi{\theenumi)}
\begin{enumerate}
\item $L$'s minimal genus surface is a Hopf plumbing,
\item $L$'s minimal genus surface is a fiber surface,
\item $\Md\Dl(L)=1-\chi(L)$ and $\Mc\Dl(L)=\pm 1$,
\item $L$'s Xu normal form is \em{not} $R$, with syllable
length of $R$ divisible by $3$. In other words, the Xu normal
form is not an extension (in the sense of definition
\ref{dfext}) of $[(123)^k]$ for $k>0$.
\item Some minimal band form of $L$ contains
$\sg_1^k\sg_2^l\sg_1^m$ or $\sg_2^k\sg_1^l\sg_2^m$
as subword for $k,l,m>0$, or is $\sg_1^k\sg_2^l$.
\end{enumerate}
\end{theorem}

\proof  $(1)\So (2)\So (3)$ are clear.

$(5)\So (1)$. Assume after deplumbing, all letters occur in single
power, and up to conjugacy the word starts with $\sg_1\sg_2\sg_1$
or $\sg_2\sg_1\sg_2$. By adjusting one of the two, we can have a
$\sg_1^2$ or $\sg_2^2$ if the next letter is $\sg_1$ or $\sg_2$. In
that case we deplumb a Hopf band. If the next letter is $\sg_3$, then
we have $\sg_1\sg_2\sg_1\sg_3=\sg_1\sg_2^2\sg_1$, can can also deplumb
a Hopf band. Then we reduce the surface for that of $\sg_1\sg_2\sg_1$
which is the Hopf band.

$(3)\So (4)$. We prove the contrary. Assume $(4)$ does not hold. 
Band-positive surfaces are always of minimal genus, so that the
properties we investigate are invariant under Hopf (de)plumbings.
Under applying skein relations at non-trivial syllables we are
left with powers of $\sg_1\sg_2\sg_3$. Apply the skein relation
for $\Dl$ at the last band. Then $\bt_{-}$ and $\bt_0$ are both
of minimal length. We already proved, in $(5)\so (1)$, that
$\bt_0$ is fibered, so $\Md\Dl=1-\chi(\bt_0)$ and $\Mc\Dl=\pm 1$.
The same holds for $\bt_{-}$ by \cite[proposition 2]{3br}, since 
$\bt_{-}$ is of Xu's minimal form $L^{-1}R$ and is not positive.
So the terms in degree $1-\chi(\bt)$ of $\Dl(\bt)$ either cancel, or
give $\pm 2$. 

$(4)\So (5)$. We prove the contrary. Assume $(5)$ does not hold.
If we do not have a word $\bt=R$ in Xu's form of length divisible by
$3$, another option would be to have a word with cyclically decreasing
indices. But note that $[321]=[121]$ is the Hopf band, while
$[321321]=[211211]$ contains a $121$, too, and so we
are done. Otherwise the index array of
the syllables of $\bt$ must contain the same entry with distance $2$,
and such a word is conjugate to the ones excluded.
\qed

\begin{theorem}\label{xLR}
Any closed braid of Xu's form $L^{-1}R$ is fibered.
\end{theorem}

The proof uses some work of Hirasawa-Murasugi.
A consequence of their result is the following lemma,
which we require. It is proved in appendix \reference{HM}.

\begin{lemma}\label{lki}
The links $L_k$, given by the closed 3-braids $[(123)^k\,-2]$ for $k>0$,
are fibered.
\end{lemma}

% \proof Consider the knots $L_k'$ given by the closed 3-braid
% words $[(12)^{3k+1}(-2)^{3k+l}]$, where $l=2$ for
% $k$ even, and $l=3$ for $k$ odd. The Xu normal form is found using
% the identity $2^3(123)=(12)^3=(123)2^3$ to be
% \[
% (12)^{3k+1}(-2)^{3k+l}\to 12(123)^k(-2)^l\to 21223(123)^{k-1}(-2)^l\to
% 1223(123)^{k-1}(-2)^{l-1}\,,
% \]
% and this surface is Hopf plumbing equivalent to that of $L_k$.
% 
% We know from \cite{3br} (as explained further in \cite{ntriv})
% that if $\bt\in B_3$ is not strongly quasi-positive, then
% $2\Md\Dl=1-\chi$ and $\Mc\Dl=\pm 1$.
% 
% Now, since $L_k'$ has Xu normal form of the type $L^{-1}R$, we
% find that $\Mc\Dl(L_k')=\pm 1$. Then, by the result of
% Hirasawa-Murasugi $L_k'$ is fibered, and hence so is $L_k$. \qed

\proof[of theorem \reference{xLR}]
We use induction on the length of $L^{-1}R$ and for fixed
length on the exponent sum. Under Hopf deplumbings assume all
syllables in $L$ and $R$ are trivial.
Assume up to mirroring that $L$ is not shorter than $R$.
Permute by conjugacy $R$ to
the left, and permute the indices so that $L$ starts with $-3$. 

The following transformations also offer a Hopf deplumbing
\[
1-3-2-1\to -2 1 -2 -1\to -2 -2 -1 2\to -2 -1 2
\]
and
\[
1 2 -3\to -2 1 1
\]
These reductions fail both if either $L$ and $R$ have length at most
$2$, or $R$ has length $1$. (Remember $L$ is not shorter than $R$.)
In former case one checks directly that one has a disk, Hopf band or
connected sum of two Hopf bands. We consider latter case.

By conjugacy permute the indices so that $R=1$; also assume $L$ has
length at least $3$. If $L$ starts and ends with $-2$ then we do
\[
1 -2 -1 \dots -2 \to -2 3 -1 \dots -2\to 3 -1 \dots -2 -2
\]
(where the right transformation is a conjugacy) and deplumb a
Hopf band. If $L$ starts with $-3$ then we transform as before
\[
1 -3 -2 -1\to -2 1 -2 -1\to -2 -2 -1 -2
\]
and deplumb a Hopf band. So $L$ starts with $-2$ and ends on $-3$.
Then the mirror image of $\bt$ is up to conjugacy of the form
$[(123)^k\,-2]=[(1 2 2 1)^k (-2)^{k+1}]$, which we dealt with in
the lemma before.  \qed

Combinedly, we obtain

\begin{corr}
Let $W$ be a 3-braid link. Then the following are equivalent:
\def\labelenumi{\theenumi)}
\begin{enumerate}
\item $W$ is fibered,
\item $\Md\Dl(W)=1-\chi(W)$ and $\Mc\Dl(W)=\pm 1$,
\item $W$'s Xu form is not of the type $L^{-1}$ or $R$,
with syllable length of $L$ or $R$ divisible by $3$. \qed
\end{enumerate}
\end{corr}

\begin{exam}
The routine verification, with Mikami Hirasawa, of the tables in
\cite{KnotScape} for fibered knots, has shown non-fibered knots 
with monic Alexander polynomial of degree matching the genus
start at 12 crossings. One of these knots, $12_{1752}$, has
braid index $4$, so that (as expected) the corollary does not
hold for 4-braids.
\end{exam}

\begin{corr}
Any non-split closed 3-braid $\hat\bt$ with $|[\bt]|\le 2$ is fibered,
in particular so is any slice 3-braid knot. \qed
\end{corr}

Again $6_1$ and $9_{46}$ show the second part is not true for 4-braids.
The first statement in the corollary extends simultaneously
the property of amphicheiral knots, which follows also
from \cite{BirMen}, since amphicheiral 3-braid
knots are closed alternating 3-braids.

Originally the insight about a unique minimal genus surface motivated
the fibered 3-braid link classification. Still this insight lacks
asset as to the somewhat improper way it emerges. Note for example,
that the stronger version of \cite{BirMen2} of Bennequin's theorem
enters decisively into \cite{Xu} and the proof of theorem \ref{wu},
but then paradoxically latter imply that Birman-Menasco's formulation
is actually equivalent to, and not really an improvement of,
Bennequin's theorem. Also, Birman-Menasco's classification
\cite{BirMen} priorizes Schreier's conjugacy algorithm and lacks
any geometric interpretation of Xu's normal form, while such an
interpretation becomes evident in our setting. This provided strong
motivation for theorem \reference{thuq}. Its proof is completed
by dealing with the non-fibered cases in appendix \reference{HI}.
It requires also a part of the further detailed consideration of the
Alexander polynomial that follows next.

\section{Polynomial invariants}

\subsection{Alexander polynomial}

In \cite{Birman}, Birman proposed (but considered as very difficult)
the problem to classify 3-braid
links with given polynomials. In \cite{3br} we dealt with the skein
polynomial. Now we can extend our results to the Alexander polynomial
(with the convention in the beginning of \S\reference{ppp}).
The following discussion gives a fairly exact description how to
find the 3-braid links, if such exist, for any possible admissible
(as specified in \S\reference{brD}) polynomial.

A solution for the Jones polynomial is presented afterwards.
% (though not explained in detail)

\begin{lemma}
If $\bt$ is strongly quasi-positive and fibered, then $\Mc\Dl=+1$.
\end{lemma}

(Here it is essential to work with the leading, not trailing
coefficient of $\Dl$ and with strongly quasi-positive links and
not their mirror images.)

\proof We know $\md_vP=\Md_zP=1-\chi$. So $[P]_{z^{1-\chi}}$ has
a term in degree $v^{1-\chi}$. The coefficient must be $+1$
because of \eqref{p21} and because it is the only coefficient that
contributes to the absolute term in $P(v,v^{-1}-v)$. Now from
the classification of leading $z$-terms of $P$ in \cite{3br}
it follows that $[P]_{z^{1-\chi}}$ can have at most one further term,
with coefficient $\pm 1$. If such term exists, the substitution
\eqref{DP} would either cancel the terms
in degree $t^{(1-\chi)/2}$ in $\Dl(t)$, or give coefficient $\pm 2$,
so our link cannot be fibered. Thus a second term does not exist,
and the claim follows from \eqref{DP}. \qed

\begin{lemma}\label{ty}
Let $\bt'=[(123)^k]$ be an even power $k$ of $[123]$.
Assume $\bt$ has Xu normal form $R$, and after syllable
reduction becomes $\bt'$, but $\bt\ne\bt'$ itself (i.e. some
syllable in $\bt$ is non-trivial). Then $\hat\bt$ satisfies
$\Md\Dl=-1-\chi$; moreover $\Mc\Dl>0$ and is equal to the 
number of non-trivial syllables in $\bt$.
\end{lemma}

\proof We use for fixed $k$ induction on the exponent sum.
If exactly one syllable is non-trivial with exponent $2$,
then applying the skein relation at the exponent-2
syllable shows that $\bt=\bt_+$ inherits the Alexander
polynomial of $\bt_-$ with positive sign (since $\bt_0=\bt'$,
whose closure has zero polynomial). 
Now $\bt_-$ is positive and has an index array with a subsequence
of the form $xyx$, and so is fibered by theorem \reference{Thq}.
Then by lemma \reference{ty}, $\Mc\Dl=+1$ and
$\Md\Dl=1-\chi(\hat\bt_-)=-1-\chi(\hat\bt)$.

If $[\bt]-[(123)^k]>1$, then applying the skein relation
at any non-trivial exponent $2$ syllable gives (with positive
sign) the Alexander polynomials of two closed braids
$\bt_0$ and $\bt_-$, former of which is fibered and
latter of which has the requested property by induction.
Then the maximal terms in degree $1-\chi(\hat\bt_-)=
-\chi(\hat\bt_0)$ are positive and do not cancel.
The one of $\bt_-$ is $+1$, while the one of $\bt_0$
by induction one less than the number of non-trivial 
syllables in $\bt$ (since in $\bt_0$ one more syllable
becomes trivial).

If some syllable in $\bt$ has exponent $>2$ then $\bt_-$ is
not minimal, and the degree and leading coefficient of $\Dl$
are inherited (with positive sign) from $\bt_0$. \qed

\begin{lemma}\label{tz}
Assume $\bt$ has Xu normal form $R$, which after syllable
reduction becomes an odd power $k$ of $[123]$.
Then $\Md\Dl=1-\chi$, and $\Mc\Dl=+2$.
\end{lemma}

\proof First we prove the claim if $\bt=[(123)^k]$.
We know that $\Md\Dl=1-\chi$ and $\Mc\Dl=\pm 2$,
so we must exclude $\Mc\Dl=-2$. Applying the skein relation gives
the polynomials of $\bt_0$ and $\bt_-$ with positive sign.
$\bt_0$ is positive and fibered as before, so $\Mc\Dl=+1$.
Then clearly $\bt_-$ (which is of the form $L^{-1}R$ and also
fibered) cannot have $\Mc\Dl=-3$. 

If $\bt\ne [(123)^k]$, then it has a non-trivial syllable.
Applying the skein relation at a letter in that syllable
we find that $\bt_-$ reduces. So the leading term comes
from $\bt_0$, and with positive sign. \qed

\begin{theorem}\label{THA}
Fix some admissible Alexander polynomial $\Dl$. Then
\def\labelenumi{\theenumi)}
\begin{enumerate}
\item If $\Dl=0$, then the 3-braid links with such polynomial
are the split links and the closures of (incl.\ negative)
even powers of $[123]$.
\item If $\Dl\ne 0$, there are only finitely many 3-braid links
with this $\Dl$. They all have $1-\chi=2\Md\Dl$ or $-1-\chi=2\Md\Dl$.
In latter case they are up to mirroring strongly quasi-positive.
\item If $\Mc\Dl\le -2$, then no 3-braid knot or 3-component
link has such Alexander polynomial, and any 2-component link
is strongly quasi-negative.
\item If $\Mc\Dl\ge +2$, any 3-braid knot or 3-component link
with such Alexander polynomial is strongly quasi-positive or a
mirror image of it. Any 2-component link is strongly quasi-positive.
\item If $|\Mc\Dl|>2$, then any 3-braid link 
with such Alexander polynomial has $2\Md\Dl=-1-\chi$.
\end{enumerate}
\end{theorem}

\begin{defi}\label{DEFB}
Let for $\bt\in B_3$ with $6\mid [\bt]$, the \em{Birman dual}
$\bt^*$ be defined by $\bt^{-1}\dl^{2[\bt]/3}$, where
$\dl=[121]$ and $\dl^2$ generates the center of $B_3$.
\end{defi}

Birman \cite{Birman} shows that $\hat\bt$ and $\hat\bt^*$
have the same skein polynomial. This observation relates to
our explanation at a couple of places, for example, in remark
\reference{rqw}, and also in the below arguments.

\proof[of theorem \reference{THA}]
Let us exclude \em{a priori} trivial and split links.  The claims
follow from the discussion of $2\Md\Dl$ and $\Mc\Dl$ in cases.

We know by theorem \reference{xLR} (or by \cite{3br,ntriv}, as noted
before) that if $\bt\in B_3$ is not (up to mirroring) strongly
quasi-positive, then $2\Md\Dl=1-\chi$ and $\Mc\Dl=\pm 1$.

It remains to deal with the Xu form $R$. The form $L$ is just the
mirror image, and mirroring preserves the Alexander polynomial for
knots and 3-component links and alters the sign for 2-component
links.

If making trivial all syllables in $\bt$, the new word
$\bt'$ is not a power of $[123]$, then we proved that $\hat\bt$
is fibered, so again $2\Md\Dl=1-\chi$ and $\Mc\Dl=\pm 1$.

If $\bt'$ is an even power of $[123]$, then by Birman duality we
conclude that $\Dl(\hat\bt')=0$. If $\bt\ne \bt'$, then one uses
lemma \reference{ty}.
% induction on the syllables exponents. If $\bt$ and $\bt'$
% differ only by a letter, then applying the skein relation at
% the exponent-2 syllable shows that $\bt=\bt_+$ inherits the Alexander
% polynomial of $\bt_-$ (since $\bt_0=\bt'$). This has form $L^{-1}R$
% and so $2\Md\Dl=1-\chi(\hat\bt_-)=-1-\chi(\hat\bt)$ and $\Mc\Dl=\pm 1$.
% Then the degree and leading term are successively propagated.
Thus $2\Md\Dl=-1-\chi$.
% and also here $\Mc\Dl=\pm 1$.

If $\bt'$ is an odd power of $[123]$, then use the observation in
lemma \reference{tz}  (or make one syllable to exponent $4$ and
use Birman duality) to conclude that $2\Md\Dl=1-\chi$ and
$\Mc\Dl=+2$. \qed

\begin{exam}
In certain situations this theorem gives the most rapid test to
exclude closed 3-braids. For example, the knot $13_{6149}$ has
$MWF=3$, but seeing that $\Dl$ has $\Mc\Dl=-2$ we immediately
conclude that it cannot be a closed 3-braid.
\end{exam}

\begin{corr}
There are only finitely many 3-braid links with given $\Md\Dl$
(provided $\Dl\ne 0$). \qed
\end{corr}

It is actually true (as we will prove below) that for links with
\em{any} bounded braid index there are only finitely many different
Alexander polynomials of given degree admitted. However, such
Alexander polynomials may be admitted by infinitely many different
links (of that braid index). See \S\reference{s4.2} for some remarks.
% Using the non-faithfulness of
% the Burau representation on $\ge 5$ strands, it should not
% be too hard to find such families. The 4-strand case however
% seems open.
% 
% \begin{question}
% Are there only finitely many 4-braid links with given $\Dl\ne 0$
% (or equivalently $\Md\Dl$)?
% \end{question}

\begin{corr}\label{cbd}
Non-split 3-braid links bound no disconnected Seifert surfaces (with
no closed components). In particular there are no non-trivial 3-braid
boundary links.
\end{corr}

\proof For such links the Alexander polynomial is zero. The only
non-split 3-braid links $L$ of zero polynomial are closures of
$[(123)^k]$ for even $k$. These cases are easily ruled out by linking
numbers. Any pair of components of $L$ has non-zero linking number.
So a connected component of any Seifert surface $S$ of $L$ must have
at least two boundary components of $L$, and if $S$ is disconnected,
$L$ has at least 4 components, which is clearly not the case. \qed

The examples that falsify this claim for 4-braids are again
easy: consider the (closure links of) words in $[1\pm 2-1]$ and
$[2\pm 3-2]$.

\begin{corr}
If $\Dl$ is the Alexander polynomial of a 3-braid link, then
$|\Mc\Dl|\le \Md\Dl+2$. For a knot $|\Mc\Dl|\le\Md\Dl+1$.
\end{corr}

\proof If $|\Mc\Dl|\le 2$ then we are easily done. 
($\Dl=\pm 2$ cannot occur for a knot.) Otherwise we have
up to mirroring a strongly quasi-positive word $\bt$ reducing to
an even power of $[123]$. Now $|\Mc\Dl|$ counts non-trivial syllables,
so $[\bt]=3-\chi\ge 2|\Mc\Dl|$, and $-1-\chi=2\Md\Dl$. For a knot
the number of syllables with exponent $\ne 2$ is at least $2$,
so $[\bt]-2\ge 2|\Mc\Dl|$. \qed

Note in particular that the proof shows that $|\Mc\Dl|$
can be any given natural number, and how to find the link
that realizes this number. We emphasize this here, because
later we will prove a contrary statement in the case of
\em{alternating} links for \em{every} arbitrary braid index
(see corollary \reference{cSG}).

\begin{corr}
A 3-braid \em{knot} is fibered if and only if $\Mc\Dl=\pm 1$.
\end{corr}

\proof It remains to explain why no 3-braid knot has $\Mc\Dl=\pm 1$
but $\Md\Dl<1-\chi$. Latter condition would imply that we have
up to mirroring a strongly quasi-positive word $\bt$ reducing to
an even power of $[123]$, and former condition that $\bt$ has only one
non-trivial syllable. But $[(123)^{2k}]$ has 3-component closure,
and making one syllable non-trivial cannot give a knot. \qed

In particular, it is worth noting

\begin{corr}
No non-trivial 3-braid knot has trivial Alexander polynomial. \qed
\end{corr}

Again, the two 11 crossing knots immediately show that this is
not true for 4-braids.

\begin{exam}\label{x29}
We can also easily determine the 3-braid links for some small degree
Alexander polynomials. For example, we see that no other 3-braid
knot has the polynomial of $3_1$, $4_1$ or $5_2$. Similarly we can
check that no 3-braid knot has the polynomial of $9_{42}$ and
$9_{49}$ (which shows that these knots have braid index 4), a
fact we will derive in the last section using entirely different
representation theory arguments. 
\end{exam}

\subsection{\label{Js}Jones polynomial}

The control of the Jones polynomial on 3-braid links was the
object of main attention in \cite{Birman}. We can accomplish
this with a similar argument to $\Dl$. The result we obtain
can be conveniently described in our setting and is as follows:

\begin{theorem}\label{Hh}
Let $L$ be a non-split 3-braid link, and $L=\hat\bt$ with $\bt\in
B_3$. Then
\begin{eqn}\label{g}
\spn V(L)\,\le\,4-\chi(L)\,.
\end{eqn}
Equality holds if and only if $L$ is strongly quasi-signed (i.e.
-positive or -negative, or equivalently $|[\bt]|=3-\chi(L)$), and
not fibered. More specifically, the following holds:
\def\labelenumi{\theenumi)}
\begin{enumerate}
\item If $L$ is strongly quasi-positive, then 
$\ds\quad
\md V\,=\,\frac{1-\chi}{2}\quad\mbox{and}\quad\mc V\,=\,\pm 1\,.
$\\
Analogously, if $L$ is strongly quasi-negative, then
$\ds\quad
\Md V\,=\,\frac{\chi-1}{2}\quad\mbox{and}\quad\Mc V\,=\,\pm 1$\,.
\item If $L$ is strongly quasi-positive and fibered, then
$\spn V(L)\,\le\,3-\chi(L)$\,.
If $L$ is strongly quasi-positive and not fibered, then \eqref{g}
is an equality and $\Mc V\,=\,\pm 1$. (The properties for strongly
quasi-negative are analogous.)
\item If $L$ is not strongly quasi-signed and $|[\bt]|<1-\chi(L)$,
then 
\[
\mc V(L)\,=\,\pm 1\,, \quad\Mc V(L)\,=\,\pm 1\quad\mbox{and}\quad
\spn V(L)\,=\,3-\chi(L)\,.
\]
\item If $L$ is not strongly quasi-signed and $[\bt]=1-\chi(L)$, then
\[
\Mc V(L)\,=\,\pm 1\quad\mbox{and}\quad
\Md V(L)\,=\,\frac{5}{2}-\frac{3}{2}\chi(L)\,.
\]
Moreover, $\md V\ge \ds\frac{-1-\chi}{2}$, and if equality holds,
then $\mc V=\pm 2$. (The case $[\bt]=\chi(L)-1$ is analogous.)
\end{enumerate}
\end{theorem}

Apart from solving Birman's problem how to determine 3-braids
with given Jones polynomial, theorem \ref{Hh} easily implies that
no non-trivial 3-braid link has trivial (i.e. unlink) polynomial. We
defer the discussion of the non-triviality of the Jones polynomial
to \cite{ntriv}, where we work in the much more general context of
semiadequate links. In that paper we will show that semiadequate
links have non-trivial Jones polynomial. This result in fact motivated
theorem \reference{thq}, which then provides a different conclusion
about the non-triviality of the polynomial. A further application
will be the classification of the 3-braid links with unsharp
Morton-Williams-Franks inequality (mentioned in remark \ref{rqw}).

\begin{corr}\label{fintg}
For a given Jones polynomial $V$ (actually the pair
$(\md V,\Md V)$ is enough), there are at most three
values of $\chi(L)$ of a 3-braid link $L$ with $V(L)=V$.
If $\md V\cdot \Md V\le 0$, then $\chi(L)$ is unique.
\end{corr}

\proof The theorem shows that the value of $\chi$ is determined
by one of $\Md V$, $\md V$ or $\spn V=\Md V-\md V$.
In particular for a pair $(\md V,\Md V)$ there exist at most three
values of $\chi$ of 3-braid links with a Jones polynomial realizing
this pair. If $\md V\cdot \Md V\le 0$, then the options that a 3-braid
$\bt$ with $V(\hat\bt)=V$ is strongly quasi-signed or almost
quasi-signed are excluded (up to a few simple cases that can be
checked directly), so $\chi(L)$ is determined (unambiguously) by
$\spn V$. \qed

In particular, since 3-braid links of given $\chi$
are only finitely many, we have

\begin{corr}\label{finth}
There are only finitely many closed 3-braids with the same
Jones polynomial, actually with the same pair $(\md V,\Md V)$. \qed
\end{corr}

For example, one
easily sees that no 3-braid knot has the polynomial of $9_{42}$.
Similarly, no other 3-braid knot has the polynomial of the figure-8-knot
(there is, however, a 4-braid knot with such polynomial, $11_{386}$).

That there are only finitely many closed 3-braids with the same
skein polynomial was known from \cite{3br}. For the Jones polynomial
one should note that infinite families were constructed
by Traczyk \cite{Traczyk} if one allows polynomials up to units.
Traczyk's examples show that for fibered
(strongly quasi-)positive links $\spn V$ may remain the same
while $\chi\to-\infty$, so one cannot expect a full (lower) 
control on $\spn V$ from $\chi$. From these links, one obtains
by connected sum infinite families with the same
polynomial for 5-braids. (The status of 4-braids remains unclear.)
Also Kanenobu \cite{Kanenobu2} constructed finite families of 3-braids
of any arbitrary size, so that our result is the maximal possible.

A question that surfaces naturally with these remarks in mind is 

\begin{question}
Does for 3-braid links the Jones polynomial (or Alexander polynomial)
determine the skein polynomial? In other words, do any two
3-braid links with the same $V$ (or $\Dl$) have also equal $P$?
\end{question}

We have at least the following partial result, whose proof we
postpone after the proof of theorem \reference{Hh}.

\begin{corr}\label{CVP}
A given Jones polynomial $V=V(L)$ is realized for 3-braid links
$L$ by at most three different skein polynomials $P(L)$. If
$\md V\cdot \Md V\le 0$, then $V(L)$ determines $P(L)$.
\end{corr}

In the general case one cannot expect a positive answer to the
above question. At least for links there is now the method
of \cite{EKT} available, which should yield large families of
links with the same Jones but different skein (or Alexander)
polynomial. For constructing families with the same Alexander
(but different Jones or skein) polynomial, further techniques
are available, applicable also for knots, like the non-faithfulness
of the Burau representation (see \cite{Bigelow}) or tangle surgeries
(see \cite{Bleiler}).

% \begin{corr}\label{fintg}
% For a given Jones polynomial $V$ (actually the pair
% $(\md V,\Md V)$ is enough), there are at most three
% values of $\chi(L)$ of a 3-braid link $L$ with $V(L)=V$.
% If $\md V\cdot \Md V\le 0$, then $\chi(L)$ is unique.
% \end{corr}
% 
% \proof The theorem shows that the value of $\chi$ is determined
% by one of $\Md V$, $\md V$ or $\spn V=\Md V-\md V$.
% 
% In particular for a pair $(\md V,\Md V)$ there exist at most three
% % values of $\chi$ of 3-braid links with a Jones polynomial realizing
% this pair; such links are only finitely many. \qed

For the proof of theorem \ref{Hh} we use the previous work in
\cite{3br} (that in
particular answered Birman's question) on the skein polynomial. We
apply again the result in \cite{3br} that $\Md_zP=1-\chi$ for closed
3-braids. As in that paper, we distinguish the cases of band-%
positive, band-negative 3-braids, and such of Xu's form $L^{-1}R$. 

By the Morton-Williams-Franks inequalities \eqref{mq},
a 3-braid $\bt$ of exponent sum (writhe) $[\bt]$ has $v$-degrees of
$P$ in $[\bt]-2$ we call \em{left degree}, $[\bt]$ we call \em{middle
degree} and $[\bt]+2$ we call \em{right degree}. The terms of
$[P]_{z^k}$ for some $k$ in these degrees will be called left,
middle and right terms\footnote{Note that here the brackets for
polynomials and for braids have a completely different meaning.}. 

\begin{lemma}\label{lemg}\mbox{}\\
(a) If $\bt$ is band-positive, then $\Mc_zP$ has $v$-terms in the left
degree and possibly in the middle degree. If $\bt$ is band-negative,
$\Mc_zP$ has $v$-terms in the right degree and possibly in the
middle degree. In either situation a term in the middle degree occurs
if and only if $\hat\bt$ is not fibered.\\
(b) If $\bt$ is $L^{-1}R$
then $\Mc_zP$ has $v$-terms in the middle degree only.\\
In both cases all non-zero occurring coefficients are $\pm 1$.
\end{lemma}

\proof In (a) we prove only the first claim (the second
claim is analogous). Let $\chi=\chi(\hat\bt)$. Since
$\md_lP\le \Md_zP=1-\chi$ and $1-\chi=[\bt]-2\ge 
\md_lP$ by MWF, we have $\md_lP=\Md_zP$, and then \eqref{p21}
implies that $\Mc_zP$ has a $v$-term in the left
degree, with coefficient $\pm 1$. From \cite[Theorem 3]{3br} we
have then that it has no right-degree term, and that if it has
a middle-degree term, the coefficient is $\pm 1$. Now the
previous work and the substitution $v=1$, $z=t^{1/2}-t^{-1/2}$
for $\Dl$ in \eqref{DP} easily show that the middle term
occurs if and only if $\hat\bt$ is not fibered.

Now consider (b). If $\Mc_zP$ has a left-degree term, then using
the substitution \eqref{p21}, we saw that $\bt$ is band-positive.
Otherwise, $[\bt]<3-\chi$, so the contribution of the coefficient
of $z^{1-\chi}v^{[\bt]-2}$ in \eqref{p21}, which is not cancelled,
is not in degree $0$. Analogously one argues if $\Mc_zP$ has a
right-degree term. \qed

Note that keeping track of $[\bt]$ and distinguishing between
left and middle term is important here. By remark \ref{rqw},
we have links with equal polynomials such that the left
term of the one is the middle term of the other.

\begin{lemma}\label{lemldg}
(a) If $\bt$ is $L^{-1}R$ and $L$ has exponent sum $[L]>1$, then
$[P]_{z^{-1-\chi}}$ has a left-degree term, and the coefficient
is $\pm 1$. \\
(b) If $[L]=1$ and
$[P]_{z^{-1-\chi}}$ has a left-degree term, then the coefficient
is $\pm 2$. \\
The two analogous statements hold replacing $L$ by $R$
and left-degree term by right-degree term.
\end{lemma}

\proof We proved already that $\Mc_zP$ has only a middle-degree term.
So it must remain under the substitution of \eqref{p21} if and only
if $[\bt]=1-\chi$. The cases (a) and (b) occur when this term must
be cancelled, and complemented to $1$, respectively. \qed

\proof[of theorem \reference{Hh}]
The statements follow mainly by putting together the last two
lemmas and looking in which degrees the
non-cancelling contributions of the coefficients of $P$
occur under the substitution \eqref{VPsub}.

If part (a) of lemma \reference{lemg} applies, we established already
(in theorem \reference{thqp})
that the positive band form of $\bt$ is equivalent to the strong
quasi-positivity of $\hat\bt$. In this case $\md V$ comes from
the left-degree term in $\Md_zP$. If $\hat\bt$ is not
fibered, then $[P]_{z^{-1-\chi}}$ must have a right-degree
term (with coefficient $\pm 1$) to cancel the middle-degree
term of $\Md_zP$ under \eqref{p21}.

If part (b) of lemma \reference{lemg} applies, we use the
further information of lemma \ref{lemldg}.
In case (a) of lemma \ref{lemldg}, the left and right
terms in $[P]_{z^{-1-\chi}}$ determine the degrees
and edge coefficients of $V$.

In case (b) of lemma \ref{lemldg},
the maximal term in $V$ comes from the right-degree term in
$[P]_{z^{-1-\chi}}$. A coefficient in $t^{(-1-\chi)/2}$
may come only from a left-degree term in $[P]_{z^{-1-\chi}}$,
which, if occurring, is with coefficient $\pm 2$. \qed

This proof underscores the significance of \eqref{VPsub}
as a tool for studying the Jones polynomial. So far it
seems to have been useful just for calculating specific
Jones polynomials from $P$. In \S\reference{Sfg} we will
see further results that come out of considering this
substitution.

We note the following equalities that follow from the proof of
theorem \ref{Hh}. These will be needed in the study of the Q
polynomial, and are also helpful for corollary \reference{CVP}.

\begin{lemma}\label{lemma_}
If $\bt$ is a not strongly quasi-signed 3-braid of exponent sum $e$
and $V$, $\chi$ the Jones polynomial resp. Euler characteristic
of its closure, then
\begin{eqnarray}
\label{one_}\Md V-e & = & \frac{1-\chi}{2}+1 \\
\label{two_}\md V-e & = & -\frac{1-\chi}{2}-1 
\end{eqnarray}
\em{Exceptions} are \eqref{two_} if $\bt$ is strongly
almost quasi-positive, and \eqref{one_} if $\bt$ is strongly
almost quasi-negative. \qed
\end{lemma}

\proof[of corollary \reference{CVP}]
Depending on whether a 3-braid representation $\bt$ is strongly
quasi-signed, almost quasi-signed, or none of both, $V$ determines
$e=[\bt]$ via $\chi(\hat\bt)$ and/or \eqref{one_} or \eqref{two_}.
With $V$ and $e$, one can recover the trace of the Burau matrix
$\psi_2(\bt)$, and from that also $P(\hat\bt)$ (see \cite{Birman} or
\eqref{V_3} below). Again, if $\md V\cdot \Md V\le 0$, then the
options that $\bt$ is strongly quasi-signed or almost quasi-signed
are ruled out easily. \qed

\subsection{Q polynomial}

We extend the scope of the previous results to the
Brandt-Lickorish-Millett-Ho polynomial. The main aim here is to prove

\begin{theo}\label{thQ}
Only finitely many non-trivial 3-braid links have given $Q$
polynomial, and none has trivial (i.e., unlink) polynomial.
\end{theo}

Note that by Kanenobu's work \cite{Kanenobu3}, again
there are finite families of arbitrary large size,
so we claim again a sort of contrary result.

While this result may be considered less relevant
that its analoga for $\Dl$ and $V$, its proof displays
the largest variety of tools necessary to apply, and 
shows how the various approaches to link invariants
(skein relations, state models and representation theory), 
which are often considered in isolation, can usefully 
complement each other. Indeed, $Q$ seems in general 
more difficult to treat than $V$. Apart from Kidwell's 
results \cite{Kidwell} for alternating links, and those
in \cite{pos} for positive \em{knots}, neither the 
non-triviality nor the finiteness property seem to have 
been known previously for any other class of links.

The proof makes use of the full extent of the study of 
Xu's form. Various additional ingredients will be necessary.
One such is J. Murakami's $Q$-$V$ formula, found in
\cite{Murakami} (by representation theory, and proved also
later by Kanenobu \cite{Kanenobu} using skein relations).
It allows to recur in many cases the problem from $Q$ to $V$.
When dealing with $V$, we use beside the previous discussion
the (Kauffman bracket based) formulas of \S\ref{Ssaq}.

In cases Murakami's formula is not helpful, we apply the
Polyak-Viro formula \eqref{pv2} for Casson's knot invariant
$v_2=\myfrac{1}{2}\Dl''(1)$,
and the following formula, (likewise skein theoretic and) due
to Kanenobu \cite{Kanenobu4}, that relates $v_2$ to $Q$.

\begin{theo}(\cite{Kanenobu4})
For a link $L$ with $n$ components $K_i$, $i=1,\dots,n$, we have
\begin{eqn}\label{Kanf}
Q_L'(-2)=3(-2)^n\sum_iv_2(K_i)+3(-2)^{n-2}\sum_{i<j}lk(K_i,K_j)^2
+(n-1)(-2)^{n-3}\,,
\end{eqn}
where $lk$ is the linking number.
\end{theo}

Using this formula, we prove first a special case of theorem
\ref{thQ}.

\begin{prop}\label{pgw}
Only finitely many (non-trivial) strongly quasi-positive 3-braid
links have given $Q$ polynomial, and none has trivial polynomial.
\end{prop}

For the case of knots, we have the following more specific statement,
that has also some independent meaning.

\begin{prop}\label{pgq}
If $K$ is a strongly quasi-positive 3-braid knot, then $v_2(K)\ge
g(K)$.
\end{prop}

Note that in \cite{pos} we proved the same inequality for a
general \em{positive} knot. Since the pretzel knots with
trivial Alexander polynomial are strongly quasi-positive, nothing
like this holds for a general strongly quasi-positive knot, though.

We need some preparations. Let us first fix a convention. We number
braid strands by 1,2,3 from left to right in the bottom of the
braid by and compose words from bottom to top. We propagate strand
number through the crossings. In particular strand( number)s may 
appear permuted in the middle or on top of the braid. To refer to the 
ordering of the strands with regard ot the local position in the braid
diagram, we speak of left, middle and right strand. So, for example, 
if strands 1 and 2 enter below a $\sg_1$ (strand 1 as left strand and
strand 2 as middle strand), then they exit above in the order 2\ 1 
(strand 2 as left strand and strand 1 as middle strand). 

Let for a (not necessarily pure, and not closed) braid word the
linking numbers of strand $i$, $j$ numbered as explained above be
the sum of the writhes of all crossings these two strands pass.
(Note that for links this differs by the additional factor 1/2,
which is not relevant, though.)

\begin{lemma}\label{lmlk}
Let $\bt$ be a strongly quasi-positive 3-braid word
(not necessarily pure, and not closed). Then the
linking numbers $l_{ij}$ between pairs $(i,j)$ of 
strands satisfy $l_{13}+l_{23}>0$, unless $\bt$ is split.
\end{lemma}

\proof The sum of writhes of crossings a strand passes
through each of $\sg_{1,2,3}$ is non-negative. So
$l_{13}+l_{23}\ge 0$. If $l_{13}+l_{23}=0$, then
strand 3 passes only from below as left or middle
strand into a $\sg_3=\sg_1^{-1}\sg_2\sg_1$. But since it
starts (at the bottom of the braid) as a rightmost 
strand, this means that is passes no crossing, so $\bt$ is split. 
\qed

We prove proposition \ref{pgq}
by induction on the length of the positive word
in $\sg_{1,2,3}$. Assume $\bt$ is written as such a word.
We can w.l.o.g. cyclically permute the indices. 

\begin{lemma}
If $\hat\bt$ is a strongly quasi-positive 3-braid, then a positive 
word $\bt$ can be chosen so that it contains a non-trivial syllable,
\em{unless} $\bt$ is a power of $[123]$. In particular, it is 
always possible if $\hat\bt$ is a knot, or a fibered link.
\end{lemma}

\proof Choose Xu's form. It is $[12]^k\bt'$, where $k\ge 0$ and $\bt'$ 
has cyclically non-decreasing indices. This word $\bt$ has a non-trivial
syllable up to cyclic permutations of the letters, unless $\bt'$
is a power of $[123]$. In this case, if $k>0$, we apply a YB relation 
at the initial $121$ in $\bt$, and are done. \qed

\proof[of proposition \ref{pgq}] Let $\hat\bt$ be a knot. 
We can w.l.o.g. now assume $\bt$ is written as a word with a 
non-trivial syllable $s$. Now we apply the skein relation
of $\Dl$ at a letter of $s$. Let $K=K,K_-,K_0$ be the skein
triple, and $\bt_{\pm,0}$ the corresponding braids. Also we can 
cyclically permute the indices in $\bt_{\pm,0}$. We choose them so 
that the strand fixed by the permutation of $\bt_0$ is number 3. We 
know from the skein relation of $\Dl$ 
that $v_2(K)-v_2(K_-)=lk(L,M)$ where $L,M$ are the components 
of $K_0$. Now that linking number is positive by lemma \ref{lmlk} 
(unless $\bt_0$ is a split braid, but this is impossible because $s$
is non-trivial). By induction on $g(K)$, the claim of proposition
\ref{pgq} follows. \qed

\proof[of proposition \reference{pgw}] We distinguish three
cases depending on the number of components of $\hat\bt$.

\begin{caselist}
\case 1 component. In this case $Q'(-2)$, by Kanenobu's formula
\eqref{Kanf},
is a multiple of Casson's invariant $v_2$. The claim we wish to 
show then is implied by the estimate in proposition \reference{pgq}.

\case 2 components. In this case Kanenobu's formula 
involves only the square of the linking number of both components 
$K,L$ of the closure link of $\bt$ and $v_2$ of the 2-string subbraid 
component $K$. Looking at Xu's normal form, we see that if $\bt$ is 
a strongly quasi-positive 3-braid, then writing $\bt$ as a word
in $\sg_{1,2,3}$, and then expanding $\sg_3=\sg_1^{-1}\sg_2\sg_1$,
we obtain a word in $\sg_{1,2}$ of some length $c$ with at most 
$c/5$ negative letters. Since thus $\bt$ has at most $c/5$ negative
crossings, it is easy to see that not both $v_2(K)$ and $lk(K,L)$
can be zero, and one grows unboundedly when $c$ grows. So we are done.

\case 3 components. In this case Kanenobu's formula \eqref{Kanf} 
reduces to the square sum of the linking number of pairs components.
Let $l_{12},l_{13},l_{23}$ be these linking numbers. Again it is 
clear that not all of $l_{12},l_{13},l_{23}=0$, and that when 
$l_{12}+l_{13}+l_{23}\ge 2c/5$ grows unboundedly, then so will
$l_{12}^2+l_{13}^2+l_{23}^2$, so again we are done.
\end{caselist}
\qed

Now we move on to settle the remaining cases in theorem \ref{thQ}.
Here the tools differ considerably. Beside the study of Xu's form,
we need the help of the following theorem, proved in \cite{ntriv}.

\begin{theorem}\label{thq}
Any $3$-braid is up to conjugacy $A$-adequate or $B$-adequate.
\end{theorem}

We proved also the following

\begin{theorem}
A $3$-braid word is $B$-adequate if and only if it is (a) a
negative word, or (b) it contains no $[121]$ as subword, and negative
entries in the Schreier vector appear isolated. (That is, in
cyclic order they are preceded and followed by a positive entry.)
\end{theorem}

Using this criterion we prove now

\begin{lemma}\label{oia}
Assume $\bt$ is strongly $\le 2$-almost positive and of exponent 
sum $e\ge 0$. Then $\bt$ is $B$-adequate.
\end{lemma}

\proof
If $\bt$ is strongly quasi-positive or almost strongly
quasi-positive, then it is of Xu's form $L^{-1}R$
where $L$ is a single letter. By direct observation
we can verify that when permuting indices cyclicly
properly and writing out $\sg_3$ in $\sg_{1,2}$, then
the resulting word has no $[121]$ subword, and negative
entries in the Schreier vector appear isolated. (In case
of $[(312)^k31-2]$, we must cancel 2 crossings first.)

If $\bt$ is 2-almost strongly quasi-positive, then then
it is of Xu's form $L^{-1}R$ where $L$ has 2 letters.
We can assume that $L=[12]$ or $L=[1^2]$, then we write down
the subword of $L^{-1}R$ consisting of the first and last
letters of $R$ and $L$ (herein $k\ge 1$):
\[
\begin{array}{ll}
[[23] \dots 33-2-1] & [[23] \dots 22-1-1] \\{}
[[23] \dots 23-2-1] & [[23] \dots 31^k2-1-1] \\{}
[[23] \dots 31^k-2-1] & [[23] \dots 33-1-1] \\{}
& [[23] \dots 23-1-1]
\end{array}
\]
(Here `[23]' means as in \S\ref{brs} that the word is to begin with
`2' or `3'.) The reductions and $B$-adequacy test are done
case-by-case. \qed

Using the formulas in \S\reference{Ssaq}, we prove

\begin{lemma}\label{yu}
If $L=\hat\bt$ is a link, which is the closure of an $A$-semiadequate 
$3$-braid $\bt$, and $V_0V_1=-1$, then $V_0V_2=1$ or $2$.
\end{lemma}

\proof By \cite{ntriv} a A-semiadequate $3$ braid is either
positive, in which case $V_1=0$, or has no $[-1-2-1]$ and
positive entries are isolated in the Schreier vector.
In \cite{ntriv} we described the words  for latter
type which have $V_1=0$. They are (up to braid relations
and cyclic letter permutations) of the form%
\footnote{Note that these are exactly the 3-braids which are
reducible in the terminology of dynamic properties;
this reducibility has nothing to do, though, with the
reducibility in Markov's theorem.}
\[
[1^{-2}2^{-2}1^{-2}\dots 1^{-2}2^{-1}1^p2^{-1}]\,.
\]

Analogously to that study, we can see that if $V_0V_1=-1$,
then $\bt$ has a word which is obtained from the above type
by replacing some $1^{-2}$ by $1^{-k}$ 
for $k\ge 3$. The intertwining graph of
the $A$-state is a path (all vertices have valence 2, except two
of valence 1), so $\chi(IG)=1$. We have $\bigtriangleup
=1$ if $k=3$, and otherwise $\bigtriangleup=0$. So by
theorem \reference{T02}, we have $V_0V_2=1$ for
$k=3$, and $V_0V_2=2$ otherwise. \qed

To explain what such a property has to do with $Q$, 
now we must introduce J. Murakami's formula \cite{Murakami}.
(See also Kanenobu \cite[Theorem 2]{Kanenobu}.)

Let $i=\sqrt{-1}$, $u=\sqrt{-t}$ and $x=u+u^{-1}$. Let further%
\footnote{Note that $V$ is here what is written as $J$, and not $V$, 
in \cite{Kanenobu}.} for a braid $\bt$ of exponent sum $e$,
\begin{eqn}\label{chi}
\chi(\bt,t)=i^eu^{-2e}V_{\hat\bt}(t)+u^{-e}(x^2-2)\,.
\end{eqn}

Then Murakami's formula is

\begin{theo}(J. Murakami)\label{TMF}
If $L$ is a closure of a 3-braid $\bt$ of exponent sum $e$,
then
\begin{eqn}\label{sonjo}
Q(L,x)\,=\,\chi(\bt,\sqrt{t})^2-1\,+\,\frac{2(x^2+x-1)}{x^2(x^2-3)}
(u^e+e^{-e})+\frac{-x^4-2x^3+3x^2+4x-4}{x^2(x^2-3)}\chi(\bt,t)\,.
\end{eqn}
\end{theo}

\proof[of theorem \ref{thQ}] We can w.l.o.g. (taking the mirror image)
assume that $e\ge 0$, and (excluding trivial to check special cases)
that $\chi\le 0$. Clearing denominators and absolute terms in
\eqref{sonjo}, we find
\begin{eqnarray}
\label{Q1} Q_1(L,x):=x^2(x^2-3)(Q(L,x)+1) & = &
\bigl[ i^e(-t)^{-e}V(t)+(-t)^{-e/2}(-t-\frac{1}{t}) \bigr]
\cdot (-x^4-2x^3+3x^2+4x-4) \\
\nonumber & &  + \bigl[ i^e(-t)^{-e/2}V(-\sqrt{-t})+(-t)^{-e/4}x \bigr]^2
\cdot x^2(x^2-3) \\
\nonumber & & +2(x^2+x-1)(\sqrt{-t}^e+\sqrt{-t}^{-e})\,.
\end{eqnarray}

Rearranging, we need to show that the sum of the following 5
expressions,
regarded as a polynomial in $t^{\pm 1/2}$, has either arbitrarily
small minimal degree or arbitrarily large maximal degree.
(Note that by \eqref{sonjo} this sum must be self-conjugate in 
$t$ up to coefficient signs, which is not at all evident directly.)

\begin{eqnarray*}
T'_1 & = & i^e(-t)^{-e}V(t)\cdot (-x^4-2x^3+3x^2+4x-4) \\
T'_2 & = & t^{-e}\bigl[V(-\sqrt{-t})\bigr]^2\cdot x^2(x^2-3) \\
T'_3 & = & 2i^e(-t)^{-3e/4}V(-\sqrt{-t})x^3(x^2-3) \\
T'_4 & = & (-t)^{-e/2}\bigl[(x^2-2)(-x^4-2x^3+3x^2+4x-4)+x^4(x^2-3)+2(x^2+x-1)
 \bigr] \\
T'_5 & = & (-t)^{e/2}\bigl[2(x^2+x-1) \bigr] 
\end{eqnarray*}

Switching $-t\to t$ to simplify the expressions, and regrouping terms, 
we have
\begin{eqnarray*}
T_1 & = & i^et^{-e}V(-t)\cdot (-1\ -2\ -1\ -2\ [-4]\ -2\ -1\ -2\ -1) \\
T_2 & = & (-t)^{-e}\bigl[V(-\sqrt{t})\bigr]^2\cdot (1\ 0\ 1\ 0\ [0]\ 0\ 1\ 0\ 1) \\
T_3 & = & 2i^et^{-3e/4}V(-\sqrt{t})\cdot (1\ 0\ 2\ 0\ 1\ [0]\ 1\ 0\ 2\ 0\ 1) \\
T_4 & = & t^{-e/2}\bigl[-2(x-1)(x^2-3)(x^2-1) \bigr] \\
T_5 & = & t^{e/2}\bigl[2(x^2+x-1) \bigr] 
\end{eqnarray*}
Now $u=t^{1/2}$ and $x=t^{1/2}+t^{-1/2}$. The last factors for $T_
{1,2,3}$ are given as a list of coefficients, with the absolute term
put in brackets. Note that these are polynomials in $\sqrt{t}$, so for
example, the first two polynomials have degree $2$ in $t$.

Since we are concerned with cancellations of the leading and trailing
coefficients of $T_k$, let us compile their minimal $m_k=\md_tT_k$ and
maximal degrees $M_k=\Md_tT_k$ (note that by taking degrees w.r.t.
$t$, and not $\sqrt{t}$, the $m_k$, $M_k$ are only half-integers).
\[
\begin{array}{c|c|c}
\raisebox{-1em}{\ry{2.3em}} k & m_k & M_k \\
\hline
\raisebox{-1em}{\ry{2.5em}}1 & \md V-e-2 & \Md V-e+2 \\[2mm]
2 & \md V-e-2 & \Md V-e+2 \\[2mm]
3 & -\ffrac{3e}{4}+\ffrac{\md V}{2}-\ffrac{5}{2} & -\ffrac{3e}{4}+
    \ffrac{\Md V}{2}+\ffrac{5}{2}\\[2mm]
4 & -\frac{e}{2}-\frac{5}{2} & -\frac{e}{2}+\frac{5}{2} \\[2mm]
5 & \frac{e}{2}-1 & \frac{e}{2}+1\raisebox{-1em}{\ry{1.7em}}
\end{array}
\]
(Here again $\md V=\md_t V(t)$ is a half-integer, and similarly
$\Md V$.)

Our aim will be to determine what of the $M_k$ is the largest, and 
to show that the coefficients of the contributing $T_k$ do not cancel.
An important case where we exclude (problematic) cancellation is

\begin{lemma}\label{lmcan}
Assume that 
\begin{eqnarray}
\label{GG} M_1=M_2 & \ge & M_k+3/2\quad\mbox{ and } \\
\label{gg} m_1=m_2 & \le & m_k-3/2
\end{eqnarray}
for $k=3,4,5$.
Then $\Md_t Q_1(L,x)\ge M_1-1$. (Here $Q_1(L,x)$ refers to \eqref{Q1}.)
If $L$ the closure of a 3-braid which is $A$-adequate resp. $B$-adequate
then assuming \eqref{gg} resp.\ \eqref{GG} alone is sufficient.
\end{lemma}

\proof By theorem \ref{thq}, $L$ is the closure of a 3-braid $\bt$ which 
is either $A$-adequate or $B$-adequate. We consider the $B$-adequate case 
and show that $M_1=M_2\ge M_k+3/2$ implies $\Md_t Q_1(L,x)\ge M_1-1$.
If $\bt$ is $A$-adequate, we obtain similarly $\md_t Q_1(L,x)\le 
m_1+1$, and the result follows by the (anti)symmetry of $Q_1$.

Now if $M_1=M_2\ge M_k+3/2$, the last three coefficients of $M_{1,2}$
are not cancelled from other $M_k$ (remember these are polynomials
in $\sqrt{t}$, and coefficients are meant for such polynomials), 
so it is enough to show that they do not (completely) cancel 
among each other. Now a look at the formulas for $T_{1,2}$ shows
that if these three coefficients are all to cancel, then (apart from
proper sign coincidences) we must have $\bar V_1\bar V_0=-1$ and 
$\bar V_2=0$. (Take also into account the switch $-t\to t$.)
However, this situation was ruled out by lemma \ref{yu}. \qed

We assumed $e\ge 0$; we also excluded the case $e=3-\chi$ of
strongly quasi-positive braids. To apply the lemma, we need to
establish \eqref{gg} and/or \eqref{GG}. For this we use now
lemma \reference{lemma_}.

Clearly, a strongly almost quasi-negative braid should not 
be considered for $e\ge 0$; the almost quasi-positive braid
require a small extra argument, which is given below. Now
\begin{eqn}\label{hj}
\Md V-e\,=\,\frac{1-\chi}{2}+1\,\ge\,\frac{e}{2}+1\,,
\end{eqn}
so $M_3\ge M_4$, and we can discard $M_4$. Also from \eqref{hj}
we obtain
\[
M_1-M_5-1=\Md V-\frac{3e}{2}\ge 1\,,
\]
so $M_1>M_5+1$, and can neglect $M_5$. Similarly, if
$\bt$ is not almost strongly quasi-positive, we have with 
\eqref{two_} and $\md V<e/2$ that
\begin{eqn}\label{mmm}
m_3< m_4<m_5\,,
\end{eqn}
so $m_{4,5}$ are also irrelevant.

So we need to deal only with $m_k,M_k$ for $k\le 3$. Moreover, we 
see from \eqref{one_} that $M_1\ge (1-\chi)/2+3$. So whenever we can 
apply lemma \ref{lmcan}, we have for $\chi\le 0$ that $\Md Q_1\ge 
M_1-1>2$, so $\Md_z Q>0$, and $Q$ is not an unlink 
polynomial. Moreover, we have for a sequence of links
from $\chi\to -\infty$ also $\Md M_1\to \infty$, so
$\Md Q\to \infty$, and a given (even just degree of a)
$Q$ polynomial occurs only finitely many times, as desired.
A similar argument applies if we use \eqref{two_}.

To apply lemma \ref{lmcan}, we need to check $M_3<M_1-1$ and $m_3>m_1+
1$. Note that if $\bt$ is almost strongly quasi-positive, then \eqref
{two_} may not hold. So may not be able to apply lemma \ref{lmcan}
directly. However, we remedied this by showing in lemma \ref{oia}
that $\bt$ is $B$-adequate. Then \eqref{one_} holds, and it is enough
to use $M_3<M_1-1$ (the condition \eqref{mmm} becomes also irrelevant).

Now using \eqref{one_} and \eqref{two_}
\begin{eqnarray}
\nonumber
m_1-m_3+1 & = & \md V -e -2 +\frac{3e}{4}-\frac{\md V}{2}+\frac{5}{2}+1
\,=\,\frac{\md V}{2}-\frac{e}{4}+\frac{3}{2} \\
\nonumber
& = & \frac{e}{2}-\frac{1-\chi}{4}-\frac{1}{2}-\frac{e}{4}+\frac{3}{2}\,=\,
\frac{e}{4}-\frac{1-\chi}{4}+1\,, \\
\label{iua}
M_1-M_3-1 & = & \frac{\Md V}{2}-\frac{e}{4}-\frac{3}{2} \\
\nonumber
& = & \frac{e}{2}+\frac{1-\chi}{4}+\frac{1}{2}-\frac{e}{4}-\frac{3}{2}\,=\,
      \frac{1-\chi}{4}+\frac{e}{4}-1\,.
\end{eqnarray}

Then \eqref{GG} becomes equivalent to
\[
\frac{1-\chi}{4}+\frac{e}{4}>1\,.
\]
So, assuming $\chi\le 0$ and $e\le 1-\chi$, for \eqref{GG} 
we are left to consider only $(e,\chi)=(0,-1)$, $(0,-3)$,
$(1,0)$, $(1,-2)$, or $(2,-1)$, which are trivial cases to check. 

For \eqref{gg} we must deal with $A$-adequate but not $B$-adequate
braids with $e\ge -3-\chi$, i.e. $\le 3$-almost strongly
quasi-positive braids. For $\le 2$-almost strongly
quasi-positive braids we can use lemma \reference{oia}.
While one may extend the argument there to $3$-almost strongly
quasi-positive braids, it would require a longer
case-by-case analysis. We provide instead a different argument.

We try to modify the argument for proposition \reference{pgw}.
In case $\hat\bt$ is a link of 2 or 3 components, we can argue again
as there, and are left with the cases $e=0,1$ of $6$ and $7$ bands,
which are easy to deal with.

So assume $\hat\bt$ is a knot.
Now from \eqref{iua}, we have $m_3=m_1+1$, so if still some of the
first two coefficients of $T_1$, $T_2$ do not cancel, we have
$\md_t Q_1\le m_1+\myfrac{1}{2}$, and are done.
Otherwise again $V_1V_0=-1$. So by the proof of lemma \reference{yu},
\[
\bt=[1^{-k}2^{-2}1^{-2}\dots 1^{-2}2^{-1}1^l2^{-1}]\,,
\]
with $k\ge 3$ and $l\ge 1$. Since $\hat\bt$ is a knot, $k$ and $l$
are easily observed to be odd. Now a $\le 3$-almost positive
word in Xu's form can be chosen, when written out as a word in
$\sg_{1,2}$ of $c$ letters, to have $c_-\le c/5+3$ negative
letters/crossings. If we reduce the word by cancelling
crossings, this inequality remains true.

We estimate the Casson invariant $v_2$ of $\hat\bt$ using the Gau\ss{}
diagram formula of Polyak-Viro \eqref{pv2}. Let us put the
basepoint right after the group of $l$ positive crossings (which
of either strands is immaterial). 

Now $l$ is odd and the $l$ positive crossings are pairwise linked,
so their contribution to the Gau\ss{} diagram sum \eqref{pv2} is,
independently on the location of the basepoint, equal to $v_2(T_{2,l})
=\ffrac{l^2-1}{8}$, where $T_{2,l}$ is the $(2,l)$-torus knot. The
analogous claim is true for the group of $k$ negative crossings.

If some of the positive crossings is linked with a negative crossing,
it must be a negative crossing of a $\sg_2^{-1}$. Now by our choice
of basepoint, in
each syllable $\sg_2^{-2}$ only one of the two crossings gives a pair
that contributes to \eqref{pv2}, and the contribution is always
$-1/2$. Putting this together, we have
\begin{eqn}\label{oiz}
v_2\,\ge\,\frac{l^2-1}{8}+\frac{k^2-1}{8}-\frac{c_--k+2}{4}\cdot l
\cdot \frac{1}{2}\,.
\end{eqn}
Assuming $c\ge 10$ (the other cases are checked directly), we have
\[
l\,\ge\frac{4c}{5}-3\,\ge\,\frac{c}{5}+3\,\ge\,c_-=c-l\,.
\]
So the r.h.s. in \eqref{oiz} is minimized for $l\,\ge\ffrac{4c}{5}-3$
when putting $l=\ffrac{4c}{5}-3$. Then we have, with $k\ge 3$ and
$c_-=c-l$, the estimate
\[
v_2\,\ge\,\frac{3c^2}{50}-\frac{29c}{40}+\frac{11}{4}\,,
\]
which is positive for $c\ge 10$, and grows when $c\to \infty$. 

With this theorem \ref{thQ} is proved. \qed

\begin{rem}
It is interesting whether $\Md Q\to\infty$ also for the strongly 
quasi-positive braids, but in Murakami's formula massive cancellations 
become possible and $Q_1$ cannot be easily controlled.
At least one can prove using theorem \ref{thq} and
some results of Thistlethwaite in \cite{Thistle} that
there are only finitely many 3-braid links with given
$\Md_zF$, where $F$ is the Kauffman polynomial.
\end{rem}

\section{Positivity of 3-braid links}

\subsection{Positive braid links\label{S2.2}}

\subsubsection{The Morton-Williams-Franks bound}

For the proof of theorem \reference{thg} we will need to study
the behaviour of the bound \eqref{mwfb} in the Morton-Williams-Franks
inequality (we abbreviate as MWF) on positive
braids. This was begun by Nakamura \cite{Nakamura},
who settled the case $MWF=2$ in the suggestive way: such braids
represent only the $(2,n)$-torus links. (The case $MWF=1$ is
trivial.) We will introduce a method that considerably simplifies
his proof (but still makes use of some of his ideas), and
then go on to deal with $MWF=3$. The example of non-sharp
MWF inequality, $13_{9365}$ in \cite{KnotScape} (the connected
2-cable of the trefoil), given in \cite{MorSho}, is in fact
only among a small family of exceptional cases.

\begin{theorem}\label{thmwf}
If $b(\hat\bt)>MWF(\hat\bt)=3$, then $\bt$ reduces to a 4-braid,
and is given by one of the following forms (assuming that $3*$
denotes a sequence of at least one letter $3$, and $11*$
resp.\ $22*$ sequences of at least two letters $1$ or $2$):\\
\[
[22*3*122*11*23211*], \quad\mbox{and}\quad [22*3122*11*23*211*].
\]
\end{theorem}

We will reduce the proof to a finite number of words to
check, which is done by calculation using the program of
\cite{MorSho}. Since a direct computation is more reliable
than an increasingly difficult mathematical argument, we have
not tried to minimize the calculation by all means. However,
we point out that for the sake of theorem \reference{thg}
alone (rather than its refinement, theorem \reference{thmwf}),
the following weaker statement is sufficient, for which
a considerable part of the case-by-case calculations can
be dropped. This corollary requires the notion of
semiadequacy \cite{LickThis}, and can be deduced from theorem
\reference{thmwf} by direct check of the exceptional words.
(We will sometimes write $MWF(\bt)$ for $MWF(\hat\bt)$.)

\begin{corr}\label{cort}
If $\bt$ is a positive braid word, and $MWF(\bt)=3$, then $\bt$
reduces (up to Markov equivalence) to a positive 4-braid word $\bt'$,
and the diagram $\hat\bt'$ is not $B$-adequate.
\end{corr}

Here the notions of $A$-adequate and $B$-adequate for diagrams
and braids are as explained in \S\ref{brD} and \S\ref{brt}. We
note, as a consequence of \cite{Thistle}, that a braid is ($A$/%
$B$-)adequate if and only if some, or equivalently any,
minimal length word of a braid in its conjugacy class it is so.

We will thus prove theorem \reference{thg} only using corollary
\reference{cort}, and indicate in the proof of theorem \reference{thmwf}
the point where the corollary follows (and the rest of the argument
is not needed). The argument that elegantly replaces the remaining
case-by-case checks requires theorem \reference{thq}. Since (by
taking again the full extent of our proof) theorem \reference{thq}
is nonetheless not indispensable, we permit ourselves to defer its
proof to a separate paper \cite{ntriv}.

% \proof It is easy to see that each positive irreducible $\bt\in B_4$
% with a $\sg_2$ can be written as $\sg_2^{k_1}\sg_1^{l_1}\sg_3^{m_1}\dots
% \sg_2^{k_n}\sg_1^{l_n}\sg_3^{m_n}$, where $m_i>0$ and $l_i+m_i>0$.
% 
% Now if $\bt$ contains no YB relation, then all if $l_i+m_i>1$,
% and $n>1$ or $l_1,m_1>1$.
% 
% Now if all syllables of $\sg_2$ are non-trivial, we can remove
% in the positive resuloution tree and have the split union of
% two non-trivial (because $\bt$ irreducible) $(2,k)$-torus links,
% so $MWF=4$$. If ....
% \proof Check directly. ($B$-adequacy is invariant under regular isotopy,
% so any other 4-braid word giving the same knot is not $B$-adequate
% either.) \qed

\proof[of theorem \reference{thg}]
If $MWF(\bt)\le 2$, then we are done. So assume $MWF(\bt)=3$.
Assume first $L$ has an $A$-adequate 3-braid. Since $A$-adequate
diagrams minimize the number of negative crossings, and $L$ has
a positive (braid) diagram, the $A$-adequate 3-braid diagram is
positive, and we are done. So let $L$ have a $B$-adequate 3-braid
diagram $\hat\bt'$. Now by corollary \reference{cort}, we find
that $L$ reduces to a positive 4-braid $\bt$
with $c_+(\bt)$ positive crossings, and $\bt$ is not $B$-adequate.
Since $\bt'$ is $B$-adequate and $B$-adequate diagrams
minimize the number of positive crossings, $c_+(\bt')<c_+(\bt)$.
On the other hand, $MWF(L)=3$ and the inequalities of MWF for the
$v$-degree of $P$ show that $\bt'$ must have exponent sum
$[\bt']=c_+(\bt)-1$. Since $[\bt']\le c_+(\bt')$, we must have
equality, so $\bt'$ is positive. \qed

Note the following easy and useful consequence of theorem 
\reference{thg}:

\begin{corr}
A link which is a closure of a positive braid of at most 4 strings,
has a minimal crossing diagram as a closed positive braid, and a
minimal string positive braid representation.
\end{corr}

\proof The case of the braid representation is straightforward,
and it implies the minimal crossing diagram statement by
looking at $\Md_zP$ and using Morton's inequalities. \qed

The examples in \cite{posex_bcr}, mentioned after theorem \ref{thg},
show that the corollary is not true in case of positive 5-braids at
least for the positive minimal braid representation. So far no examples
are known where no minimal crossing positive braid diagram exists
(it was known to exist from \cite{WilFr,Murasugi} for closed
positive braids with a full twist, which include the torus links,
and from \cite{posex_bcr} for positive braid knots of at most
16 crossings), but the pathologies for minimal strings hint to caution.
It was shown in \cite{posex_bcr} that one, and in \cite{ntriv} that
infinitely many fibered positive knots
have no minimal crossing positive diagram.

\subsubsection{Maximal subwords}

Here we start the technical considerations needed to prove 
theorem \reference{thmwf}. We consider the form \eqref{wd},
now with all $l_i>0$.

\begin{defi}
We define \em{summit syllables} in \eqref{wd}:
\def\theenumi{\alph{enumi}}
\def\labelenumi{\theenumi)}
\begin{enumerate}
\item All $\sg_{n-1}^{l_i}$ are summit syllables, and
\item if $\ap=\sg_k^{l_i}$ and $\ap'=\sg_k^{l_j}$ are summit
syllables, with no $\sg_{k'}^{l'}$ for $k'\ge k$ occurring between
$\ap$ and $\ap'$, then all $\sg_{k-1}^{l'}$ occurring between $\ap$
and $\ap'$ are summit syllables.
\end{enumerate}
\end{defi}

Note that, according to definition \reference{dfext}, we
consider syllables in cyclic order. The relation ``between'' in
the above definition should also be understood in that sense:
a syllable occurring after the last index $i$ syllable
$\ap$ and/or before the first index $i$ syllable $\ap'$ is
considered to be between $\ap$ and $\ap'$.

For the following considerations it is (not necessary but) helpful
to visualize $\bt$ by the braid scheme explained in \cite{bind}.

Separate $\bt$ in
\eqref{wd} into subwords $\ap_1\dots \ap_n$, such that $\ap_i$
contains only syllables of odd or even index, and this parity
changes between $\ap_i$ and $\ap_{i+1}$. Then for a syllable
$\sg_k^l$ occurring in $\ap_i$, put the integer $l$ at the point
$(k,i)\in \bN\times \bN\subset \bR^2$ in the plane. Here $(k,i)$
is the point in the $i$-th row and $k$-th column, with rows
numbered (as in Cartesian coordinates) from bottom to top and
columns from left to right.

One obtains a certain checkerboard pattern of integers we call
\em{braid scheme} of $\bt$. (If we do not put any integer on
a point $(k,i)$, we assume its ``content'' is zero, or it is
``empty''. So for all non-empty points $(k,i)$ in the scheme,
$i+k$ is always even or always odd.)

One can \em{reduce} the
scheme by moving an integer $l$ at $(k,i)$ to $(k,i-2)$ if $i>2$
and the points $(k\pm 1,i-1)$ are empty. We call the scheme
reduced if it does not admit any such move. Then in a reduced
scheme, summit syllables of $\bt$ are those, whose entries in
the scheme are ``on top'' when viewing the scheme from the left.
{}From this viewangle the following ``geographic'' choice of
terminology becomes more plausible.

\begin{defi}
Summit syllables still have a cyclic order from \eqref{wd}. We call
the subword $\bt'$ of $\bt$ in \eqref{wd} made of summit syllables
the \em{maximal subword}. The subword made of non-summit syllables
(i.e. the subword obtained by deleting in $\bt$ all syllables in
$\bt'$) is called \em{non-maximal subword}. 

Note that neighbored summit syllables have indices $k_i$ differing
by $\pm 1$. We say that a summit syllable is \em{minimal} resp.
\em{maximal} if its both neighbors have higher resp.\ lower index.

We call $\bt$ \em{summit reduced} if all its minimal summit syllables
are non-trivial.  We call $\bt$ \em{index reduced} if it is
non-singular and its index sum $\sum_{i=1}^n k_i\cdot l_i$ cannot be
reduced by a Yang-Baxter relation, i.e. $\bt$ contains no
$\sg_{i+1}\sg_i\sg_{i+1}$ as subword.
\end{defi}

\begin{lemma}
Index reduced $\So$ summit reduced. In particular a summit reduced
form always exists.\qed
\end{lemma}

Recall that a positive resolution tree is a rooted tree with directed 
edges,
whose vertices (nodes) contain positive braid words, the root labelled
by $\bt$. Every vertex has exactly one incoming edge, except the
root that has none, and zero or two outoutgoing edges. In former
case it is labelled by an unlink (terminal node).
In latter case it is labelled by
a word of the form $\ap\sg_i^2\ap'$, with $\ap,\ap'$ positive
words, and the two vertices connected by the outgoing edges 
are labelled by $\ap\sg_i\ap'$ and $\ap\ap'$, or positive words
obtained therefrom by Markov equivalence (isotopy of the closure link).

In \cite{Nakamura} the following fact was observed, and used
decisively, and we shall do the same here.

\begin{theorem}(Nakamura \cite{Nakamura})
$MWF(\bt)$ is the maximal number of components of a (link in a) node
in a positive resolution tree for $\bt$.
\end{theorem}

In particular, $MWF$ is monotonous (does not decrease) under word
extension, and does not depend on the exponent of non-trivial
syllables.

\begin{lemma}\label{lemh}
If $\bt$ is \em{summit reduced}, and $\bt'$ is obtained from $\bt$
by removing all summit syllables, then there is a positive resolution
tree for $\bt$ that contains $\bt'$ as a node.
\end{lemma}

\proof Since minimal syllables are non-trivial, one can delete them
in the resolution tree. The two neighbors in the maximal subword
join to a non-trivial new minimal syllable, and so one iterates the
procedure. \qed

Since all $\sg_{n-1}^{l'}$, $\sg_{n-2}^{l'}$ in \eqref{wd} occur as
summit syllables, $\bt'$ has split last two strands, and so we have
a quick proof of Nakamura's main result.

\begin{corr}(Nakamura \cite{Nakamura})
Any summit reduced positive word on $n\ge 3$ strands has $MWF\ge 3$.
\qed
\end{corr}

In particular any positive braid representation of a $(2,n)$-torus
link can be reduced to the standard one by index-decreasing YB
relations and removals of nugatory crossings.  

\subsubsection{The proof of Theorem \reference{thmwf}:
Initial simplifications\label{4.0}}

The following fact is well-known:

\begin{theorem}\label{th4.1}
$MWF(\bt)=1$ if an only if $[\bt]_i=1$ for $i=1,\dots,n-1$.
\end{theorem}

Now for $MWF=3$ it suffices to ensure that (either we can reduce the
braids) or can find words, whose non-maximal subwords do not
give an unknot. \em{For the rest of the section we assume that
$\bt'$ gives the unknot.}

We will work by induction on the number of strands, and for
fixed number of strands on the index sum. So we consider a
positive braid word $\bt$, and assume w.l.o.g. it has the
smallest index sum among positive braid representatives of
its closure link for the same number of strands. For such
$\bt$, we will either reduce it (by at least one strand or
crossing), or show $MWF\ge 4$.

Most braids $\bt$ will be easily dealt with, but there remain certain
families of words that require a case-by-case study. We decided not
to omit too many of the (tedious) details of this part, in order
to keep the proof followable, even if it may not contribute to
its (esthetic) appearance.

Note that in order to prove $MWF(\bt)\ge 4$, it suffices to
go over to a (link in a) suitably chosen node in a positive
resolution tree for $\bt$ and show $MWF\ge 4$ for this node. In
particular, we can remove from $\bt$ all syllables of index $\le
k$ and the $k$ resulting left isolated strands.

The case of reductions is more delicate. In some situations
we can describe them directly, but this is not always the case.
Then we proceeded as follows. First we took generic examples, in
making all syllables non-trivial whose triviality we have not argued
about. We adjust parities so that the closure is a knot, and checked
using KnotScape \cite{KnotScape} that the braid reduces (by at least one
strand/crossing).

Later we wrote a computer program that seeks reductions by keeping
given crossings rigid. Such reductions would commute with replacing
rigid crossings by any tangle, in particular by any non-trivial
braid word syllable. (Non-trivial syllables behave
similarly to  rigid vertices, and suggest that the reduction
is likely to work in general.) This way we can find reductions for
infinite families of braids on a given number of strands. By turning
all summit index-1 syllables into rigid crossings, one can also
handle the braids that occur for an increasing number of strands.

The technical details of the application of this program are,
however, tedious and little insightful. Instead we content
ourselves with giving the examples we processed with KnotScape.
% The alternative proof of theorem \reference{thg} given later
% is another reason why we felt not motivated to make this
% way the approach here even more technical.

\begin{defi}
A \em{valley} resp. \em{mountain} is a subword of the maximal word
that starts with same index syllables and contains only one minimal
resp. maximal syllable. This syllable is called the \em{bottom} of
the valley resp.\ \em{summit} or \em{top} of the mountain.
The index of the bottom/top is the \em{depth} resp.\ \em{height}.
\end{defi}

We assume there are at least two mountains of maximal height (i.e.
$n-1$). Otherwise we have a split component or a $(2,n)$-torus
connected component, or a reducible braid and can work by induction
on the number of strands. Similarly at least one valley has
depth $1$, otherwise $\sg_1$ in $\bt'$ remains reducible in $\bt$.

The following operation will be somewhat important, and we will
call it ``filling the valley''.

\begin{lemma}(``filling the valley'')
Any valley can be removed from the maximal subword in the positive
resolution tree.
\end{lemma}

\proof Same as for lemma \reference{lemh}. \qed

\begin{lemma}
If a mountain $M$ is not of maximal height
(i.e. $n-1$), then $MWF(\bt)\ge 4$.
\end{lemma}

\proof Let $k<n-1$ be the height of $M$.
We assumed there are at least two mountains of maximal height. So
now w.l.o.g. assume some, say
the left, of the neighbored mountains of $M$ has height $k'>k$.
Fill the two valleys around $M$ starting with $\sg_{k-1}$.
Then the maximal subword has a syllable index sequence
$k+1,k,k-1,k,k-1,k$. Make the second and third syllable trivial
(if not already), and apply YB relations, moving the fourth
syllable to the left: $k+1,\ul{k-1},k,k-1(,k-1),k$. The result is a
summit reduced word, in which a new syllable of index $k-1$
(the underlined one) was removed and it became non-maximal. Hence the
non-maximal subword has exponent sum $[\bt']_{k-1}>1$, and so
$MWF\ge 4$. \qed

\begin{lemma}\label{lem2val}
If $\bt$ has $>2$ valleys of depth at most $n-3$, then $MWF(\bt)\ge 4$.
\end{lemma}

\proof % It suffices to check it for 4-braids.
It suffices to check for 3 mountains (as one can fill
separate valleys) and 4-braids (as one can fill valleys by
levels as in the proof of lemma \reference{lemh} and the remark
after theorem \reference{th4.1}). This is just
the word [1232112321121321], which is easily checked (to
have $MWF=4$). \qed

\subsubsection{Two mountains\label{S4.1}}

We assume in \S\reference{S4.1} and \S\reference{S4.2} that $n\ge 5$.
The case of $4$-braids is considered later in \S\reference{S4.3}.
We refer to \S\reference{brs} for the use of notation we will employ.

We assume first $\bt$ has two mountains. By the previous remarks they
are both of height $n-1$.

So now consider words with syllable index sequence
\begin{eqnarray}\nonumber
& 1,2,\dots,n-2,n-1,\quad
p_1@,\dots,p_k@,\quad
n-2,n-3,\dots,g+1,g!,\quad \\
\label{xx} &
g+1,\dots,n-2,n-1,\quad
q_1@,\dots,q_l@,\quad
n-2,n-3,\dots,2,1
\end{eqnarray}
such that $l+k=n-3$ and
$\{p_1,\dots,p_k,q_1,\dots,q_l\}=\{1,\dots,n-3\}$.
We will distinguish only between non-trivial and trivial syllables
(in former case exponent is immaterial). For non-trivial syllables
we write an exclamation mark
after the index, for trivial ones an `at' (@) sign.
If none of $!$ and @ is specified, we do not exclude explicitly any of
either types. We write $\bt_1,\dots,\bt_6$ for the subwords separated by
space in \eqref{xx}.

Assume w.l.o.g. (up to reversing the braid's orientation) that
some $p_i$ is $1$, and let up to commutativity the $p_1,\dots,p_k$
subword be written as $h,h-1,\dots,1,p_1',\dots,p_{k-h}'$ (the 
$p_i'$ contain the indices above $h$ occurring as $p_i$).

Assume the maximal syllable of $\sg_{h+1}$ in $\bt_1$ is non-trivial.
Then one can write $\bt$ as word with a subword of index sequence 
\begin{eqn}\label{siks}
1,\dots,h+1!,h@,\dots,1@\,. 
\end{eqn}
If now $h<n-3$, we can make in \eqref{siks} all syllables
trivial except $h+1$, which we make of exponent $2$, then
split a loop (component of a link in a node of a positive
resolution tree) by removing all the syllables in \eqref{siks}
together with the terminating `$1$' in \eqref{xx}. We obtain
\begin{eqnarray*}
& h+2,\dots,n-2,n-1,\quad
p_1'@,\dots,p_{k-h}'@,\quad
n-2,n-3,\dots,g+1,g!,\quad \\
& g+1,\dots,n-2,n-1,\quad
q_1@,\dots,q_l@,\quad
n-2,n-3,\dots,2.
\end{eqnarray*}
(Here no syllables of index 1 occur, and the split loop is the
isolated leftmost strand.) Then we fill the valley starting with the
index-$n-1$-summits, splitting another loop (the rightmost strand),
\[
h+2,\dots,n-3,n-2,\quad p_1'@,\dots,p_{k-h}'@,\quad q_1@,\dots,q_l@,
\quad n-2,n-3,\dots,2,
\]
and are left with a
word that has at least two $\sg_{n-2}$. So $MWF\ge 4$. 

\begin{caselist}
\case Assume the smaller valley has depth $g>1$.
The $\bt$ has two index-1 syllables, a trivial (non-summit)
and a non-trivial (summit) one. By a flype one can exchange
them, and so have a non-summit reduced word. Then one
can change $\bt$ to a word of smaller index sum, and so
we are done by induction.

\case Now $g=1$. We write $\bt_1,\dots, \bt_6$ for the
6 subwords separated by spacing in \eqref{xx}.
By a similar argument as after \eqref{siks} we can argue that
if one can reorder the syllables in $\bt_{2,5}$ so that
$\bt$ has a subword with an index sequence 
\begin{eqn}\label{(eight)}
k,k+1,\dots,h-1,h!,h-1,\dots,k\,,
\end{eqn}
with the first or last
$h-k$ syllables belonging to $\bt_{2,5}$ and the others
to $\bt_{1,3,4,6}$, and $h<n-2$, then $MWF\ge 4$.

W.l.o.g. assume $1\in \bt_5$ (which is meant to abbreviate that
$\bt_5$ contains an index-$1$ syllable).

\begin{caselist}
\case
Now if $2\not\in\bt_5$ (so $2\in\bt_2$) then the
$2$-index syllables in $\bt_{4,6}$ are trivial
(because we have otherwise \eqref{(eight)} with $k=1$, $h=2$).

We distinguish several cases by the subwords of
(non-summit) syllables of index $1,2$ and $3$ in $\bt_{2,5}$.
We separate the subwords between $\bt_{2}$ and $\bt_5$
by a vertical line `$|$'. (We assume here that $n\ge 6$.
The case $n=5$ must be handled by a separate, but
simplified, argument.) Note that by symmetries we can
exchange the words of $1,2$ and $3$ left and right from `$|$'
and also (simultaneously) reverse both, and can also use the
commutativity of $1$ and $3$. Then we are left with
the following cases.

\begin{caselist}
\case
$2|31$. If $3\in \bt_5$, then both $3$ in $\bt_{1,3}$ are trivial.
(Otherwise, we would have \eqref{(eight)} for $k=2$ and $h=3$.)
% \begin{verbatim}
Below we give pairs of words, the first obtained by extending all
admissible syllables to be non-trivial, and the second one by
extending the first word to one with knot closure, which was then
checked to reduce (by at least one strand, not necessarily to a
3 braid).
\\
$[1122344  2 55665544322111   2334455 3441 66554433211]$, \\
$[11222344 2 5556655444322111 2334455 3441 666554433211]$ reduces. \\
% \end{verbatim}
\case
If $3\in \bt_2$, then one of both $3$-index syllables in
$\bt_{1,3}$ must be trivial. (It is the syllable in $\bt_3$ if `$3$'
occurs before `$2$' in $\bt_2$, or the syllable in $\bt_1$
otherwise.)

\begin{caselist}
\case
$ 23|1$.
% \begin{verbatim}
\\
$[1122344 23 556655443322111    2334455 441 66554433211]$, \\
$[1122344 23 556665554443322111 2334455 441 66554433211]$ reduces. \\
% \end{verbatim}
\case $ 32|1$
% \begin{verbatim}
\\
$[11223344  32 55665544322111   2334455 441 66554433211]$, \\
$[112233444 32 5556665544322111 2334455 441 66554433211]$ reduces.  \\
% \end{verbatim}
\end{caselist}
\end{caselist}

\case
Now assume $2\in\bt_5$. Since $1\in\bt_5$, now the exclusion
of \eqref{(eight)} shows that only one of the $2$-index syllables
in $\bt_{4,6}$ is trivial, but we will show that also one in
$\bt_{1,3}$ is. Namely, by making the proper index-$1$ (summit)
syllable to exponent 2, and one of the $2$-index syllables of $\bt_
{1,3}$ trivial, one can slide by braid relations the $1$-index
syllable $\sg_1=X$ from $\bt_5$ to $\bt_2$.
By applying the previous argument, both
$2$-index syllables in (the now modified) $\bt_{1,3}$ are trivial.
One of them was previously made trivial to slide $X$ in, but the
condition on the other one persists for the original braid.

In all situations,
make all the other syllables in $\bt_{1,3,4,6}$ non-trivial and
check using KnotScape that the braid reduces. 

% MAKE THIS EXACT!!!

\begin{caselist}
\case $ 3| 12$. In this case the above argument shows that the
syllable $2\in\bt_1$ is trivial.
% \begin{verbatim}
\\
$[11233445566  344 55443322111 234455   12 66554432211]$, \\
$[112334455666 344 55443322111 23444555 12 665544322111]$ reduces. \\
% \end{verbatim}
\case $ 3| 21$. Here $2\in\bt_3$ is trivial.
% \begin{verbatim}
\\
$[112233445566 344 5544332111 2234455 21 6655443211]$, \\
$[112233445566 344 5544332111 2234455 21 66655544432111]$ reduces. \\
% \end{verbatim}
\case $| 123$. As before $2\in\bt_1$ is trivial.
% \begin{verbatim}
\\
$[11233445566 55443322111 234455 123 665544332211]$, \\
$[11233445566 55443322111 234455 123 666555444332211]$ reduces. \\
% \end{verbatim}
\case $| 132$.
% \begin{verbatim}
\\
$[11233445566 55443322111 2334455 132 66554432211]$, \\
$[11233445566 55443322111 2334455 132 66655544432211]$ reduces. \\
% \end{verbatim}
\end{caselist}
\end{caselist}

\end{caselist}

\subsubsection{More than two mountains\label{S4.2}}

% VERY VAGUE AND MUCH NEEDS TO BE DONE USING COMPUTATIONS !!!
To deal with the general case, 
now we make the following modifications. We call a
summit syllable sequence with indices $n-2$ and $n-1$
terminated on both sides by $n-1$'s a modified
mountain or \em{plateau}. We have again by lemma
\reference{lem2val} only two valleys of depth $<n-2$,
or alternatively only two plateaus (now instead of mountains).
The case of more than two mountains is thus mainly a
adaptation of the case of two mountains, replacing
mountains by plateaus.

% Also we allow $\bt$ and the maximal summit subword
% to throw in occurrance of non-trivial $\sg_{n-1}$ syllables.

% The words that admit reductions apply as in the previous
% two-mountains-case,
% 
% 2 deep, 1 shallow valley
% 123211232213211 reduces
% 123121123223211 reduces 
% 
% CHECK CAREFULLY !!!
% 
% , but for proving $MWF\ge 4$ the
% following modifications of the argument are needed.

Again we may assume non-maximal subwords have exactly one (and trivial)
syllable per index. The elimination of the maximal subwords can be done
similarly. % only that when arriving at $\dots,n,n-1@,n,n-1@,n$

We distinguish two cases as in the above study of the 2-mountain
words, depending on the depth $g$ of the second valley (the
other valley has depth $1$ by the same argument as above).

\begin{caselist}
\case $g>1$. We use the previous flyping argument.
% In the first case above the same arguments show $g=2$ and (after
% deleting the index-$n-2$ syllables in the plateaus) that
% the syllables in $\bt_{3,5}$ in the form \eqref{seven} are trivial,
% except the index-$n-2$ syllables in both and the index-$2$ and
% $n-1$ syllable in $\bt_{3}$.
% (Otherwise, fill one valley, remove a loop from \eqref{hh}, and
% something non-trivial remains.) Also $n-2\in\bt_1$ may now be trivial.
% 
% {$112233\hat 4554455\ \ 321\ 5544\ul 32\ 2233445544\ 5544
% \ul 3\ul 2$}
% 
% % CHECK ANYTH REDUCES !!!
% One has to check that such braids reduce. For example, \\
% % \begin{verbatim}
% 11223345544455 321 55544322 22333445544 5544432 reduces, and \\
% 1122334445544455 321 55544322 22333445544 5544432 reduces. 
% \end{verbatim}

\case $g=1$.
In the second case we had restrictions on exponents of syllables
with index $2$ and $3$ occurring in the maximal subwords from
the position of syllables with index $1$ and $2$ occurring in
the non-maximal subwords. The restrictions on $2$-(index) syllables
from $1$-(index) syllables remain. So do the restrictions on
$3$-syllables from $2$-syllables unless we have $\le 5$
strands. The argument is the same: one can still pull
out two loops and has at least two letters of $\sg_{n-2}$.
% (additionally discarding all additional non-trivial $\sg_{n-1}$
% syllables)

Now we check restrictions on $3$-index syllables for $5$ strands
and reducibility. We have up to extensions a finite number of
special braids to verify. Clearly, extensions are never
admissible for non-summit syllables, and always admissible for
summit index-1-syllables (since they are all non-trivial). 
Also extensions of subwords $(n-2\ n-2\ n-1)^k$ occurring
repeatedly are redundant, since doubling the letter $n-1$ is the
same as deleting the $n-2\ n-2$ for the next $k$. We
will verify that the property $MWF\ge 4$ resp. reducibility
does not depend on the value of $k$ as soon as $k>0$.
We will thus deal only with the other extensions.

\begin{caselist}
\case The non-summit $2$-index syllable is between two valleys
of depth $n-2=3$ - this is handled as before: \\
$[1234334123343211234321]$ has $MWF=4$, \\
$[1234334213343211234321]$ has $MWF=4$, \\
$[1234334233432112314321]$ has $MWF=4$.

\case The non-summit $2$-index syllable is between
one depth $n-2=3$ valley and one depth $1$ valley. These are words
of the form $[123412(334)^k3211234(334)^l321]$, \\
$[123421(334)^k3211234(334)^l321]$,\\
$[12342(334)^k32112341(334)^l321]$ \\
for $k>0$, $l\ge 0$, and their extensions.

\begin{caselist}
\case $[123412(334)^k3211234(334)^l321]$ and extensions.

\begin{caselist}
\case $l=0$. To display the extendability of syllables,
in the following notation the necessarily trivial syllables are
hatted, while a possibly trivial syllable is underlined.

$[1\hat 2\hat 3\ul 412(334)^k\ul 3\ul 211\hat 2\hat 3\hat 4\ul 3\ul 21]$. 

When $l=0$, then making non-trivial any single of the hatted syllables
makes $MWF=4$ already for $k=1$, while for any $k>0$ making any
combination (possibly all) of the underlined syllables non-trivial
gives $MWF=3$.

Now again check that braids reduce, for example: \\ 
$[1234123343343211234321]$ has $MWF=3$, \\
$[12341233433433432111234321]$ reduces to 3 strands, \\
$[1234441233433433221123433221]$ reduces to 3 strands.

\case 
When $l=1$, the already for $k=1$ we have without extensions
$[1234123343211234334321]$ and $MWF=4$.
\end{caselist}

\case $[123421(334)^k3211234(334)^l321]$ and extensions.

\begin{caselist}
\case When $l=0$, we have for all $k$ the extendability

$[1\ul 2\hat 3\ul 421(334)^k\ul 3\hat 211\ul 2\hat 3\hat 4\ul 3\hat 21]$,

with the same explanation as before.

For example, $[122344213343343343321122343321]$ reduces to 3 strands.

\case When $l\ge 1$, then already for $k=1$ and no extensions
$MWF([1234213343211234334321])=4$.

\end{caselist}

\case $[12342(334)^k32112341(334)^l321]$ and extensions.

\begin{caselist}
\case When $l=0$, we have

$[1\ul 2\hat 3\ul 42(334)^k\ul 3\ul 211\hat 2\hat 3\hat 41\ul 3\hat 21]$.

For example, $[122344233433433433221123413321]$ reduces to 3 strands.

\case When $l>0$, already for $k=1$, $l=1$ without extension
we have $[1234233432112341334321]$ and $MWF=4$.
\end{caselist}
\end{caselist}

\case The non-summit $2$-index syllable is between
two depth $1$ valleys. These are (up to symmetry) words of the form\\
$[1234(334)^k321123421321]$ and \\
$[12341(334)^k32112342321]$, \\
for $k>0$, and their extensions (here necessarily $l=0$).

\begin{caselist}
\case $[1234(334)^k32112342132]$ and extensions:

$[1\ul 2\ul 3\ul 4(334)^k\ul 3\hat 211\ul 2\hat 3\hat 421\hat 3\hat 21]$.

Again for $k=1$ making non-trivial all underlined syllables 
gives $MWF=3$, while making non-trivial any of the hatted
syllables gives $MWF=4$. 

For example,
$[12 3  4 3343       2112 3 214321]$ has $MWF=3$, and its extensions\\
$[122333443343       211223 214321]$ and \\
$[12233 4433433433433211223421 321]$ \\
were checked to reduce to 3 strands.

The same is the outcome for $k>1$.

\case $[12341(334)^k32112342321]$ and extensions:

$[1\h 2\u{34}1(334)^k\u 3\h 211\u 2\h 3\h 42\h 3\u 21]$

The reducibility cases follow analogously. % CHECK !!
For example, $[123344133433433433211223423221]$ reduces to 3 strands.

% For $l>0$ we have always $MWF=4$.

\end{caselist}

\end{caselist}
\end{caselist}

\subsubsection{4-braids\label{S4.3}}

If some mountain is not of height 3, or $>2$ valleys of
depth $1$ exist, then we are done as before (see \S\reference{4.0}).

So the maximal subword is of the form
$1,2,3,(2!,3)^p,2,1,1,2,3,(2!,3)^n,2,1$, and the non-maximal
subword is a single $\sg_1^1$. We separate the summit syllables
and their letters by the summit syllables of index $1$ into
a \em{left} and \em{right plateau}. Assume w.l.o.g. the (non-summit)
$\sg_1^1$ is in the left plateau. The word `in' is to mean that
in cyclic order of the syllables of $\bt$ the syllable $\sg_1^1$
can be written to occur just before or after a syllable with
index $n-1$ that belongs to the left plateau. This means that
we can write $\bt$ as
\begin{eqn}\label{blmn}
\bt_{l,m,n}=[1\vtbox{2\\.}3(223)^n1(223)^m\vtbox{2\\..} 1123(223)^l21]
\end{eqn}
with $n,m,l\ge 0$, or some of its extensions. Note that MWF will
be monotonous in $m,n,l$, i.e. $MWF(\bt_{l,m,n+1})\ge
MWF(\bt_{l,m,n})$ etc. Using symmetry assume $n\ge m$.

One can check already at this stage that such words are not
$B$-adequate. So we obtain corollary \reference{cort}, and
for the proof of theorem \reference{thg} the rest of the
argument here can be replaced by the application of theorem
\reference{thq}. Note that $B$-adequacy is invariant
under isotopy preserving writhe and crossing number, so
any other positive 4-braid word giving the same link is
not $B$-adequate either.

If $n+l>0$, then the $2$-index syllable $\vtbox{2\\..}$
in \eqref{blmn} must be trivial.
% there are $2*$ of the left plateau on the left (right) of
% $\sg_1^1$ or some in the right plateau, then the right
% (left) $2$-index syllable of that mountain must be trivial.
Otherwise remove all $2!$ in $(2!,3)^m$ (if any), and split two loops
as explained after \eqref{siks}. The `$2!$' in $(2!,3)^l$ or
$(2!,3)^n$ remain, and so $MWF=4$.
With a similar argument we see that if $m+l>0$, then the $2$-index
syllable $\vtbox{2\\.}$ is trivial.

We distinguish 3 cases depending on whether these arguments
apply or not.

\begin{caselist}
\case Both non-triviality arguments apply. So we have a family of
words $\bt_{l,m,n}=[123(223)^n1(223)^m21123(223)^l21]$ with $n,m>0$,
or $l>0$ and their extensions, and both $\vtbox{2\\..}$ and
$\vtbox{2\\.}$ are trivial.

\begin{caselist}
\case $l=0$. We assumed $n,m>0$, and already for $n=m=1$,
the word $[12322312232112321]$, we have $MWF=4$.

\case $l>0$.

\begin{caselist}
\case $m=n=0$. These are extensions of $[1231211123(223)^l21]$,
and the admissibility is found to be:

$[1\h 2\h 31\h 211\ul{23}(223)^l\ul 21]$.

The two letters `$2$' in the left plateau cannot be doubled
($[123122112322321]$, $[122312112322321]$), neither the `$3$' \\
($[123132112322321]$), since $MWF=4$ already for $l=1$.

Without extension, $l=1$ ($[123121112322321]$) and $l=3$
($[123121112322322322321]$) reduce.
% 
% $[1231211223223223221]$ reduces.
% 
% % So so do (sloppy) $123121123[23]+21$ in regexp notation.
% 
% $[123121112223223223223221]$ has $MWF=3$. \\
% $[1231211122232232232232221]$ reduces. \\
So we find that $[12312112*3[23]+2*1]$ reduce.
(Recall that, while `$2*$' in
a braid word should mean at least one letter `$2$', the term `$[23]+$'
should mean a possibly empty sequence of letters `$2$' and `$3$'.
We distinguish braid words from index sequences by not putting
commas between the numbers.)

\case $m+n>0$; this reduces to the case of $m=n=0$ with some of the
twos or the three in the left plateau doubled, where we found $MWF=4$.

\end{caselist}

\end{caselist}

\case In the case one of the non-triviality conditions on
$\vtbox{2\\..}$ and $\vtbox{2\\.}$ does not apply,
we have $123(223)^n12112321$ with $n>0$ and its extensions.
(Now the right plateau is a mountain.)

% \begin{caselist}
% \case
The case $n=3$ ($[123223223223121112321]$) simplifies. \\
% So so do (sloppy) $123[23]+312112321$ in regexp notation.
% Also 12223223223223121112321, so so do (sloppy) $12*[23]+312112321$.
$[122332232232232231211222321]$ simplifies.

However, $MWF([123223122112321])=4$, so the right $2$-index syllable
in the left plateau ($\vtbox{2\\..}$ in \eqref{blmn}) must be
trivial.

The right 2 and the 3 in the right plateau must be trivial: \\
$[1232231211123221]$, \\ $[1232231211123321]$ have $MWF=4$.

But the left `$2$' of the left and right plateau may not be trivial:
$[12223223223223121112222321]$ has $MWF=3$. \\
$[122232232232231211122222321]$ simplifies. \\
So $[12*[23]+312112*321]$ simplifies. % ?

We arrive at the form $[1\u{23}(223)^n1\h 211\u 2\h 3\h 21]$.

% \dots
% \end{caselist}

\case
In case both non-triviality conditions do not apply, we have
$[12312112321]$ and its extensions. (So both plateaus are
mountains.)

Since non-triviality is nonetheless possible, we may have non-trivial
$2$-index syllables in the left plateau.
We distinguish three cases again according to whether
$\vtbox{2\\..}$ and $\vtbox{2\\.}$ are trivial or not.

\begin{caselist}
\case Both $2$-index syllables are trivial: $[12312112321]$.

% 12312112321

By direct check: \\
% 122312112321 $MWF=3$
$[123312112321]$ has $MWF=3$, \\
% 123122112321 $MWF=3$
$[123121122321]$ has $MWF=3$, \\
$[123121123321]$ has $MWF=3$, \\
$[123121123221]$ has $MWF=3$. \\

% 12233122112321 $MWF=3$
% 122331221122321 $MWF=4$
% 122331221123321 $MWF=4$
% 122331221123221 $MWF=4$

$[12312112233221]$ has $MWF=3$, \\
% 122312112233221 $MWF=4$
$[123312112233221]$ has $MWF=4$. \\
% 123122112233221 $MWF=4$

So one can extend the right mountain's `$2$'s and one of the
left or right mountain's `$3$'s, but not both `$3$'s.

$[123121122233221]$ reduces, \\
$[123312112223221]$ reduces, \\
$[1233121122233221]$ has $MWF=4$.

So a reducing check is to be made on $[11*231211*22*3*22*]$.
% and 11*22*3*122*11*232

$[123121112233221]$ reduces, \\
$[123121112222333322221]$ reduces, \\
$[1233312111222322221]$ reduces.

% \begin{caselist}
% \case
% 11*231211*22*3*22*
% 
% \case 11*22*3*122*11*232
% 
% 122331222112321 DOES NOT REDUCE !!!
% 122331222112321 connected-2-cabled MWF =7 !!!!!
% 
% So we have exceptions 11*22*3*122*11*232 !
% (Or 22*3*122*11*23211*, which is the
% first family in theorem \reference{thmwf}.)
% 
% 
% 12331222112321 2-cabled MWF =6 !!!!!
% 112331222112321 reduces
% 11233331222112321 reduces
% 
% 12231222112321 2-cabled MWF =7 !!!!!
% 1223122112321 2-cabled MWF =7 !!!!!
% 
% So we have exceptions 11*22*3122*11*232 !
% (Or 22*3122*11*23211*, which is a special case of the
% first and second family in theorem \reference{thmwf}.)
% 
% 1223312112321 2-cabled MWF =6 !!!!!
% 12222333312112321 reduces

\case
One $2$-index syllable is non-trivial. This is the word
$[123122112321]$, with $MWF=3$ (the case
$[122312112321]$ is symmetric).

We have the following extensions:

$[1233122112321]$ has $MWF=3$ (and reduces), \\
$[1231221122321]$ has $MWF=4$, \\
$[1231221123221]$ has $MWF=3$ (and reduces), \\
$[1231221123321]$ has $MWF=3$ (and reduces).

The following combined extensions are to check:
$[12331221123321]$ has $MWF=4$, \\
$[12331221122321]$ has $MWF=4$, \\
$[12312211223321]$ has $MWF=4$. \\

Thus we are left to deal with
$[1233*122*112321]$,\\
$[123122*112322*1]$,\\
$[123122*11233*21]$,\\
and check that they all reduce:

$[12333312222112321]$ reduces,\\
$[12312222112322221]$ reduces,\\
$[12312222112333321]$ reduces.

\case Both $2$-index syllables are non-trivial:
$[1223122112321]$. This is a braid word $\bt_0$ for $13_{9465}$.
We have the following extensions:

$[12233122112321]$ has $MWF=3$, \\
$[12231221122321]$ has $MWF=4$, \\
$[12231221123221]$ has $MWF=4$, \\
$[12231221123321]$ has $MWF=3$.

The only common extension of the two $MWF=3$ extensions is \\
$[122331221123321]$, which has $MWF=4$.

So it remains to verify that (the closures of)
$[22*3*122*11*23211*]$ and $[22*3122*11*23*211*]$ have braid index 4.

% 1222233312211232111 no reduce
% 12223122211123321 no reduce

For this we use the two-cabled MWF inequality. Let for a braid
$\bt\in B_n$, the ``two-cabled'' braid $(\bt)_2\in B_{2n}$ be
obtained from $\bt$ by replacing in \eqref{wd} each $\sg_i$
by $\sg_{2i}\sg_{2i-1}\sg_{2i+1}\sg_{2i}$. Then $(\bt)_2$
is a braid representation of the (blackboard framed) two-cable
link $L_2$ of the closure $L=\hat\bt$ of $\bt$. We consider now
for $\bt$ the above braid $\bt_0$.

We know, from the computations described in \cite{MorSho,WilFr},
that $MWF((13_{9465})_2)=7$. This is in fact true also for the
connected cable (the one with braid representation $(\bt_0)_2
\cdot\sg_1$). Now we claim that reducing or resolving a clasp
(changing a $\sg_i^2$ into a $\sg_i$ or deleting it) does not
reduce the two-cabled MWF bound. The 2-cable of a clasp can be
resolved by resolving 4 clasps. The 2-cable of a crossing in a
2-cabled clasp can be resolved by resolving one clasp and changing
twice $\sg_i^2\to \sg_i$. Finally, the 2-cable of an isolated 
$\sg_i$ can be reduced into two internal twists of the
doubled original strand. Such twists can be collected for
every doubled component, resolved for each doubled component to
one, and joined if doubled components are joined by reducing
a doubled crossing in a doubled clasp. So the two-cabled MWF
reduces to the one of the connected cable of $13_{9465}$ and
we are done.
\end{caselist}
\end{caselist}

The proof of theorem \reference{thmwf} is now completed.

\subsection{Positive links\label{Spos}}

In this section, we will refine the arguments proving
theorem \reference{thqp} to restrict the possible 
3-braid representations of positive links. 
% and then, % in an alternative way, classify
% positive braid links of braid index 3.
% 
Our positivity considerations will make use of the criterion
of Yokota \cite{Yokota}, and the Kauffman polynomial $F$.
We recall the properties \eqref{wseven}~--~\eqref{wseven.7}\,
that determine $F$ and its writhe-unnormalized version $\Lm$.

Note that for $P$ one can similarly define a regular isotopy
invariant
\begin{eqn}\label{dPtl}
\tl P(D)(a,z)=(ia)^{-w(D)}P(D)(ia,iz)\,,
\end{eqn}
with $i=\sqrt{-1}$. Then $\tl P$ satisfies similar relations to
\eqref{wseven}~-- \eqref{wseven.7}. The difference to $\Lm$ is that
$\tl P$ is defined on oriented link diagrams, and that the term making
orientation incompatible on the right of \eqref{wseven} is missing.

\begin{theorem}\label{thF}(Yokota \cite{Yokota})
If $L$ is a positive link, then 
\[
\md_a F(L)\,=\,\md_v P(L)\,=\,1-\chi(L)\,,
\]
and
\[
[F(L)]_{a^{1-\chi(L)}}=[P(L)(ia,iz)]_{a^{1-\chi(L)}}\,.
\]
\end{theorem}

We also require an extension of braids to the context of $F$.
This was described in \cite{BW} and \cite{Murakami}, but we use
only the generators of the algebra defined there. Strings will be
assumed numbered from left to right and words will be composed
from bottom to top. We write $\sg_i$ for a braid generator, where
strand $i$ from the lower left corner, passing over strand $i+1$,
goes to the upper right corner.

We add elements $\dl_i$ of the following form:
\[
\diag{4mm}{6}{2.5}{
  \pictranslate{0 0.5}{
  \picline{0 0}{0 2}
  \picmultigraphics{3}{0.2 0}{
    \picfilledcircle{0.55 1}{0.03}{}
  }
  \picline{6 0}{6 2}
  \picmultigraphics{3}{0.2 0}{
    \picfilledcircle{5.05 1}{0.03}{}
  }
  \picline{1.5 0}{1.5 2}
  \picline{4.5 0}{4.5 2}
  \piccirclearc{3 0}{0.5}{0 180}
  \piccirclearc{3 2}{0.5}{180 0}
  \picputtext{2.5 -.3}{\footnotesize $i$}
  \picputtext{3.5 -.3}{\footnotesize $i+1$}
  }
}\es.
\]
By hat we denote the usual closure operation.

A word in the described generators gives rise to an unoriented tangle
diagram that turns into an unoriented link diagram under closure. 
If the word has no $\dl_i$ then this diagram can be oriented to
give a(n oriented) closed braid diagram. We will assume this
orientation is chosen. Otherwise, a coherent orientation is not
generally possible. In particular, the sign of exponents of $\sg_i$
in this context may not coincide with the sign of the corresponding
crossings after some (or even any) orientation choice of the diagram.
For kinks (the diagram fragments occurring in \eqref{wseven.5} on the
left hand-sides),
however, a sign is definable since any possible component
orientation chosen gives rise to the same (skein) sign. So we will be
able (and we will need) to distinguish between positive (in the
left equation of \eqref{wseven.5}) and negative
kinks (in the right one).

\begin{lemma}\label{lemma1}
Let 
\begin{eqn}\label{Df}
D=\wh{\es}\,\left(\,
\sg_1^{k_{1,1}}\sg_2^{k_{1,2}}\dots \sg_1^{k_{1,n_1}}\sg_2^{-1}\quad
\sg_1^{k_{2,1}}\dots \sg_1^{k_{2,n_2}}\sg_2^{-1}\quad\dots
\quad \sg_l^{k_{l,1}}\dots \sg_1^{k_{l,n_l}}\dl_2\,\right)
\end{eqn}
where $n_i\ge 1$ odd and $k_{i,j}\ge 2$ when $1<j<n_i$ and
$k_{i,j}\ge 1$ when $j=1$ or $n_i$.
Then $\md_a\Lm(D)=-l$.
\end{lemma}

\proof For $l=1$ see lemma \reference{lemma4} below. Then use induction
on $l$. Change a crossing of a $\sg_2^{-1}$ in $D=D_-$. Then:

$D_+$ can be reduced by the same argument as in the proof of
lemma 5.2 of \cite{ntriv} until we have a form \eqref{Df} with
smaller $l$, and a non-zero number of negative kinks added.
Then since negative kinks shift the $a$-degree of $\Lm$ up,
we have $\md_a\Lm(D_+)>-l$.

$D_0$ has $\md_a\Lm(D_0)=1-l>-l$ by induction.

$D_\infty=D_{l_1}\# D_{l_2}$ with $l_1+l_2=l$, so $\md_a\Lm(D_
\infty)=-l_1-l_2=-l$, and $\md_a\Lm(D)$ is inherited from
$\Lm(D_\infty)$.
\qed

\begin{lemma}\label{lemma2}
Let 
\begin{eqn}\label{Df2}
D=\wh{\es}\,\left(\,\sg_1^{k_{1,1}}\sg_2^{k_{1,2}}\dots
\sg_1^{k_{1,n_1}}\sg_2^{-1}\quad
\sg_1^{k_{2,1}}\dots \sg_1^{k_{2,n_2}}\sg_2^{-1}\quad\dots
\quad \sg_l^{k_{l,1}}\dots \sg_1^{k_{l,n_l}}\sg_2^{-1}\,\right)
\end{eqn}
where $n_i\ge 1$ odd and $k_{i,j}\ge 2$ when $1<j<n_i$ and
$k_{i,j}\ge 1$ when $j=1$ or $n_i$.
Then $\md_a\Lm(D)\ge -2$ when $l\le 2$ and $\md_a\Lm(D)=-l$
when $\l\ge 3$.
\end{lemma}

\proof Assume first we proved the result for $l\le 2$, and that
$l>2$. We argue by induction on $l$.

Apply the $\Lm$-relation at a $\sg_2^{-1}$ crossing in $D=D_-$.
Then $D_0$ has $\md_a\Lm=1-l$ by induction. $D_+$ simplifies
as in the proof of lemma 5.2 of \cite{ntriv}. This simplification
only removes negative letters. It can be iterated until one
of two situations occurs. It can (a) happen that all negative
letters disappear. Then we have a positive braid and $\md_a\Lm=-2>-l$
by Yokota's result. Or it can (b) occur that no $\sg_1\sg_2\sg_1$
or $\sg_2\sg_1\sg_2$ occur as subwords. Then $k_{i,j}\ge 2$
for $1<j<n_i$, and the number $l$ of negative crossings has
decreased strictly. So we have by induction $\md_a\Lm(D_+)>-l$.
Finally we must deal with the $D_\infty$ term. This follows
from lemma \reference{lemma1}.

It remains to justify the claim $\md_a\Lm(D)\ge -2$ when $l\le 2$.
This is done exactly with the same argument, only that now we
observe that all of $\Lm(D_{0,+,\infty})$ have $a$-degree $\ge -2$. \qed

\begin{lemma}\label{lemma3}
The closed braids in \eqref{Df2} are not positive for $l\ge 2$.
\end{lemma}

\proof The representation \eqref{Df2} clearly gives rise to
a positive band representation by replacing $\sg_2\sg_1^{k_{i,n_i}}
\sg_2^{-1}$ by $\sg_3^{k_{i,n_i}}$. So $[\bt]=3-\chi(\hat\bt)$.
If $l>2$, then we have
\[
\md_a F(D)\,=\,w(D)\,+\,\md_a \Lm(D)\,=\,3-\chi(\hat\bt)+\md_a \Lm(D)\,
<\,1-\chi(\hat\bt)\,
\]
but $\md_v P(D)=1-\chi(\hat\bt)$ as before, so $\md_a F\ne \md_v P$
and we are done by theorem \reference{thF}.

If $l=1$, then if $n_1=1$, we have a $(2,k_{1,1})$-torus
link, and for $n_1=3$ we have the $(1,k_{1,1},k_{1,2},k_{1,3})
$-pretzel link. Otherwise we apply the relations for $\tl P$
and $\Lm$ at the negative crossing. Then $\md_a\Lm(D_{+,0})=
\md_a\tl P(D_{+,0})=-2$, and
$[\Lm(D_{+,0})]_{a^{-2}}=[\tl P(D_{+,0})]_{a^{-2}}$
by theorem \reference{thF} since $D_{+,0}$ are positive braids.
That the extra term $\Lm(D_\infty)$ has no contribution to
$a^{-2}$ follows from lemma \reference{lemma4} below. So
$D$ satisfies Yokota's conditions in theorem \reference{thF}.
(We do not always know if $D$ depicts a positive link; see
remark \reference{rty}.)

If $l=2$ then resolve a negative crossing in $D=D_-$ via the
relation $\tl P_-=z \tl P_0-\tl P_+$ and via the $\Lm$-relation
$\Lm_-=z\Lm_0+z\Lm_\infty-\Lm_+$. Now $[\tl P_{+,0}]_{a^{-2}}=
[\Lm_{+,0}]_{a^{-2}}$. This follows from the above argument
for $l=1$. The additional term $\Lm_\infty$ has $a$-degree
$-2$ by lemma \reference{lemma1}, so $[\Lm_-]_{a^{-2}}\ne
[\tl P_-]_{a^{-2}}$. Again by theorem \reference{thF},
$D=D_-$ can therefore not belong to a positive link.  \qed

\begin{rem}
Ishikawa asks in \cite{Ishikawa} whether strongly
quasi-positive knots satisfy the equality $TB=2g_s-1$,
where $TB$ is the maximal Thurston-Bennequin invariant.
The proof of the lemma shows that there are infinitely many
strongly quasi-positive 3-braid knots with $TB<2g_s-1$:
take any of the knots with $l>2$. (There are also many other
such knots, like the counterexamples to Morton's conjecture
in \cite{posex_bcr}.)
\end{rem}

\begin{lemma}\label{lemma4}
If 
\[
D=D_{[n]}=\wh{\es} (\sg_1^{k_1}\sg_2^{k_2}\dots\sg_1^{k_n}\dl_2)
\]
with $n\ge 1$ odd, $k_1,k_n\ge 1$ and $k_i\ge 2$ when $i=2,\dots,n-1$,
then $\md_a \Lm(D)=-1$. Writing $\ol{k}=(k_1,\dots,k_n)$ and
\[
\tl w(\ol{k})\,=\,(k_1-1)+\sum_{l=2}^{n}(k_l-2)\,, 
\]
we have $\Md_z [\Lm(D)]_{a^{-1}}\,=\,\tl w(\ol{k})$.
\end{lemma}

\proof If $n=1$ we check directly (we have a reduced diagram
of the $(2,k_1)$-torus link), so let $n\ge 3$.

Consider first three special forms of $\ol{k}$.

If $\ol{k}=(12^*1)$ (with $2^*$ being a sequence of `$2$'),
then $D$ is regularly isotopic to a trivial
2-component link diagram. If $\ol{k}=(12^*)$ or $(2^*1)$, then $D$
is regularly isotopic to an unknot diagram with one positive kink.
In both situations the claims follow directly.

Now let $\ol{k}=(2^*)$. We orient $D$ so as to become negative.
Then $D$ depicts the $(2,-n-1)$-torus
link. We can evaluate $\Lm$ on $L$ from $F(L)$ by
normalization. For its mirror image $!L$,
we can use theorem \reference{thF} to conclude that
$\md_a F(!L)=n$. Now it is also
known that $\spn_a F(!L)=c(!L)$, and $c(!L)=n+1$, so
\[
\Md_a F(!L)=\md_a F(!L)+\spn_a F(!L)=n+(n+1)=2n+1\,.
\]
Thus $\md_a F(L)=-1-2n$, and since $w(D)=-2n$, we have
\[
\md_a \Lm(D)=-w(D)+\md_a F(D)=-(-2n)-1-2n=-1\,.
\]
That $\Md_z[\Lm(D)]_{a^{-1}}=1$ can also be obtained by
direct calculation.

If $\ol{k}$ is not of these special types, then $k_l\ge 3$ for some
$1\le l\le n$. We resolve a positive crossing in $\sg_j^{k_l}$.
We have (with $D=D_+$)
\begin{eqn}\label{*}
\Lm (D_+)\,=\,z\Lm (D_0)\,+z\Lm(D_\infty)-\Lm(D_-)\,.
\end{eqn}
Here $D_0$ has $\tl w$ by one less and comes with a $z$-factor,
so it is enough to show that $\Lm(D_-)$ and $\Lm(D_\infty)$
do not contribute.

If $k_l>3$, then $D_-$ is of the required form and has $\tl w$
by two less, so by induction $\Lm(D_-)$ has too small $z$-degree
in $[\Lm]_{a^{-1}}$.

Now consider $k_l=3$. As in the proof of lemma 5.2 of \cite{ntriv} we
can move
by braid relations (regular isotopy) a $\sg_1$ in $\sg_1\sg_2\sg_1$
until we get a $\sg_i\sg_i^{-1}$ (and then cancel). Here we
replaced $\sg_i^{-1}$ by $\dl_i$, so that such $\sg_i$ right before
a $\dl_i$ becomes a kink, which is negative. By repeating this
transformation, we obtain a form that has no $\sg_i\sg_{i\pm 1}\sg_i$
(i.e. $k_i>1$ for $1<i<n$) and a certain non-zero number of
negative kinks collected at both ends. The negative kinks
shift the degree in $a$ of $\Lm$ up, so by induction $D_-$
has no contribution to $[\Lm]_{a^{-1}}$.

It remains to deal with the term of $D_\infty$ in \eqref{*}.
It is the connected sum $D_{[n_1]}\# D_{[n_2]}$ with $n_1+n_2=n-1$,
and with $k_l-1>1$ negative kinks. So by induction $\md_a\Lm(D_\infty)
=k_l-1+(-1)+(-1)>-1$, and $\Lm(D_\infty)$ gives no
contribution to $[\Lm(D_+)]_{a^{-1}}$. \qed

\begin{theorem}\label{thposq}
If a 3-braid link is positive, then it is the closure
of a positive or almost positive 3-braid. Along these links
the non-fibered ones are exactly the $(1,p,q,r)$-pretzel
links.
\end{theorem}

\proof Consider the first claim.
% If a 3-braid is almost positive, and does not
% reduce to a positive one, then it is of the form
% \eqref{Df2} with $l=1$. % We observed in the proof
% of lemma \reference{lemma3}, that such links are positive.
% 
Since $L$ is positive, it is strongly
quasi-positive, and so has a positive 3-braid band representation
$\bt$ by theorem \reference{thqp}. So $\bt$ is of Xu's form $R$
or $(21)^kR$ (with $k>0$) up to extensions. If $\bt$ contains
$\sg_1\sg_2\sg_1$ or $\sg_2\sg_1\sg_2$ we can reduce it
as before, until (a) it becomes positive, or (b)
it still has a positive band representation, but it
does not contain $\sg_1\sg_2\sg_1$ or $\sg_2\sg_1\sg_2$.
In case (a) we are done. In case (b) we observe that
$\bt$ is of the form \eqref{Df2}, apply lemma
\reference{lemma3}, and conclude that $l\le 1$.

It remains to argue which links are not fibered.
Positive braids are always fibered, and by direct
observation for an almost positive braid we have
Xu's form $R$ or $(21)^kR$ depending on whether
in \eqref{Df2} (with $l=1$) we have $n_1=3$ or
$n_1>3$. We proved in theorem \reference{Thq}
that the forms $(21)^kR$ (for $k>0$) give fibered
closure links. For $n_1=3$ we have as before the
$(1,p,q,r)$-pretzel links. \qed

\begin{rem}\label{rty}
Unfortunately, we cannot completely determine which of the almost
positive braid representations with $l=1$ and $n_1>3$ in \eqref{Df2}
give positive links. For example the Perko knot $10_{161}$ has
such a representation, and it is positive. The Perko move (see
\cite{HTW}), turning the closed braid diagram of $10_{161}$
into a positive diagram, applies for more general examples.
However, some knots, like $14_{46862}$, do not seem subjectable
to this or similar moves, and their positivity status remains
unclear at this point.
\end{rem}

Similarly, one would hope to prove that
among these links none is a positive braid link (and this way
to obtain a different, but much more insightful proof of theorem
\reference{thg}). Using the polynomials, one can exclude certain
families, for example all links of $n_1=5$, but a complete argument
again does not seem possible.

\section{Studying alternating links by braid index\label{Sfg}}

The combination of the identity \eqref{p21} and the skein-Jones
substitution \eqref{VPsub} was already used in \S\reference{Js}
to translate the determination of the 3-braid link genus from
$P$ to $V$. A similar line of thought will now enable us to
extend the other main result in \cite{3br}, the description
of alternating links of braid index 3. This result was
motivated by the work of Murasugi \cite{Murasugi},
and Birman's problem in \cite{MortonPb} how to
relate braid representations and diagrammatic properties
of links. We will see how via \eqref{VPsub} and 
the famous Kauffman-Murasugi-Thistlethwaite theorem
\cite{Kauffman2,Murasugi3,Thistle2} the Jones polynomial
enters in a new way into the braid representation picture.
The argument will lead to the braid index 3
result surprisingly easily, and then also to the
classification for braid index 4 (which seems out of scope
with the methods in \cite{3br} alone). We also obtain a
good description of the general (braid index) case.

Our starting point is the following general result concerning the
MWF-bound \eqref{mwfb}. A diagram is called special if it
has no separating Seifert circles; see \cite{Cromwell}.
The number of Seifert circles of $D$ is denoted by $s(D)$.

\begin{theorem}\label{_ty}
Assume $L$ is a non-trivial non-split alternating link, and $MWF(L)=
k$. Then an alternating (reduced) diagram $D$ of $L$ has $s(D)\le 2k-
2$ Seifert circles, and equality holds only if $D$ is special.
\end{theorem}

\proof As $L$ is non-trivial and non-split, we have $1-\chi(L)> 0$.
It is also well known, that $\Md_zP(L)=1-\chi(L)$ (see
\cite{Cromwell}). Now $\spn_v P(L)\le 2k-2$ by assumption. Under
the substitution in \eqref{mwfb} this translates to 
\begin{eqn}\label{sV}
\spn V(L)\le 1-\chi(L)+2k-2\,.
\end{eqn}
On the other hand, by \cite{Kauffman2,Murasugi3,Thistle2},
$\spn V(L)=c(D)$, and also $1-\chi(L)=c(D)-s(D)+1$.
So
\begin{eqn}\label{ff}
c(D)\le c(D)-s(D)+1+2k-2\,,
\end{eqn}
and then $s(D)\le 2k-1$. Now if this is an equality, then
so is \eqref{sV}. Then one easily sees that $P(L)$ must have
non-zero coefficients in both monomials $z^{\Md_zP}v^{\Md_vP}$
and $z^{\Md_zP}v^{\md_vP}$. Now under the substitution in \eqref{p21},
both these monomials give a non-cancelling contribution,
and one of them is not in ($v$-)degree 0, so the identity
\eqref{p21} cannot hold. Moreover, if \eqref{sV} fails
just by one, then still one of the two coefficients must
be non-zero. In order its (non-cancelling) contribution on
the left of \eqref{p21} to be in degree $0$, we see that
either $\Md_zP=\md_vP$, or $\Md_zP=-\Md_vP$. Using \cite
{Cromwell}, one concludes then that these are precisely the cases of a
(positively or negatively) special alternating link. \qed

Using \cite{Yamada} and \cite{Vogel} we have some simple
estimate on the unsharpness of MWF for alternating links.

\begin{corr}
For an alternating non-trivial non-split link $L$, we have
$b(L)\le 2MWF(L)-2$. \qed
\end{corr}

Of course, for many (in particular alternating) links
$b(L)=MWF(L)$ (that is, MWF is exact), or at least $b-MWF$
is small. So the above estimate should be considered as
a worst-case-analysis. Even if not strikingly sharp,
it is still far from trivial, in view of what we
already know can occur for non-alternating knots. Namely,
using the construction in \cite{Kanenobu2} and the work
in \cite{BirMen3} (see remark \reference{rfn} below),
one can find sequences of knots $(K_i)$ for which $MWF$
(in fact the full $P$ polynomial) is constant, but
$b(K_i)\to\infty$.

Another observation is that one can now extend the outcome
of the work in \cite{SV} by replacing crossing number of
an alternating knot by its braid index.

\begin{corr}
Let $n_{g,b}$ denote the number of alternating knots of
braid index $b$ and genus $g$. Then $n_{g,b}$ is finite.
Moreover, for $g$ fixed, we have that
$\lim_{b\to \infty} n_{g,b}/b^{6g-4}=C_g$ is a constant.
\end{corr}

\proof The finiteness of $n_{g,b}$ follows by \eqref{ff}.
When $\chi(L)$ is fixed, and $c(L)\to\infty$, one easily
sees from \eqref{ff} that $MWF(L)$ behaves asymptotically
(up to an $O(1)$, i.e. bounded, term)
at least like $c(D)/2$. Now Ohyama's result \cite{Ohyama}
implies that $b(L)$ (and $MWF(L)$ as well) behave asymptotically
exactly as $c(L)/2$. So from \cite{SV} we have the result.
\qed

\begin{rem}\label{rfn}
Of course we could also gain, as in \cite{SV}, an estimate on the
$C_g$ and its asymptotics for $g\to \infty$, it would just multiply
by $2^{6g-4}$. One should also note that the finiteness of
$n_{g,b}$, which one sees from \eqref{ff}, is not necessarily
clear \em{a priori}. In fact, however, Birman-Menasco proved
\cite{BirMen3} that $n_{g,b}$ is a finite number even for \em{general}
(i.e. without restriction to alternating) knots. Their methods seem,
though, quite unhelpful to estimate these numbers properly.
\end{rem}

Theorem \ref{_ty} immediately leads to the first slight sharpening
of the description of alternating links of braid index 3 in
\cite{3br}. (The case of $MWF=2$ is even more obvious, and omitted.)

\begin{corr}\label{opi}
An alternating link has $MWF=3$ if and only if it has braid index 3.
\end{corr}

\proof If $MWF(L)\le 3$, then the alternating diagram has at
most 4 Seifert circles, and exactly 4 only if it is special.
Apart from connected sums (which are easily handled), we obtain
the diagrams of closed alternating 3-braids and the $(p,q,r,s)$-
pretzel diagrams. By direct calculation of $P$ we saw in
\cite{3br} that if $\min(p,q,r,s)\ge 2$, then $MWF\ge 4$,
and that otherwise the pretzel link has braid index 3. \qed

The case of 4-braids is now not too much more difficult.
% (even though
% it was completely out of scope with the method in \cite{3br}).

\begin{figure}[htb]
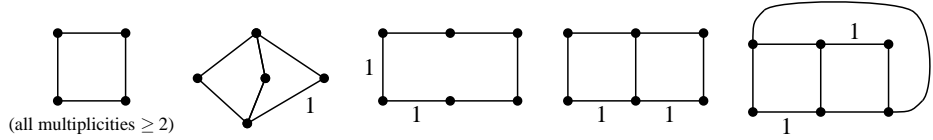

{\small
\[
\begin{array}{c}
  \diag{6mm}{3}{2}{
  \pictranslate{0 0.5 x}{
    \diag{6mm}{3}{2}{
      \cycl{0.5 2}{0.5 0.5}{2 0.5}{2 2}
    }
  }
  \picputtext{1.75 0}{\scriptsize (all multiplicities $\ge 2$)}
  }
  \qquad
  \diag{6mm}{3}{2}{
    \cycl{1.3 2}{0 1}{1.1 0}{1.5 1}
    \cycl{1.3 2}{2.8 1}{1.1 0}{1.5 1}
    \picputtext{2.5 0.4}{$1$}
  }
  \qquad
  \diag{6mm}{3}{2}{
    \pictranslate{0 0.5}{
    \xcycl{0 0}{1.5 0}{3 0}{3 1.5}{1.5 1.5}{0 1.5}
    \picputtext{0.75 -0.3 x}{$1$}
    \picputtext{0.75 -0.3}{$1$}
    }
  }
  \qquad
  \diag{6mm}{3}{2}{
    \pictranslate{0 0.5}{
    \xcycl{0 0}{1.5 0}{3 0}{3 1.5}{1.5 1.5}{0 1.5}
    \picline{1.5 1.5}{1.5 0}
    \picputtext{2.25 -0.3}{$1$}
    \picputtext{0.75 -0.3}{$1$}
    }
  }
  \qquad
  \diag{6mm}{4}{2.5}{
    \pictranslate{0 0.5}{
    \xcycl{0 0}{1.5 0}{3 0}{3 1.5}{1.5 1.5}{0 1.5}
    \picline{1.5 1.5}{1.5 0}
    \picputtext{2.25 1.8}{$1$}
    \picputtext{0.75 -0.3}{$1$}
    \opencurvepath{0 1.5}{0 2}{0.5 2.5}{3 2.5}{4 2}{4 0.5}{3.5 0}{3 0}{}
    }
  }
  \\[5mm]
\end{array}
\]
}
\caption{The reduced Seifert graphs of the alternating diagrams
of special alternating links of braid index 4. Simple edges
have their multiplicity (1) attached, and the other edges (of
multiplicity one or more) are unlabelled. In the first graph,
exceptionally, all edges have multiplicity at least 2.\label{figSG}}
\end{figure}

\begin{theorem}\label{tty}
Let $L$ be a prime non-split alternating link. The following 3
conditions are equivalent:
\begin{enumerate}
\item\label{Xa} $MWF(L)=4$
\item\label{Xb} $b(L)=4$
\item\label{Xc} $L$ is one of the links, whose reduced 
alternating diagrams are described (up to mirror images) as follows
\begin{enumerate}
\item\label{Ya} The Murasugi (or connected) sum of three
  $(2,n_i)$-torus links (with $|n_i|>1$),
\item\label{Yb} The Murasugi (or connected) sum of a $(2,n)$-torus
  link with a $(p,q,r,s)$-pretzel link, with one of $p,q,r,s$ equal
  to 1, or its mirror image, or
\item\label{Yc} a special diagram whose Seifert graph
(see \S\ref{Sdg}) is as shown in figure \reference{figSG}.
\end{enumerate}
\end{enumerate}
\end{theorem}

\proof $\reference{Xb}\So \reference{Xa}$. This follows from
corollary \ref{opi}.

$\reference{Xc}\So \reference{Xb}$. That the links in
\reference{Xc} have braid index at least 4 follows from the
description of the links with $b(L)\le 3$ in \cite{3br}
(which also comes out of the proof of corollary \eqref{opi}).
It is also not too hard to check that for all these links
$b=4$ by exhibiting a diagram $D$ with $s(D)=4$. The most
systematic way seems to apply the graph index inequality of
Murasugi-Przytycki \cite{MurPrz2} (see also \cite{Ohyama}).

$\reference{Xa}\So \reference{Xc}$. By applying theorem
\reference{_ty}, we need to deal with non-special
diagrams of at most 5 and special diagrams of at most 6
Seifert circles.

First consider the non-special diagrams.

For the fibered links (the reduced Seifert graph is a tree),
Murasugi's result \cite{Murasugi} leads directly to case
\reference{Ya}. For non-fibered links, the Seifert graph
must have a cycle, which must be of length at least 4 (the
Seifert graph is bipartite). Then the only option
that remains is case \reference{Yb}. We must still
argue why one of $p,q,r,s$ must be $\pm 1$. This can
be done using Murasugi-Przytycki's work, but one easily
sees it also by a direct skein theoretic argument, which
we explain.

Let $L(p,q,r,s)$ be the corresponding link and $P(p,q,r,s)$
its skein polynomial. Look first at $p=0$, $q,r,s\ge 2$. Then
$L(p,q,r,s)$ a connected sum of two $(2,n)$-torus
links, and a closed alternating 3-braid. So
$MWF(L(p,q,r,s))=5$. Now since the case $p=1$, $q,r,s\ge 2$ has
$MWF(L(1,q,r,s))\le b(L(1,q,r,s))=4$, the skein relation \eqref{srel}
easily shows that the maximal degree coefficient $\Mc_vP(p,q,r,s)$
of $v$ in $P$ (which is a polynomial in $z$) for $p=2$ is inherited
from $p=0$. Then further applications of \eqref{srel} show that
the $z$-degree of $\Mc_vP(p,q,r,s)$ increases with $p$, so 
in particular this term never vanishes, and so $MWF(L(p,q,r,s))=5$.

Now consider the special diagrams. For them one considers the
Seifert graph, and needs to write
down all bipartite planar graphs on at most 6 vertices, which
have no cut vertex. Since the diagram $D$ is special, the placement
of multiple copies of an edge give diagrams equivalent up to
flypes, so it is enough to consider simple graphs (the reduced
Seifert graph), and have the multiplicities of an edge written
as its label. The graphs can be easily compiled using the observation
that they must contain a cycle of length 4 or 6; see figure
\reference{figSG}. By direct inspection we see that the edge
multiplicities must be as specified in the figure. (In fact if an
edge is multiple it turns out irrelevant what its multiplicity
is, so in this case we just omit the label.) We rule out the
remaining multiplicities by a skein theoretic calculation,
similar to the one explained for case \reference{Yb}.

Let w.l.o.g. (up to mirroring) $D$ be positive. For each edge
$e$ in $\Gm(D)$ of variable multiplicity $\ge i$ (where $i=1,2$)
we calculate the skein polynomial of the diagram that
corresponds to $\Gm$ for multiplicities $i$ and $i+1$ of $e$.
That is, if $\Gm(D)$ has $l$ edges of variable multiplicity,
we have $2^l$ polynomials to calculate. Then we check for
each such set of $2^l$ polynomials that $Q=P_{v^{9-\chi(D)}}$
is non-zero, and $\Md_zP-\Md_zQ$ as well as $\Mc_zQ$ is
constant within this set of $2^l$ polynomials. Then
by \eqref{srel} this property is inherited to diagrams
$D$ whose $\Gm(D)$ have edges of higher multiplicity,
and in particular $MWF\ge 5$.

There is one more graph,
\[
\diag{6mm}{3}{2}{
  \cycl{2.0 2}{0 1}{1.8 0}{1.2 1}
  \cycl{2.0 2}{2.4 1}{1.8 0}{3.5 1}
}
\]
not included in figure \reference{figSG}. In that
case (by the method we just explained we verify that)
$MWF\ge 5$ for all non-zero edge multiplicities. \qed

In \cite{MurPrz2}, Murasugi-Przytycki define a certain
quantity $\inx(D)$, assigned to a link diagram $D$,
called \em{index}. (We omit here the detailed discussion;
one can consult also \cite{Ohyama} or \cite{gener}.)
Their motivation was to give an upper estimate 
\begin{eqn}\label{mp}
b(L)\le s(D)-\inx(D)
\end{eqn}
for the braid index of the underlying link $L$. Their origin
of \eqref{mp} consists in an appropriate move (see
figure 8.2 in \cite{MurPrz}) which reduces
the number of Seifert circles of the diagram.
Murasugi-Przytycki conjectured that for an alternating
diagram $D$, the inequality \eqref{mp} is exact. This
conjecture is also confirmed for alternating links up
to braid index 4. 

Based on the Murasugi-Przytycki procedure, we can re-enter
the Bennequin surface topic.

\begin{corr}\label{cfT}
Alternating links of braid index at most 4 carry a Bennequin
surface on a minimal string braid.
\end{corr}

\proof In \cite{gener} it was explained how to apply
restrictedly the Seifert circle reduction move of
Murasugi-Przytycki (and then those of Yamada \cite{Yamada})
so as to obtain a braided surface. This modified reduction
is easily checked to lead to the minimal number of strings for the
links in question. \qed

% (It was verified in \cite{gener}
% for alternating knots of genus up to 4.)

\begin{corr}
Let $L$ be an alternating link of braid index $4$. Then
$|\Mc\Dl(L)|\le 5$, and if $|\Mc\Dl(L)|>2$, then $L$ is
special alternating.
\end{corr}

\proof It is well known, that $|\Mc\Dl(L)|$ is multiplicative
under Murasugi sum and for a special diagram depends only
on the reduced Seifert graph. The result follows by calculation
for the specific types. \qed

The proof of this corollary, and the extension of the multiplicativity
of $\Mc\Dl$ to $\Mc_zP$ for diagrammatic Murasugi sum \cite{MurPrz}
demonstrates also the following more general principle:

\begin{corr}\label{cSG}
For any given braid index, there are only finitely many values of
$\Mc\Dl$ and $\Mc_zP$ among Alexander and skein polynomials of
alternating links of that braid index. \qed
\end{corr}

We saw that for $\Dl$ this statement is wrong for non-alternating
links even among 3-braids. On the other hand, by \cite{3br},
it is true for $P$, and we do not know if it remains true for
(closed) braids on more strings. (One could also ask if
infinitely many leading coefficients of $\Dl$ occur if $\Md\Dl
=1-\chi$, but we se no deeper meaning in this question, so will
not dwell further upon it here.)

The following is also worth observing. Call a subclass
$\cC'$ of a class $\cC$ of links \em{generic} in $\cC$ if
\[
\lim_{n\to\infty}\,\frac{\#\,\{L\in \cC'\,:\,c(L)=n\,\}}
{\#\,\{L\in \cC\,:\,c(L)=n\,\}}\,=\,1\,.
\]

\begin{corr}
The number of special alternating links of given braid index
grows polynomially in the crossing number. In particular,
a generic alternating link of given braid index is not
special alternating.
\end{corr}

\proof A special alternating link is determined by the Seifert
($=$ checkerboard) graph of its alternating diagram. The number
of such graphs with a fixed number of vertices grows polynomially
in the number of edges. The second claim in the corollary
follows because it is easy to see that in contrast the number
of non-special alternating links grows exponentially (due to
exponentially many, in the crossing number, ways to perform the
Murasugi sum at a separating Seifert circle of the alternating
diagram). \qed

\begin{rem}
Note that in contrast we showed in \cite{SV} that a generic
alternating knot (and the case of links is analogous) of given
\em{genus is} special alternating. This shows from yet another
point of view the opposition between genus and braid index.
\end{rem}

Another immediate and useful consequence of theorem \ref{tty}
and the preceding remarks is

\begin{corr}\label{cft}
If an alternating link $L$ has $MWF(L)\le 4$ (in particular
if $b(L)\le 4$), then the MWF inequality is exact
(i.e. an equality) for $L$. \qed
\end{corr}

This gives a nice complement to the MWF exactness results in
\cite{Murasugi,MurPrz2}. (As another such amplification, we 
proved the case of \em{knots} and genus $\le 4$ in \cite{gener}.)

Note that at MWF bound 5 we hit already at the Murasugi-Przytycki
examples \cite{MurPrz2} of non-exact MWF (with $b=6$). So corollary
\reference{cft} is not true for $MWF\ge 5$ or $b\ge 6$. We do
not know about the case $b=5$. However, ruling out braid index
$5$ for the Murasugi-Przytycki family is a serious computational
problem (only two specific members were dealt with, a 4-component
15 crossing link and an 18 crossing knot; see \cite[\S 19]
{MurPrz2}). Already with this circumstance in mind,
one cannot expect to easily extend the corollary (or theorem
\reference{tty}) for $b=5$ either, even if it may be true.

On the other hand, leaving these troublesome exceptions aside,
the above discussion should fairly clearly explain how
the general picture continues for alternating links with MWF
bound 5 and more.

\section{Applications of the representation theory\label{SBurau}}

So far the representation theory behind the skein, Jones, and
Alexander polynomial was not used. We will give some applications
of it now. The theory is well explained in \cite{Jones}. We will
use Jones's conventions, unless otherwise specified.

\subsection{The Jones conjecture\label{SJS}}

There is a conjecture, often attributed to Jones (who speculated
on it at least for knots and in it its weaker form, as given below;
see 357 l.-6 of his paper), stating that a minimal
string braid representation has a unique exponent sum.

\begin{conj}(Jones's conjecture)\label{jcc}
\def\labelenumi{\theenumi)}
\begin{enumerate}
\item (weaker version)\label{weV}
If $\bt,\bt'\in B_n$ satisfy $\hat\bt=\hat\bt'=L$ and $n=b(L)$,
then $[\bt]=[\bt']=:w_{\min}(L)$.
\item (stronger version) Part \reference{weV} holds, and if
$\bt''\in B_n'$ for $n'>n$ has $\hat\bt''=L$, then 
\begin{eqn}\label{stz}
\big|[\bt'']-w_{\min}(L)\big|\le n'-n\,.
\end{eqn}
\end{enumerate}
\end{conj}

It was observed in \cite{posex_bcr} that counterexamples
to the conjecture would make MWF and \em{all} of its cabled
versions unsharp. (In \cite{posex_bcr} the weaker version was
focussed on, but the same arguments address the stronger version
too.) Thus, for example, corollary \reference{cft} can also
be regarded as a partial solution to Jones's conjecture.
Similarly, the work in \S\reference{S2.2} shows:

\begin{corr}
If a link $L$ is the closure of a positive braid, and $b(L)\le 4$,
then the strong Jones conjecture is true for $L$. 
\end{corr}

\proof We use theorem \reference{thmwf}, with the remark that
the exceptional braid words therein were shown to have braid
index 4 using the 2-cabled MWF. \qed

For 3-braids the weaker version was known to be true again
from Birman-Menasco's classification result. We will show now
the stronger version, as a consequence of the description in
\cite{3br} of 3-braid links with unsharp MWF inequality (and thus
with a much simpler proof than appealing to Birman-Menasco). See
also remark \reference{rJC}.

\begin{theorem}\label{thJC}
The stronger version of Jones's conjecture holds for 3-braid links.
\end{theorem}

\noindent{\bf Convention.}
The letter $\Dl$, with an integer subscript $n$, is used in
\S\ref{SJS},
in deviation from other sections, exclusively for the half-twist
braid on $n$ strings, and \em{not} for the Alexander polynomial.
(The Alexander polynomial will not appear in \S\reference{SJS}.)

\proof[of theorem \reference{thJC}]
As explained, it suffices to deal only with the 3-braid links
of unsharp MWF inequality. In \cite{3br} these links were
described fully. If $\mu=\sg_1^{6a\pm 1}\sg_2^{\mp 1}$,
then for the Birman dual $\bt=\mu^*$ of $\mu$ (see definition
\ref{DEFB}) we have $\hat\mu^*\ne \hat\mu$ when $6a\pm 1>6$, but
$K=\hat\bt$ has the skein polynomial of the $(2,6a\pm 1)$-torus
knot $\hat \mu$. In \cite{3br} is was proved that these knots $K$
form the full list of 3-braid links with unsharp MWF inequality.
 
Consider first %(with a slight change of notation)
the knots $K=\hat\bt$ with $\bt=\mu^*$ for
$\mu=\sg_1^{6a+1}\sg_2^{-1}$ and $a>0$. It is easy
to bring $\bt$ into Xu's normal form, and to observe
that it has exactly one negative band. This exhibits
a Seifert ribbon (see \cite{Rudolph}) of smaller genus
than $g(K)$, i.e. $g_s(K)<g(K)$ (see also the proof of
theorem \reference{thqp}). Now from Xu's form one sees $w_{\min}
(K)=[\bt]=2g(K)$. Take this for a moment as a definition of
$w_{\min}(K)$; this writing will be justified when we
show that the writhe is unique. (One could also quote
Birman-Menasco here, but the complexity of their argument
is unnecessary.) So 
\begin{eqn}\label{iopa}
1-\chi_s(L)<w_{\min}(L)\,.
\end{eqn}
Using \eqref{mq},
we see that if \eqref{stz} is violated for some $\bt''\in B_{n'}$,
then $[\bt'']=w_{\min}(L)+2+(n'-3)$. However, by the Rudolph-%
Bennequin inequality \cite{Rudolph2} (see also \cite{pos}),
we would have then 
\[
1-\chi_s(L)\ge [\bt'']-n'+1=w_{\min}(L)\,,
\]
in contradiction to \eqref{iopa}.

The effort focusses on the other braids $\bt=\mu^*$
for $\mu=\sg_1^{6a-1}\sg_2$.
In this case we consider $!K$, whose Xu form has only negative
bands. The inequalities \eqref{mq} show that we must rule out a
braid representation $\bt''\in B_{n'}$ of $!K$ with $[\bt'']=w_{\min}
(!K)+2+(n'-3)$, where (we set in the same \em{a posteriori} to
be justified manner) $w_{\min}(!K)=-2g(K)-2$.

We prohibit $\bt''$ by evaluating the 2-cable skein polynomial
$P(\hat\gm_a)$ of $!K$ for
\begin{eqn}\label{gma}
\gm_a\,=\,[4354]\cdot [2132]^{6a-1}\,(\Dl_6^2)^{-2a}\,\in\,B_6,
\end{eqn}
with $a>1$ and $\Dl_6^2=[12345]^6$ generating the center of $B_6$,
and showing that the polynomial has non-zero terms in $v$-degree
$[\gm_a]-5$ or $[\gm_a]-3$. This may appear a banality, but in
fact requires a substantial use of Jones' work in \cite{Jones}.

Consider the more general 6-braids
\[
\bt_{k,l}\,:=\,[4354]\cdot [2132]^k\,\Dl_6^{2l}\,.
\]
(We are of course only interested in the special case $k=6a-1$,
$l=-2a$ and $a>1$, but it is useful to treat the 2-parameter family
first.)

The skein polynomial of $\hat\bt_{k,l}$ can be evaluated using
the representation theory in \cite{Jones}. We adopt the convention
that \em{all references} to lemmas, page and paragraph numbers,
equations of the form ($x.y$), etc. for the rest of this proof
are understood to be \em{to Jones's paper}, unless noted otherwise.

Jones's % paper gives a clear and detailed account on the
% calculation of the skein polynomial. His
version of the skein polynomial, 
\[
X(q,\lm)=P(v,z)\mbox{\quad with $v=\sqrt{\lm q}$\ \ and
\ $z=\sqrt{q}-1/\sqrt{q}$\,},
\]
can be evaluated on a closed
$n$-braid $\bt$ by a weighted sum, with weights $W_Y$ in $\lm$ and
$q$, of traces (in $q$) of irreducible representations (irreps) $\pi_Y$
of $B_n$, indexed by Young diagrams (or tableaux\footnote{Here
consistently `tableau' is used as a synonym for (Young) diagram,
i.e. with no additional information attached to it.}) $Y$, or
equivalently by partitions of $n$. (A partition of $n$ is a
tuple $(n_1,\dots,n_k)$ with $n_k>0$, $n_i\ge n_{i+1}$ and
$\sum n_i=n$.) We may identify the Young diagram with its
partition, counting partitions in horizontal rows. For
example the partition $(4)=\Young{4}$ means one row (trivial
representation), while $(1111)=\Young{1,1,1,1}$ means one column
(parity representation); thus $\pi_{3,1}=\pi_{\Young[1.5mm]
{3,1}}$. Since the calculation of $W_Y$ was given in
\cite{Jones}, we will not repeat it in detail. We will
also deal (mostly) with $Y$ where the calculation of
$\pi_Y$ is explained in \cite{Jones}.

Using Definition 6.1 (of \cite{Jones}), (5.5) and the formula
for $W_Y(q,\lm)$ from p. 347 top, one has for a 6-braid $\bt$
of exponent sum $e=[\bt]$,
\begin{eqn}\label{STR}
X_{\hat\bt}(q,\lm)\,=\,-\sqrt{\lm}^{e-5}\,\sum_{Y\vdash 6}
\wt W_Y(q,\lm)\,\tr \pi_Y(\bt)\,.
\end{eqn}
Here $\wt W_Y$ denote the slightly rescaled weights
\[
\wt W_Y(q,\lm)\,=\,\frac{R_{Y}(q,\lm)}{Q_Y(q)}\,\cdot\,
\frac{1-q}{1-\lm q}\,,
\]
with $R_Y$ being specified on p.\ 347 after Figure 5.6, and
$Q_Y(q)$ being the hook length product term from p. 346 middle.
The symbol `$\mathord{\vdash}$' is taken from partition theorists
and is used to mean here that $Y$ is a Young tableau of $6$ boxes
(or equivalently a partition of $6$).
In order to avoid denominators it is useful to multiply
\eqref{STR} by $(1-q^2)\cdot\dots\cdot(1-q^6)$, so
\begin{eqn}\label{STR2}
-(1-q^2)\cdot\dots\cdot(1-q^6)X_{\hat\bt}\,=\,
\sqrt{\lm}^{e-5}\,\sum_{Y\vdash 6}
\wh W_Y(q,\lm)\,\tr \pi_Y(\bt)\,,
\end{eqn}
with 
\[
\wh W_Y(q,\lm)\,=\,\frac{R_Y(q,\lm)}{1-\lm q}\,\cdot\,
\frac{(1-q)(1-q^2)\cdot\dots\cdot(1-q^6)}{Q_Y(q)}\,
\]
becoming a Laurent polynomial in $\lm$ and $q$\,.

Note that $\pi_Y$ are representations that involve only the
variable $q$, not $\lm$. The obvious question is how to
evaluate their traces.

Consider first 
\[
\bt_k=\bt_{k,0}\,=\,[4354]\,[2312]^k\,.
\]

There are 11 Young diagrams $Y$ of 6 boxes. Now $\bt_k\in B_4\subset
B_6$, and according to p. 340 top, the restriction of $\pi_Y$ to $B_4$
splits into a direct sum of representations $\pi_{Y'}$ indexed by
4-box Young diagrams $Y'$, which we name for shorthand
\begin{eqn}\label{shh}
A=\Young{1,1,1,1}\quad
B=\Young{2,1,1}\quad
C=\Young{2,2}\quad
D=\Young{3,1}\quad
E=\Young{4}\,.
\end{eqn}

The representations $\pi_{Y'}$ are clarified completely in \S 8
(see also \S\reference{SUn} in \em{this} paper below). Denote by
$\psi_{n-1}$ the (reduced) $(n-1)$-dimensional representation
of $B_n$; see \S 2. Let us also write $-1$ for the parity
representation $(-1)^{[\,.\,]}$, and $-\rho=-1\otimes \rho$ for the
direct (tensor, or Kronecker) product of a representation $\rho$ with
the parity. Then $\pi_A=-1$, $\pi_E$ is given by $q^{[\,.\,]}$,
$\pi_B=-\psi_3$ (Note 5.7; the sign disappears here, though,
because for us always $e=[\bt]$ is even), $\pi_D=-\bigwedge^2\psi_3$
and $\pi_E=\pi_{\Young[1.5mm]{2,1}}\circ\ol{\ry{0.6em}\es}$,
with $\ol{\es\ry{0.6em}}$ being the homomorphism $B_4\to B_3$
given by $\ol{\sg_{1,2,3}\ry{0.6em}}=\sg_{1,2,1}$ (p. 355).
Also $\pi_{\Young[1.5mm]{2,1}}=-\psi_2$.

Let us write below $\ol{Y}$ for the \em{transposed} (or \em{dual})
Young diagram to $Y$, given by exchanging rows and columns.
For example, for $Y=\Young{3,1}$, we have $\ol{Y}=\Young{2,1,1}$.
The relation between $\pi_Y$ and $\pi_{\ol{Y}}$ is given in Note
4.6.

Now $\psi_3(\ap)$ and $\psi_2(\bar\ap)$ for 
\[
\ap=[2132]
\]
are easy to
calculate (see \S 2). We find former's eigenvalues to be $\pm t$ and
$-t^2$, and latter's $t$ and $t^3$. In particular, both matrices are
diagonalizable (because the eigenvalues are distinct for generic $t$).
Setting $t=q$ as in Note 5.7, then we
have the following table of eigenvalues of the $\pi_{Y'}(\ap)$
\[
\begin{array}{|c|c|c|c|c|}
\hline
\ry{1.3em}A & B & C & D & E \\[1mm]
\hline
\ry{1.5em}1& -q   & q^3 & -q^3 & q^4 \\
	 &  q   & q   &  q^3 &     \\
	 & -q^2 &     & -q^2 &     \\[2mm]
\hline
\end{array}
\kern2cm
\begin{array}{c|*{7}{|c}}
\ry{1.2em}i & 1 & 2 & 3 & 4 & 5 & 6 & 7 \\[1.5mm]
\hline  
\ry{1.2em}\dl_i & 1 & q^4 & q & -q & -q^2 & q^3 & -q^3\\[1mm]
\end{array}
\]
Let us number the 7 possible eigenvalues of $\pi_Y'(\ap)$ by
$\dl_i$ as shown on the right. 

Thus $\pi_Y(\ap)$ are all diagonalizable, with eigenvalues $\dl_i$.
The multiplicities of $\pi_{Y'}$ in $\pi_Y$ are also easy to
calculate. They are the number of descending paths from $Y'$ to
$Y$ in the Figure 3.3 (continued one more row to the bottom).
Table \reference{tab1} shows the multiplicities of $\pi_{Y'}$
(in the shorthand of \eqref{shh}) in $\pi_Y$ and the resulting ones,
which we write $M(i,Y)$, of the eigenvalues $\dl_i$ in $\pi_Y(\ap)$.

\begin{table}[ptb]
\captionwidth\vsize\relax
\newpage
\vbox to \textheight{\vfil
\rottab{ 
\hbox to \vsize{%\hss%\footnotesize
\hss
% \setbox\@tempboxa=\hbox{%
% \textwidth\vsize\relax
$\displaystyle
\begin{array}{|c||*{11}{c|}}
\hline\ry{1.0cm}%
Y & \Young{1,1,1,1,1,1} & % 1
\Young{2,1,1,1,1} &       % 2
\Young{2,2,1,1} &         % 3
\Young{2,2,2} &           % 4
\Young{3,1,1,1} &         % 5
\Young{3,2,1} &           % 6
\Young{3,3} &             % 7
\Young{4,1,1} &           % 8
\Young{4,2} &             % 9
\Young{5,1} &             % 10
\Young{6} \\[6mm]         % 11
\hline\ry{1.6em}%
Y' & A & 2A+B & A+2B+C & B+C & A+2B+D & 2B+2C+2D & C+D & B+2D+E & C+2D+E
& D+2E & E \\[2mm]
\hline\ry{1.6em}%
%  PARTITION #7 & #8 ARE SWOPPED FROM MY NOTES
%      #1 #2  #3  #4  #5  #6  #8  #7  #9  #10 #11
1    & 1 & 2 & 1 &   & 1 &   &   &   &   &   &   \\
q^4  &   &   &   &   &   &   &   & 1 & 1 & 2 & 1 \\
q    &   & 1 & 3 & 2 & 2 & 4 & 1 & 1 & 1 &   &   \\
-q   &   & 1 & 2 & 1 & 2 & 2 &   & 1 &   &   &   \\
-q^2 &   & 1 & 2 & 1 & 3 & 4 & 1 & 3 & 2 & 1 &   \\
q^3  &   &   & 1 & 1 & 1 & 4 & 2 & 2 & 3 & 1 &   \\
-q^3 &   &   &   &   & 1 & 2 & 1 & 2 & 2 & 1 &   \\[2mm]
\hline\ry{1.6em}%
d    & 1 & 5 & 9 & 5 & 10& 16& 5 &10 & 9 & 5 & 1 \\
r    & 0 & 1 & 3 & 2 & 4 & 8 & 3 & 6 & 6 & 4 & 1 \\
30r/d& 0 & 6 &10 &12 & 12& 15& 18& 18& 20& 24& 30\\[2mm]
\hline
\end{array}%
$%
% }%
% \@tempdima0.3\vsize\relax\advance\@tempdima by -0.3\wd\@tempboxa\relax
% \kern\@tempdima
% \hss
% \copy\@tempboxa
\hss}%
}{This table displays the Young tableaux $Y$ of 6 boxes, the
decomposition into $\pi_{Y'}$ for 4-box Young tableaux $Y'$ of the
sub-representation of $B_4$ in $\pi_Y$ (writing $Y'$ for $\pi_{Y'}$),
the multiplicities $M(i,Y)$
of the 7 possible eigenvalues $\dl_i$ of $\pi_Y(\ap)$, and the
quantities $d$ ($=\dim \pi_Y$) and $r$ (rank of the idempotent
$e_i$) occurring in Jones's lemma 9.3.\label{tab1}}{table}
\vss}
\newpage
\end{table}

Since with $\pi_{Y'}(\ap)$ also all $\pi_Y(\ap)$ are
diagonalizable, we have
\begin{eqn}\label{XX}
\tr \pi_Y(\bt_k)\,=\,\tr \pi_Y(\ap^k\cdot [4354])\,=\,
\sum_{i=1}^7\,c(i,Y)\,\dl_i^k\,,
\end{eqn}
with $c(i,Y)$ given as follows. Consider $\pi_Y([4354])$ in
the basis of $\pi_Y$ that diagonalizes $\pi_Y(\ap)$. Then to
obtain $c(i,Y)$, sum the $(j,j)$-entries of the matrix
of $\pi_Y([4354])$ over rows/columns $j$, for which the
$(j,j)$-entry of $\pi_Y(\ap)$ is $\dl_i$.

So the problem to evaluate $\tr \pi_Y(\bt_k)$ transforms into
the one to determine $c(i,Y)$. There are \em{a priori} 77
of those, given by combining $7$ eigenvalues $\dl_i$ with
11 Young diagrams $Y$. However, one immediately notes that
clearly $c(i,Y)=0$ when $M(i,Y)=0$. This leaves 45 of the 77 values.

If we can calculate $\tr \pi_Y$ for general 6-braids, then we
can use \eqref{XX} as a linear equation for $c(i,Y)$.
If $\pi_Y(\ap)$ has $l_Y$ different eigenvalues, we can
determine $c(i,Y)$ for that $Y$ by calculating the l.h.s. of
\eqref{XX} for $l_Y$ different values of $k$ and solving for $c(i,Y)$.
(Always $l_Y\le 6$, as evident from table \reference{tab1}.)

In case $Y$ is one of the one-hook diagrams like $\Young{4,1,1}$,
then $\pi_Y$ is by p. 354 bottom (a tensor product of a parity,
which disappears at even exponent sum, with) an exterior power
of $\psi_5$. Thus, as in \S 7 (of Jones's paper), one can evaluate
$\tr \pi_Y(\bt_k)$ from the characteristic polynomial of the Burau
matrix $\psi_5(\bt_k)$. We determined this way the corresponding
$c(i,Y)$, and 
verified them with a few extra values of $k$ in \eqref{XX}. While
it is clear that all $c(i,Y)$ should be rational expressions in $q$,
we expected them to be in fact Laurent polynomials. So we were a bit
startled by the denominators $1+q^2$. However, according to p. 343
top, the Hecke algebra may degenerate at roots of unity $q$, which
justifies at least cyclotomic polynomials as denominators.

\begin{table}[ptb]
\captionwidth\vsize\relax
\newpage
\vbox to \textheight{\vfil
\rottab{ 
\hbox to \vsize{%\hss%\footnotesize
\hss
% \setbox\@tempboxa=\hbox{%
% \textwidth\vsize\relax
$\displaystyle
\begin{array}{|c||*{11}{c|}}
\hline\ry{1.0cm}%
Y & \Young{1,1,1,1,1,1} & % 1
\Young{2,1,1,1,1} &       % 2
\Young{2,2,1,1} &         % 3
\Young{2,2,2} &           % 4
\Young{3,1,1,1} &         % 5
\Young{3,2,1} &           % 6
\Young{3,3} &             % 7
\Young{4,1,1} &           % 8
\Young{4,2} &             % 9
\Young{5,1} &             % 10
\Young{6} \\[6mm]         % 11
\hline\ry{1.6em}%
%  PARTITION #7 & #8 ARE SWOPPED FROM MY NOTES
%      #1 #2  #3  #4  #5  #6  #8  #7  #9  #10 #11
1 & 1 & -\left( \frac{q^4}{1 + q^2} \right)  & \frac{q^7}
   {\left( 1 + q^2 \right) \,\left( 1 + q + q^2 \right) } & 0 & 0 & 0 & 0 & \
0 & 0 & 0 & 0 \\
q^4  & 0 & 0 & 0 & 0 & 0 & 0 & 0 & 0 & \frac{q}
   {\left( 1 + q^2 \right) \,\left( 1 + q + q^2 \right) } & -\left( \frac{q^2}{1 + q^2} \right)  & q^4 \\
q    & 0 & \frac{1}{2} & \frac{-\left( 
       \left( -1 + q \right) \,q \right) }{2} & q + \frac{q^3}{2} & \frac{q^2 + 2\,q^4}{-2 - 2\,q^2} & \frac{q\,\left( 1 - 4\,q + q^2 \right) }
    {2} - c(6,\Young[1.5mm]{ 3,2,1} ) & 0 & \frac{q^6}{2 + 2\,q^2} & \frac{q^6}
   {1 + q + q^2} & 0 & 0 \\[2mm]
-q &   0 & \frac{1}{2} & \frac{-\left( 
       \left( -1 + q \right) \,q \right) }{2} & \frac{q^3}{2} & \frac{q^2 + 
     2\,q^4}{-2 - 2\,q^2} & \frac{{\left( 1 - q \right) }^2\,q}{2} - 
   \frac{q^2}{1 + q} - c(6,\Young[1.5mm]{ 3,2,1} ) & 0 & \frac{q^6}
   {2 + 2\,q^2} & 0 & 0 & 0 \\
-q^2 &   0 & \frac{1}{1 + q^2} & \frac{q\,
     \left( 1 + \left( -1 + q \right) \,q \right) }{1 + q^2} & -q^2 & 
-q^2 & q\,\left( 1 + \left( -1 + q \right) \,q \right)  & -q^2 & -q^2 & \frac{q^3\,\left( 1 - q + q^2 \right) }{1 + q^2} & \frac{q^6}
   {1 + q^2} & 0 \\[2mm]
q^3 &   0 & 0 & \frac{1}{1 + q + q^2} & 0 & \frac{1}
   {2 + 2\,q^2} & c(6,\Young[1.5mm]{ 3,2,1} ) & \frac{q}{2} + q^3 & \frac{-\left( q^2\,
       \left( 2 + q^2 \right)  \right) }{2\,\left( 1 + q^2 \right) } & \frac{\left( -1 + q \right) \,q^2}{2} & \frac{q^4}
   {2} & 0 \\[3mm]
-q^3 &   0 & 0 & 0 & 0 & \frac{1}{2 + 2\,q^2} & \frac{q^2}{1 + q} + 
   c(6,\Young[1.5mm]{ 3,2,1} ) & \frac{q}{2} & \frac{-\left( q^2\,
       \left( 2 + q^2 \right)  \right) }{2\,\left( 1 + q^2 \right) } & \frac{\left( -1 + q \right) \,q^2}{2} & \frac{q^4}{2} & 0 \\[2mm]
\hline
\end{array}%
$%
% }%
% \@tempdima0.3\vsize\relax\advance\@tempdima by -0.3\wd\@tempboxa\relax
% \kern\@tempdima
% \hss
% \copy\@tempboxa
\hss}%
}{This table shows the values $c(i,Y)$ of \eqref{XX} for Young tableaux
$Y$ of 6 boxes and the 7 possible eigenvalues $\dl_i$ of $\pi_Y(\ap)$.
\label{tab2}}{table}
\vss}
\newpage
\end{table}

For $Y=\Young{2,2,2}$, $\pi_Y$ was written down directly on p. 362,
and for its transposed (dual) Young diagram $\Young{3,3}$, one
uses Note 4.6. 

There remain 3 representations, for $Y=\Young{3,2,1}$,
$\Young{4,2}$, and its dual $\Young{2,2,1,1}$, with a total of
15 unknown $c(i,Y)$. These were more complicated to find, since
we knew of no way to evaluate $\tr \pi_Y$ directly. To help
ourselves, first observe that we have, for fixed $Y$, the trace
identities
\begin{eqn}\label{tri}
\sum_{i=1}^7\,c(i,Y)\,=\,\tr \pi_Y([4354])\,=\,
\tr \pi_Y([2132])\,=\,\sum_{i=1}^7 M(i,Y)\dl_i\,,
\end{eqn}
which again give a linear condition on the $c(i,Y)$. To find
further identities, we used \eqref{STR2} in a ``backward'' manner.
We calculated for small\footnote{In the parametrization $P(v,z)$,
used by Morton-Short, unlike for $X$, the coefficients of
$P(\hat\bt_{k,l})$ become quickly large and produce machine
size integer overflows. In particular, we could not calculate
correctly polynomials for $|l|>2$.} $k,l$ ($k$ odd, $|k|\le 5$,
$|l|\le 1$) the polynomial on the left using Morton-Short's
program \cite{MorSho2}. We substituted the known $c(i,Y)$ on
the right of \eqref{STR2}, obtaining thus linear conditions
for the yet unknown $c(i,Y)$. To determine the coefficients,
it remains to understand the effect on $\tr\pi_Y$ of multiplying
with the full twist $\Dl_6^2$. This was, however, done also by
Jones in \S 9, lemma 9.3:
\begin{eqn}\label{l93}
\pi_Y(\Dl_n^2)\,=\,q^{rn(n-1)/d}\,Id_{\pi_Y}\,.
\end{eqn}
For given $Y$, the number $r=r_Y$ is calculated as in lemma 9.1
from figure 3.3, and $d=d_Y=\dim \pi_Y$ more easily by the hook
length formula on p. 341. Call $e_Y=n(n-1)r/d$ (with $n=6$) the
exponent\footnote{This is not to be confused with the variable
$e$, which we use for exponent sum of a braid, or with the
idempotent $e_i$ from Jones's lemma 9.1.} of $q$ on the right
hand-side of \eqref{l93}. These values are given in table
\ref{tab1}. Two simple checks are $\sum\limits_{Y\vdash n}d_Y^2=
n!$ (because the multiplicity of each irrep in the Hecke algebra
equals its dimension), and that $e_Y$ are integers (Remark on p. 358
bottom) and satisfy $e_{Y}+e_{\ol{Y}}=n(n-1)$ for all $Y$
(because of Note 4.6). So \eqref{STR2} gets
\begin{eqn}\label{STR3}
-
(1-q)\cdot\dots\cdot(1-q^6)\,X(\hat\bt_{k,l})(q,\lm)\,=\,
\sqrt{\lm}^{e-5}\,\sum_{Y\vdash 6}\,\wh W_Y(q,\lm)\,
q^{l\cdot e_Y}\,\sum_{i=1}^7\,c(i,Y)\,\dl_i^k\,,
\end{eqn}
with $e=[\bt_{k,l}]=4+4k+30l$.

Actually, each polynomial
$X(\hat\bt_{k,l})$ gives 6 equations for $c(i,Y)$, because
there are 6 relevant $\lm$-coefficients on both hand-sides
of \eqref{STR2} (in degrees $\ffrac{e-5}{2}$, $\ffrac{e-3}{2}$,
\dots, $\ffrac{e+5}{2}$\,; it is helpful to multiply again by
$(1+q^2)$ to get disposed of the denominators of the known
$c(i,Y)$). We have with 18 polynomials 111 equations (6 equations per
polynomial plus the three relevant trace equations \eqref{tri}).
Still the resulting system was too hard to solve by computer, using 
MATHEMATI\-CA\TM{} \cite{Wolfram}, since its coefficients are
(Laurent) polynomials in $q$, with dozens of terms each.

However, substituting some (rational) values of $q$, the system
can be solved immediately. We used this to check first the rank of
the matrix (i.e. which equations are linearly redundant). Again we
were surprised that for 15 variables $c(i,Y)$ the rank was only
$14$. This, however, can be explained from our restraint to odd
$k$ (which we chose for some, purely technical,
component number concerns in the calculation
with Morton-Short's program). Whenever two opposite eigenvalues
$\dl_i$ and $\dl_{i'}=-\dl_i$ occur in $\pi_Y(\ap)$, the equations
\eqref{STR2} for odd $k$ can detect only $c(i,Y)-c(i',Y)$. The
trace equations \eqref{tri}, which involve $c(i,Y)$ and $c(i',Y)$
with the same sign, remedy the shortcoming for $\Young{4,2}$ and
its transposed diagram, but for $\Young{3,2,1}$ we have two pairs
$(i,i')$ of opposite eigenvalues, so we still lose one dimension.
This is, however, not really a problem, because in the
braids $\gm_a=\bt_{6a-1,-2a}$ of \eqref{gma} we
need for our proof, $k=6a-1$ is always odd, so we need only
$c(i,Y)-c(i',Y)$ to evaluate the polynomial in \eqref{STR2}.

We used the special evaluations to select equations that give
the full matrix rank, and to guess the formula for general $q$
for some of the $c(i,Y)$. (The ones we already found suggest that
these formulas should not be so complicated.) Substituting these
(yet potential) solutions too, gives an even simpler linear system
for the still unknown $c(i,Y)$, which then could be solved in $q$.

Since we know that our matrix has rank 14, it is enough to check
our solution (up to the 1-dimentional ambiguity, which will
disappear in \eqref{STR2}) with the 111 equations we have. The
result was confirmed, and is shown in table \ref{tab2}.

With all $c(i,Y)$ determined, the main work is done.
So far we can evaluate $X(\hat\bt_{k,l})$ for odd $k$.
Multiplying \eqref{STR3} by $Z=(1+q)(1+q^2)(1+q+q^2)$ to clear
all denominators in the $c(i,Y)$, normalizing and taking
coefficients in $\lm^m$ for $m=0,\dots,5$, we have
(with $e_Y$ being the exponents of $q$ in \eqref{l93} and
$e=[\bt_{k,l}]$)
\begin{eqn}\label{OPI}
\left[-\frac{1}{\sqrt{\lm}^{e-5}}\,\cdot\,Z\,\cdot\,
(1-q)\cdot\dots\cdot(1-q^6)\,
X_{\hat\bt_{k,l}}(q,\lm)\,\right]_{\lm^m}\,=\,
\sum_{Y\vdash 6}\,\sum_{i=1}^7\,
(Z\cdot c(i,Y))\,\dl_i^k\,
q^{l\cdot e_Y}\,\bigl[\wh W_Y(q,\lm)\,\bigr]_{\lm^m}\,.
\end{eqn}

Recall that for our proof it is enough to show that this term
becomes non-zero for $m=0$ or $m=1$, when $k=6a-1$ and $l=-2a$
for an integer $a>1$ (and $\bt_{k,l}=\gm_a$ in \eqref{gma}).
Now, for odd $k$ (and fixed $Y$), we can group
the sum over $7$ terms $Z\cdot c(i,Y)$ into 5 terms $\tl c(i,Y)$
accounting for $\dl_3=-\dl_4$ and $\dl_6=-\dl_7$, and thus
excluding $i=4,7$. So the above sum in \eqref{OPI} becomes
\begin{eqn}\label{yy}
\sum_{Y\vdash 6}\,\sum_{\scbox{\shortstack{$i=1$\\$i\ne 4$}}}^6
\tl c(i,Y)\dl_i^k\,q^{l\cdot e_Y}\,\bigl[\wh W_Y(q,\lm)\,\bigr]
_{\lm^m}\,.
\end{eqn}
(Now in the $\tl c(i,Y)$, the 1-degree ambiguity of $c(i,Y)$ for
$Y=\Young{3,2,1}$, as explained above, cancels out.)

Among the 55 possible $(i,Y)$ (with $i\ne 4,7$), only 21
of the $\tl c(i,Y)$ are non-zero. It turns out that for $m=0$,
when $k=6a-1$, $l=-2a$ and $a>1$, there is a unique term among
the 21 summands in \eqref{yy} whose minimal degree in $q$ is the
smallest (it is $14-36a$). Thus $[P(\hat\bt_{k,l})]_{v^{e-5}}\ne 0$.
We calculated the polynomial with Morton-Short's program for $a=2$
(where the calculation was still feasible), and it confirmed that
all 6 $v$-terms appear. We also calulated that for $m=0$ and
$a=0,1$ there are two terms of smallest minimal degree. This
better ought to be so, because in that case $\hat\gm_{a}$ are just
a 2-cable of the unknot and $!5_1$, resp., and the coefficients
for $m=0,1$ in \eqref{yy} must be $0$, which we checked once
more separately. (Also we found that for $m=1$ there are two
terms of smallest minimal degree in \eqref{yy} for all $a\ge 0$.)
With this the proof of theorem \reference{thJC} is complete.

Let us finally say that the computer part of the calculation owed
a lot to the use of MATHEMATI\-CA. While, if done properly, it could
be carried out in a few minutes, it required a week of work to
find the way of skillfully programming MATHEMATI\-CA
to do all the separate steps in an efficient way. \qed

\begin{rem}
Note that one could handle the cases $\mu=\sg_1^{6a+1}\sg_2^{-1}$
from the beginning of the proof also using the representation
theoretic argument, by looking at $m=4,5$ for $a<0$. We waived on
this investigation, though, since the proof for $a>0$ was laborious
enough.
\end{rem}

\subsection{Unitarity of the Burau representation\label{SUn}}

In the following $q$ and $t$ are unit norm complex numbers.
We define $\arg (e^{is}):=s\bmod 2\pi$ for $s\in\bR$.
We continue using the formalism of Young tableaux,
the representations $\pi_Y$ and the notations $X(q,\lm)$
and $W_Y$ of the proof of theorem \reference{thJC}.

For the Alexander polynomial $\Dl$ (we resume the
notational convention from before \S\ref{SJS}), as well as
for $3$- and $4$-braids, the
representations $\pi_Y$ are given by Burau representations. We note
(again, and more explicitly) the following descriptions of $\pi_Y$
in terms of the Burau representation $\psi_n$ given in \cite{Jones}.
Again the indexing is chosen so that $\psi_{n-1}$ is the reduced
$(n-1)$-dimensional representation of $B_n$, and by $-\rho$ we denote
the direct product of $\rho$ with the parity representation.

\Youngunitlength1.5mm
As before, $e$ stands for the exponent sum of a braid $\bt$.
For $3$-braids we have the following properties (with
reference to the explanation in \cite{Jones}):
\begin{enumerate}
\item $\pi_{\Young{2,1}}(\bt)(q)=(-1)^e\psi_2(\bt)(q)=
q^e\psi_2(\bt^{-1})(q)$ (because of
row-column symmetry; see Note 4.6).
\item $\pi_{\Young{1,1,1}}(\bt)=(-1)^e$ and
$\pi_{\Young{3}}(\bt)=q^e$ (Note 4.7)
\end{enumerate}

For $4$-braids we have:
\begin{enumerate}
\item $\pi_{\Young{2,1,1}}=-\psi_3$. So
$\pi_{\Young{2,1,1}}(\bt)=(-1)^e\psi_3(\bt)$. (Note 5.7)
\item $\pi_{\Young{3,1}}(\bt)(q)=q^e\psi_3(\bt^{-1})(q)$ (because of
row-column symmetry; see Note 4.6). Also 
\[
\pi_{\Young{3,1}}=-\pi_{\Young{2,1,1}}\wedge
\pi_{\Young{2,1,1}}=-\psi_3\wedge \psi_3\,,
\]
where wedge denotes antisymmetric product (see p.354 bottom).
\item $\pi_{\Young{2,2}}(\bt)=(-1)^e\psi_2(\bar\bt)$, where bar denotes
the homomorphism from $B_4$ to $B_3$ given by $\bar\sg_{1,2,3}=\sg_{1,
2,1}$. (p. 355)
\item $\pi_{\Young{1,1,1,1}}(\bt)=(-1)^e$ and
$\pi_{\Young{4}}(\bt)=q^e$ (Note 4.7)
\end{enumerate}

Now Squier observes in \cite{Squier}, that $\psi_i(\bt^{-1})(q)$
and $\psi_i(\bt)(q^{-1})$ are conjugate and so have the same
trace. So by the self-symmetry of $\pi_{\Young{2,1}}$ we have
$\tr \pi_{\Young{2,1}}(t), \tr \pi_{\Young{2,2}}(t)\in (-t)^{e/2}\bR$.
Similarly $(-t)^{-e/2}\tr \pi_{\Young{2,1,1}}(t)$ and $(-t)^{-e/2}\tr
\pi_{\Young{3,1}}(t)$ are conjugate complex numbers. 
These properties will be important below.

\begin{rem}
Squier uses a different convention for
$\psi_i$ from Jones. He transposes and changes sign in
matrix entries with odd row-column sum (i.e., conjugates by
$\dig(1,-1,1,-1,\dots)$). This, however, does not affect
our arguments. 
\end{rem}

The key point in arguments below is Squier's result.
We write $M^*$ for the conjugate transposed of a matrix $M$.
(That is, $M^*_{i,j}=\ol{M}_{j,i}$.)

\begin{theorem}(Squier \cite{Squier})
For any $n\ge 1$ there exists a Hermitian matrix $J=J^{[n+1]}$ and a
regular matrix $M$, such that with $J_0=M^*JM$ we have $\psi_n^*J_0
\psi_n=J_0$.
\end{theorem}

In particular, $J$ is degenerate or definite iff $J_0$ is so.
Moreover,
\[
J_{i,j}\,=\,\left\{\,\begin{array}{c@{\quad}l}
-1 & \mbox {if\ $|i-j|=1$} \\
\sqrt{t}+1/\sqrt{t} & \mbox {if\ $i=j$} \\
0 & \mbox {otherwise}
\end{array}\right.\,.
\]

It is easy to see that if $t=1$, then $J$ is positive definite.
Now definiteness is an open condition, so for $t$ close to $1$,
it is still valid. 
% By using the criterion of positivity of the principle minors,
One can determine when $J$ loses this
property. % in simple cases. It is easy to verify that for $3$-braids
% this occurs in $e^{\pm 2\pi i/3}$, while for $4$-braids in $\pm i$.
% In general we have

\begin{prop}\label{p4.1}
The Squier form $J^{[n]}$ on $B_n$ degenerates exactly in the $n$-th
roots of unity. In particular, it is positive definite exactly
when $|\arg t|<2\pi /n$.
\end{prop}

\proof Denote by $J^{[n]}$ the form corresponding to $n$-braids,
i.e. the one given by restricting $J$ to the first $n-1$
rows and columns. It is not too hard to calculate the determinant
of $J^{[n]}$. By development in the last row,
\[
\det J^{[n]}\,=\,(\sqrt{t}+1/\sqrt{t})\det J^{[n-1]}-\det J^{[n-2]}\,,
\]
whence
\[
\det J^{[n]}\,=\,\frac{t^n-1}{\sqrt{t}-1/\sqrt{t}}\,\cdot\,
\frac{1}{(\sqrt{t})^n}\,.
\]
Then the claim follows easily. To see definiteness use the
(positivity of) the principal minor criterion. Since $e^{\pm 2\pi i/n}$
is a simple zero of $\det J^{[n]}$, the determinant must become
negative for $|\arg t|\in(2\pi/n,4\pi/n)$. Then applying this
argument to all $n'<n$ shows that $J^{[n]}$ is not positive definite
for $|\arg t|\ge 2\pi/n$. \qed

So on the arcs of $S^1$ that connect the primitive $n$-th root of
unity to $1$, we have that $J$ is positive definite. Now
if $J$ is such, it can be written as $Q^*Q$, and then
conjugating $\psi_i$ by $QM$ we obtain a $U(n-1)$-representation.
This means in particular that all eigenvalues of $\psi_i$
have unit norm. We will below derive implications of this
circumstance for the link polynomials.

\subsection{Norm estimates\label{s4.2}}

The Jones polynomial $V$ can be specified, for our
purposes, by $V(t)=X(t,t)$. In the following, which
root of complex numbers is taken is irrelevant,
important is though that it be kept fixed in subsequent
calculations. By $\Re$ we denote the real part of
a complex number.

\begin{theorem}\label{th1}
If $|t|=1$, $\Re t>0$ and $\bt$ is a $4$-braid, then 
$\big |V_{\hat\bt}(t)\big|\,\le\,(2\Re \sqrt{t})^3$.
If $\bt$ is a $3$-braid and $\Re t>-1/2$, then
$\big |V_{\hat\bt}(t)\big|\,\le\,(2\Re \sqrt{t})^2$.
\end{theorem}

\proof We have from \cite{Jones} that if $\bt\in B_4$ with
$[\bt]=e$, then
\begin{eqn}\label{Vf}
V_{\hat\bt}(t)\,=\,\left(-\sqrt{t}\right)^{e-3}
\left[\,\frac{t(1-t^3)}{1-t^2}\,\tr \psi_3+\frac{t^2}{1+t}
\tr \bar\psi_2\,+\,\frac{1-t^5}{1-t^2}\,\right]\,,
\end{eqn}
where $\bar\psi_2$ is the composition of $\psi_2$ with
$\bar{ }\,:\,B_4\to B_3$. Taking norms and using that
$\bar\psi_2$ and $\psi_3$ are unitary, we find
\begin{eqn}\label{X}
\big |V(t)\big|\,\le 3\left|\frac{1-t^3}{1-t^2}\right|+
\frac{2}{|1+t|}+\left|\frac{1-t^5}{1-t^2}\right|\,.
\end{eqn}
It is now a routine (but somewhat tedious) calculation to
verify that the r.h.s. is equal to $(2\Re \sqrt{t})^3$
for $|t|=1$, $\Re t>0$. 

For $\bt\in B_3$ we have similarly
\begin{eqn}\label{V_3}
V_{\hat\bt}(t)\,=\,\left(-\sqrt{t}\right)^{e-2}
\left[\,t\cdot \tr\psi_2\,+\,(1+t^2)\,\right]\,,
\end{eqn}
and the result follows using $|\tr\psi_2|\le 2$. \qed

This theorem generalizes Jones's result \cite[proposition 15.3]{Jones}
for $n\le 4$, where he considers $t=e^{2\pi i/k}$, $k\ge 5$. In fact
the comparison to (and established coincidence with) Jones's
estimate led to the simplification of the r.h.s. of \eqref{X}.
In \cite{posex_bcr} we noted that Jones's estimate can be
better than MWF when $MWF=3$, but by connected sum one
can give an example for $MWF=4$. Again it appears that for
$\bt\in B_4$ and $t\ne e^{\pm \pi i/3}$ the set $\{\,|V_{\hat
\bt}(t)|\,\}$ is dense in $[0,(2\Re \sqrt{t})^3]$, and similarly it
is in $[0,(2\Re \sqrt{t})^2]$ for $\bt\in B_3$ and $t\ne e^{\pm \pi
i/3},e^{\pm \pi i/5}$. (See also the remarks at the end of \S 12
in \cite{Jones}.)

\begin{conjecture}
If $\bt\in B_n$ and $|t|=1$ with $|\arg t|<2\pi /n$ then
$|V(t)|\,\le\,(2\Re \sqrt{t})^{n-1}$.
\end{conjecture}

In the case of the Alexander polynomial $\Dl(t)=X(t,1/t)$, we
can say something on general braids.

\begin{theorem}
For each $n\ge 2$, %there is an $\eps_n>0$ such that
if $|t|=1$ and $|\arg t|\le 2\pi/n$ and $\bt\in B_n$ then
\[
\big |\Dl_{\hat\bt}(t)\big|\,\le\,\frac{2^{n-1}|1-t|}{|1-t^n|}\,.
\]
\end{theorem}

\proof
% Let $\eps_n>0$ be so that $\ds\frac{|1-t^n|}{|1-t|}>n-0.1$ and
$J_0$ is positive definite when $|t|=1$ and $|\arg t|<2\pi/n$.
Then
\[
\big |\Dl_{\hat\bt}(t)\big|\,\frac{|1-t|}{|1-t^n|}\,=\,
\big |\det (1-\psi_{n-1}(\bt))\big|\,,
\]
and all eigenvalues of $1-\psi_{n-1}$ have norm $\le 2$.
The case $|\arg t|=2\pi/n$ follows by continuity. \qed

\begin{corr}
For all $n,k$ the set
\[
\big\{\,\Dl(\hat\bt)\,:\,\bt\in B_n\,,\ \ \deg\Dl\le k\,\big\}
\]
is finite. That is, among closed braids of given number of strands
only finitely many Alexander polynomials of given degree occur.
\end{corr}

\proof $\Dl$ is determined by $\Dl(t_i)$ for $k$ different
$t_i$ with $|t_i|=1$, $0<\arg t_i<2\pi/n$, by means of a 
linear transformation using the (regular Vandermonde)
matrix $M=(t_i^j)_{i,j=1}^n$.
So $||[\Dl(t)]_{t^j}||\,\le\,||M^{-1}||\cdot ||\Dl(t_i)||$\,.
\qed

This result should be put in contrast to the various constructions
of knots with any given Alexander polynomial. For example
a recent construction of Nakamura \cite{Nakamura2} allows to
realize the degree of the polynomial by the (actually braidzel)
genus of the knot. (We were subsequently independently
able to further specialize this result to canonical genus.) 
% do not know the following 
% 
% \begin{question}
% Is any Alexander polynomial admitted by a knot with canonical
% genus equal to the degree of the polynomial?
% \end{question}
% 
A different construction of Fujii 
\cite{Fujii} shows that knots with 3 bridges admit all
Alexander polynomials. So the situation between braid and
plat closures is completely different.

Compare also Birman-Menasco's result in \cite{BirMen3}, mentioned
(for knots) in remark \ref{rfn}, that there are only finitely
many closed braids of given number of strands with given genus. 
Note that we do not claim that only finitely many 
closed braids of given number of strands with given Alexander
polynomial (degree) occur. For $3$-braids it is true, but
from $5$-braids on the non-faithfulness of the Burau
representation should (in principle, modulo the evaluability
of another invariant) make it possible to construct
infinite families of links with the same (for example, trivial)
polynomial. It makes some sense to ask about the status of 4-braids.

\begin{question}
Are there only finitely many 
closed 4-braids of given Alexander polynomial (degree)?
\end{question}

When working with $\Dl$, for $3$- and $4$-braids we can be
more explicit.

\begin{corr}\label{pp1}
If $\bt\in B_4$, then
\begin{eqn}\label{Xx}
\big|\,\Dl_{\hat\bt}(t)\,\big|\,\le\,\frac{8\,|1-t|}{|1-t^4|}\,,
\end{eqn}
when $|t|=1$ and $\Re t>0$. If $\bt\in B_3$ and $\Re t>-1/2$, then
\begin{myeqn}{\qed}
\big|\,\Dl_{\hat\bt}(t)\,\big|\,\le\,\frac{4\,|1-t|}{|1-t^3|}\,.
\end{myeqn}
\end{corr}

Putting $t=e^{2\pi i/5}$ in \eqref{Xx} we have $\big|\,\Dl(e^{2\pi
i/5})\,\big|\,\le\,8\,.$ This improves the bound $10.47\dots$ in
\cite{posex_bcr} suggested to replace Jones's (incorrect) value 6.5
in \cite{Jones2}. For 3-braid \em{knots} Jones gives in
\cite[proposition 15.2]{Jones} the better bound 3 when $t=i$, using
the property $V(i)=\pm 1$.

\begin{exam}\label{4249}
The simplest knots with $MWF=4$ which can be excluded from being a
4-braid using corollary \reference{pp1} are $13_{8385}$ (where one
can use $t=e^{2\pi i/9}$) and $14_{37492}$ (with $t=e^{2\pi i/8}$).
For $3$-braids we can deal with the known examples $9_{42}$
and $9_{49}$. (This gives now an alternative proof that
there are no 3-braid knots with such Alexander polynomial, as
we explained in example \reference{x29}.)
\end{exam}

If we know the exponent sum of $\bt$ we can do better.

\begin{prop}\label{pp2}
If $\bt\in B_4$, $e=[\bt]$ and $t$ as before, then
\[
\left|\left(
-\sqrt{t}\right)^{e-3}\Dl_{\hat\bt}(t)\frac{1-t^4}{1-t}-1
+(-t)^e\,\right|\,\le 6.
\]
If $\bt\in B_3$ and $\Re t>-1/2$, then
\begin{myeqn}{\qed}
\left|\left(
-\sqrt{t}\right)^{e-2}\Dl_{\hat\bt}(t)\frac{1-t^3}{1-t}-1
-(-t)^e\,\right|\,\le 2\,.
\end{myeqn}
\end{prop}

\begin{exam}\label{lxd}
Since one can determine the possible $e$ via MWF from $P$,
one can apply proposition \reference{pp2} for given $P$. In
the $3$-braid case we can exclude $10_{150}$ this way. The
remaining two 10 crossing knots with unsharp MWF, $10_{132}$
and $10_{156}$, fail~-- understandably, since they share the
skein polynomials of (the closed 3-braids) $5_1$ and $8_{16}$
resp. (See the table\footnote{The second duplication was
noted in the remarks after Jones's table, but not referred
to correctly in its last column.} in \cite{Jones}.) For
$4$-braids several new 14 crossing knots can be ruled
out, for example $14_{21199}$.
\end{exam}

\begin{rem}
The use of $P$ to restrict the possible values of $e$ is
usually most effective, but not indispensable. There are
other conditions on $e$, originating from Bennequin's
work \cite{Bennequin}, that can be more applicable in
certain cases where the calculation of $P$ is tedious.
Also, when $t$ is a root of unity of order $n$, the tests
depend only on $e\bmod 2n$.
\end{rem}

\begin{rem}
Note that the quantity $\chi$ in \eqref{chi} is equal to $u\,
\tr \psi_2$, where $u=\sqrt{-t}+\ds\frac{1}{\sqrt{-t}}$ and
$\psi_2$ is the Burau matrix of $B_3$ (see the explanation in
\cite{Kanenobu}). Thus one can obtain
similar estimates for values of $Q$ on 3-braids.
\end{rem}

\subsection{Skein polynomial}

Now it is natural to look at the full 2-variable skein polynomial
$X$. We have, as in \eqref{STR}, for $\bt\in B_4$ of exponent
sum $e$,
\[
X_{\hat\bt}(q,\lm)\,=\,-\sqrt{\lm}^{e-3}\,
\sum_{Y\vdash 4}\,\tr\pi_Y(q)\,\tW_Y(q,\lm)\,,
\]
where the weights $\tW_Y$ are given in $\lm$-coefficients by
table \reference{table1}. (They are all polynomials in $\lm$ of
degree $3$, with coefficients being rational expressions in
$q$.) Now, with given $e$ and $P$, we have $4$ equations
(the coefficients in $\lm$) in $3$ unknowns (the traces of
$\pi_{211}$, $\pi_{31}$ and $\pi_{22}$; we use here the
partition notation for the subscripts). However, the
restriction of the matrix in table \reference{table1} to
the columns of $\pi_{211}$, $\pi_{31}$ and $\pi_{22}$
has rank $2$. This means that two of the $X_i=[X]_{\lm^
{i+(e-3)/2}}$ for $0\le i\le 3$ determine the other two.
One could believe now to use
this as a 4-braid test. However, these two relations
result from the general substitutions $\lm=1$ and $\lm=1/q^2$
that turn $X$ into the component parity count or $1$.
These substitutions kill all trace weights except of the trivial
or parity representation, and for these representations the weights
become also independent on the braid group. Thus the
relations between the $X_i$ will hold whenever $MWF\le 4$,
and are useless as a 4-braid test. In a similar vein, one has 

\newbox\@tempboxb
\setbox\@tempboxb=\hbox{
\capt The weights of the traces contributing to the
$\lm$-coefficients of $X$ for a $4$-braid. Each table entry
must be multiplied with the factor on the right in the
first row to obtain the contribution of the corresponding trace
to $-\sqrt{\lm}^{e-3}$ times the power of $-\lm$ in the
first column. In the first row the symbol $[i]$ denotes $1-q^i$.}

\begin{table}
\[
\begin{array}{c*5{|c}}
& 
\ds \tW_{1111}=\frac{1}{[2][3][4]}\times &
\ds \tW_{211}=\frac{q}{[2][1][4]}\times &
\ds \tW_{22}=\frac{q^2}{[2][2][3]}\times &
\ds \tW_{31}=\frac{q^3}{[2][1][4]}\times &
\ds \tW_{4}=\frac{q^6}{[2][3][4]}\times \\[5mm]
\hline
\ry{6mm}
\ds 1 & 1 & 1 & 1 & 1 & 1 \\[1mm]
-\lm  & q^2+q^3+q^4 & 1+q^2+q^3   & 1+q+q^2   & 1+1/q+q^2 & 1+1/q+1/q^2 \\[2mm]
\lm^2 & q^5+q^6+q^7 & q^3+q^2+q^5 & q+q^2+q^3 & 1/q+q+q^2 & 1/q+1/q^2+1/q^3 \\[2mm]
-\lm^3& q^9         & q^5         & q^3       & q         & 1/q^3 \\[2mm]
\end{array}
\]
\caption{\unhbox\@tempboxb
\label{table1}}
\end{table}

\begin{prop}\label{pp7.6}
If a braid $\bt$ has $MWF(\bt)\le 4$, then $V(\hat\bt)$ and
$\Dl(\hat\bt)$ together with the exponent sum $[\bt]=e$,
determine $P(\hat\bt)$.
\end{prop}

\proof From $e$, $V$ and $\Dl$ we know
\[
\tl X_a\,=\,\sum_{i=0}^3\,X_i(q)\,q^{a(i+(e-3)/2)}\,,
\]
for $a\in\{-2,-1,0,1\}$. (For $a=-2$ and $a=0$ we have the
trivializing substitutions, $a=-1$ corresponds to $\Dl$ and
$a=1$ to $V$.) So one can recover $X_i$ from $\tl X_a$
and $e$. (One can write down an explicit formula easily.) \qed

This condition is thus equally unhelpful as a 4-braid
test. For a similar reason, I expect (though I have not
rigorously derived) an explanation of the (experimentally
observed) failure of Jones's conditions \cite[\S 8]{Jones}
(see formula (8.10)) to obstruct to a 4-braid.

\begin{exam}
The knot $11_{386}$, known from \cite{LickMil}, has the
Jones polynomial of the figure-8-knot. So $11_{386}$ and its mirror
image show that the dependence of $P$ on $e$ in proposition
\reference{pp7.6} is essential.
\end{exam}

\begin{question}
Are there 5-braids $\bt_{1,2}$ (at least one of which has 
$MWF=5$), with the same exponent sum,
$V(\hat\bt)$ and $\Dl(\hat\bt)$, but different $P(\hat\bt)$?
\end{question}

The lack of such examples, after some check in the knot tables,
is at least not fully explainable. Only 6-braids could be found.

\begin{exam}
The knots $16_{443392}$ and $!15_{223693}$ have the same $V$
and $\Dl$ but different $P$ polynomial. They have $MWF=5$ resp. $4$,
with $5$- resp. $4$-braid representations of exponent sum $2$ resp.
$-1$, so one obtains 6-braids of exponent sum $1$ by stabilization.
\end{exam}

One can now, as before, go over to the norms in each
row in table \reference{table1}, or of arbitrary linear
combinations of such rows. Then for $q$ where $J_0$ is
definite we obtain again estimates on $|X_i(q)|$, or of
$|X(q,\lm)|$ for any non-zero complex number $\lm$. (In
particular, for $\lm=q$ we obtain theorem \reference{th1} and
for $\lm=1/q$ corollary \reference{pp1}.) Although they still
contain the Jones and Alexander polynomial conditions (which
were observed both non-trivial in comparison to MWF), such
skein polynomial norm estimates have not proved in practice,
as a 4-braid test, an efficient improvement over their special
cases and MWF. Below we will explain how to do much better.

Knots like $10_{132}$ and $10_{156}$ in example \ref{lxd} show
a disadvantage of our test, resulting from not taking into account
information of other invariants. On the opposite hand, the
mere use of $\Dl$ or $P$ reduces calculation complexity, and
excludes any potential further knot with such polynomial.
In the case of $P$ this gives yet a different way to
answer negatively Birman's question if one can realize any
skein polynomial by a link making MWF sharp (see \cite
{posex_bcr,3br}).

Using the Brandt-Lickorish-Millett-Ho polynomial, Murakami
\cite{Murakami} and later Kanenobu \cite{Kanenobu} gave with
theorem \reference{TMF} a more
efficient (in excluding examples, though less in calculation
complexity) test for a 3-braid. But the work in \cite{3br} and
in the previous sections of this paper makes the study of 
polynomials of 3-braids anyway less relevant. $4$-braids
become much more of interest, and in this case, after MWF,
the problem to find applicable conditions has been largely
unsettled for quite a while.

\subsection{Recovering the Burau trace}

\subsubsection{Conditions on the eigenvalues}

Using just norms clearly weakens the conditions considerably,
and so one would like to identify the Burau eigenvalues directly.
However, the relations between the $X_i$ do not allow to recover 
by simple algebraic means the individual traces from $P$. 

Using Squier's unitarity, there is an analytic way to recover
the Burau trace of 4-braids (at least for generic $t$ and
up to finite indeterminacy). Since by the previous remarks
the use of $P$ is not essential, we will describe the
procedure given $\Dl$, $V$ and $e$. This gives the most
significant practical enhancement to the above 4-braid tests.

In the following we fix a 4-braid $\bt$ of exponent sum $e$, whose
closure link has Jones polynomial $V$ and Alexander polynomial $\Dl$.
We let $t$ be a unit norm complex number with non-negative real part.

We use that
\[
\tr \psi_3\wedge \psi_3 = (-t)^{e/2}\ol{\tr \psi_3(-t)^{-e/2}}\,.
\]
Then we have
\[
\Dl(t)\,=\,\frac{1-t}{1-t^4}\,\left(-\frac{1}{\sqrt{t}}\right)^{e-3}
\cdot\,\left[ 1-\tr \psi_3 + (-t)^{e/2}\ol{\tr \psi_3(-t)^{-e/2}}-
(-t)^e\,\right]\,.
\]
So, with $i=\sqrt{-1}$,
\[
\dl\,:=-\frac{1}{2i}(-t)^{-e/2}
\left[ \Dl(t)\frac{1-t^4}{1-t}\left(-\frac{1}{\sqrt{t}}\right)^{3-e}
-1 + (-t)^e\,\right]
\]
is a real number. Now when $\lm_{1,2,3}$ are the eigenvalues of
$\psi_3(\bt)$, we have
\begin{eqnarray*}
A & = & \lm_1+\lm_2+\lm_3\,=\,(-t)^{e/2}(y+i\dl) \\
B & = & \lm_1\lm_2+\lm_1\lm_3+\lm_2\lm_3\,=\,(-t)^{e/2}(y-i\dl) \\
C & = & \lm_1\lm_2\lm_3\,=\,(-t)^e\,.
\end{eqnarray*}
Here $y$ is a real number we do not know, and we try to determine.
Since $|\lm_k|=1$, we have
as before $|A|\le 3$ and $|B|\le 3$. So the range
for $y$ is $[-y_0,y_0]$ with $y_0=\sqrt{9-\dl^2}$.
(If $|\dl|>3$ we are done as before.) Then
\[
\tr\psi_3\,\in\,[(-t)^{e/2}(-y_0+i\dl),\, (-t)^{e/2}(y_0+i\dl)]\,,
\]
with the interval understood lying in $\bC$.
Let $\psi_-$ and $\psi_+$ be the endpoints of this interval.

Now one can restrict the interval $[-y_0,y_0]$ for $y$ using the
Jones polynomial as follows.

Let $\rho\,:=(-t)^{-e/2}\tr\bar\psi_2\in [-2,2]$. The restriction
to the given range follows because $J^{[3]}$ (see the proof of
proposition \reference{p4.1}) is also definite when
$J^{[4]}$ is. From \eqref{Vf} we have
\[
\tr \psi_3\,=\,\frac{1-t^2}{t(1-t^3)}\,\left[\,
V(t)\,\left(-\frac{1}{\sqrt{t}}\right)^{e-3}-\frac{1-t^5}{1-t^2}
-\frac{t^2}{1+t}(-t)^{e/2}\cdot \rho\,\right]\,=:\,\tl\psi(\rho)\,,
\]
so
\[
\tl\psi(\rho)=\frac{1-t^2}{t(1-t^3)}\left [
V(t)\,\left(\frac{-1}{\sqrt{t}}\right)^{e-3}-
\frac{1-t^5}{1-t^2}\right] -\frac{\rho}{1+2\Re t}(-t)^{e/2}\,.
\]
(Just for the purpose of defining $\tl\psi$, we should regard here
$\rho$ as a formal parameter, rather than as a concrete value.)

Let $\tl\psi(\pm 2)=:\tl\psi_\pm$. Since $1+2\Re t>0$, we have
$\Re((-t)^{-e/2}\tl\psi_+)>\Re((-t)^{-e/2}\tl\psi_-)$.

Then $[\tl\psi_-,\tl\psi_+]\subset \bC$ is an interval
of the same slope as $[\psi_-,\psi_+]$, so we check if they overlap.

Let
\[
\tl y_{\pm}= (-t)^{-e/2}\tl\psi_\pm -i\delta.
\]

Then for a consistent restriction on $\tr \psi_3$ the following holds:
\begin{enumerate}
\item $\tl y_{\pm}$ are real
\item $ \tl y_+\tl y_-\le 0$ or $\min(\tl y_\pm^2+\delta^2)\le 9$.
\end{enumerate}

Potentially these conditions may be violated, but in
practice they seem always to hold. (We have not elaborated on why
this is so, though it may be worth understanding.) Then at least,
we consider 
\begin{eqn}\label{qw0}
y\in [-y_0,y_0]\cap [\tl y_-,\tl y_+]\,.
\end{eqn}

We have now the cubic
\[
x^3 +ax^2 +bx +c \,:=\,x^3-A x^2+B x-C \,=\,0\,.
\]
One solution is obtained by Cardano's formula
\begin{eqn}\label{qw}
\lm_1\,=\,-\frac{a}{3}
-\frac{\sqrt[3]{2}(-a^2+3b)}{3\Gm} + \frac{\Gm}{3\sqrt[3]{2}}\,,
\end{eqn}
where 
\begin{eqn}\label{Gm}
\Gm=\sqrt[3]{
  -2a^3 + 9ab - 27c + 
  \sqrt{
    27\,\left( -a^2b^2 + 4b^3 + 4a^3c -18abc + 27c^2
    \right)
  }
}\,.
\end{eqn}
We must have that $|\lm_1|=1$. Then we must check $|\lm_{2,3}|=1$.
For this their exact determination is not necessary.
We have $\lm_2+\lm_3=A-\lm_1$ and $\lm_2\lm_3=C/\lm_1$.
In order $|\lm_{2,3}|=1$, we must have
\[
\lm_2+\lm_3=c\cdot \sqrt{\lm_2\lm_3}\,,
\]
with $c\in [-2,2]$, so $A-\lm_1=c\sqrt{C/\lm_1}$, which is
equivalent to
\begin{eqn}\label{qw2}
\frac{(A-\lm_1)^2\lm_1}{C}\,\in\,[0,4]\,.
\end{eqn}

\subsubsection{Applications and examples}

The image of the r.h.s. of \eqref{qw0} under \eqref{qw}
will be some curve in $\bC$ that generically intersects
$S^1$ in a finite number of points. This parametric
equality can be examined numerically and allows to recover
$A=\tr \psi_3$ up to finite indeterminacy. In particular,
we have

\begin{prop}\label{4brt}
Assume for some $t\in S^1$ with $|\arg t|\le \pi/2$ we
cannot find $y$ as in \eqref{qw0}, such that $\lm_1$
given by \eqref{qw} has norm 1 and \eqref{qw2} holds. Then
there is no $\bt\in B_4$ with the given $e$ whose closure has
the given Alexander and Jones polynomial. \qed
\end{prop}

\begin{rem}
To apply the test in practice, one chooses a small stepwidth $s$
for $y$ in the interval \eqref{qw0}, and calculates the derivative
of the r.h.s. of \eqref{qw} in $y$ to have an error bound on
$||\lm_1|-1|$ in terms of $s$. When some of the radicands in \eqref{Gm}
becomes close to $0$, some care is needed. One situation where
such degeneracy occurs are the knots $15_{144634}$, $15_{144635}$,
and $15_{145731}$. They have representations $\bt\in B_4$, whose
Burau matrix is trivial for $t=e^{\pi i/3}$. This may be noteworthy
on its own in relation to the problem whether $\psi_3$ is faithful. 
(Clearly $\psi_3$ is unfaithful at any root of unity on the
center of $B_4$, but here $\bt$ is not even a pure braid.)
\end{rem}

\begin{exam}
Applying proposition \reference{4brt} we can exclude $11_{387}$, one
of the 7 prime 11 crossing knots with $b=5$ but $MWF\le 4$. Eleven
prime knots with 12 crossings, and 63 of 13 crossings where braid
index 4 is not prohibited by MWF can be ruled out. The correctness
of these examples was later verified by the 2-cabled MWF. Up to 16
crossings more than 4000 examples were obtained. (Let us note that
from them only about 100 can be identified using the norm estimates.)
\end{exam}

\begin{exam}
The check of prime 14 crossing knots to which our criterion applied,
revealed 6 knots with 2-cable MWF bound 8. One, $14_{22691}$,
can be excluded from braid index 4 (as done with $14_{45759}$ 
in \cite{posex_bcr}), by making the $v$-degrees of the
polynomial of the 2-cable contradict the exponent sum of its
possible 8-braid representation. However, for the other 5
knots, $14_{28220}$, $14_{30960}$, $14_{41334}$, $14_{41703}$,
and $14_{44371}$, the argument fails, and so our condition
seems the only applicable one. (Clearly a 3-cable polynomial
is not a computationally reasonable option, and even the 2-cable
requires up to several hours, while our test lasts a few seconds.)
\end{exam}

Still our criterion leaves open several interesting examples
of (apparent) failure of the 2-cable MWF inequality, among
them the knot $13_{9684}$ encountered with M.~Hirasawa. A more
general, and important, possible application is as follows.

\begin{rem}\label{rJC}
In relation to Jones's conjecture \reference{jcc}, we
already quoted (in the proof of theorem \ref{thJC})
the observation in \cite{posex_bcr} that counterexamples
to the conjecture would make MWF and \em{all} of its cabled
versions unsharp. In that sense, our 4-braid test may be
the first possible approach toward identifying such a
counterexample. Birman-Menasco claimed indeed a family
of 6-string potential counterexamples, and K.~Kawamuro
gave later a simpler family on 5 strings. Extensive
checks with our test of Kawamuro's knots failed to turn up
successful cases. This puzzled us a while, until Kawamuro
reported recently that in fact H.~Matsuda falsified all
Birman-Menasco (and hence also Kawamuro's) candidates.
\end{rem}

\subsection{Mahler measures}

\subsubsection{3-braids}

The material in section \reference{SBurau} originated from
a question of S.~Kamada whether the Alexander polynomial
Mahler measure $M(\Dl)$ is bounded on closed 3-braids. 

\begin{defi}
For a Laurent polynomial $p\in\bZ[t^{\pm 1}]$, and $n\in\bN_+$,
we define the $n$-norm of $p$ by
\[
||p||_n\,:=\,\sqrt[n]{\sum_{k=-\infty}^{\infty}\bigl([p]_k\bigr)^n}\,,
\]
and its Mahler measure by
\[
M(p)\,:=\,\prod_{|t|\ge 1,\,p(t)=0}\,|t|\,.
\]
This is extended to $p\in\bZ[t^{\pm 1/2}]$ in the obvious way.
See \cite{SSW,CK}.
\end{defi}

Kamada's question is related to controlling $|V(t)|$ and $|\Dl(t)|$
for $t\in S^1$, since the 2-norm $||\,.\,||_2$ and the Mahler
measure $M$ of polynomials have circle integral formulas (see
\eqref{55'} and \eqref{Ji}). One
can thus ask whether $||\Dl\cdot W_n||_2$ is bounded for proper
$W_n\in \bZ[t,1/t]$, or (weaker) whether $M(\Dl)$ is bounded
(and similarly $M(V)$ and $||V\cdot W_n||_2$) for braid
index $\le n$.

For values of $t$ where $J_0$ is not definite, however, there
seems little one can say on the range (on closed 3-braids)
of $|V(t)|$ or $|\Dl(t)|$. Likely they are dense in $\bR_+$.
% So bounded 2-norms seem unlikely, but even if $|\Dl(t)|$
% is unbounded for any $t$ with indefinite $J_0$, that is
% not sufficient.
We can conclude boundedness properties in
special cases where the indefinite $J_0$ values of $t$ are
controllable. For example the following can be proved easily.
(Note that for polynomials with integer coefficients the properties
a set of polynomials to have bounded $2$-norms, or finitely many
distinct $2$-norms, are equivalent, and also equivalent to the same
two properties for the $1$-norm.)

\begin{prop}
The set
\[
\{\,||(1-t^{2n})\Dl_{\hat\bt}(t)||_2\,:\,
\bt\in B_{2n},\,\Dl_{\hat\bt}\in\bZ[t^{\pm n}]\,\}
\]
is finite for any $n\ge 2$.
\end{prop}

\proof We have 
\begin{eqn}\label{55'}
||X||_2^2\,=\,\int_{0}^{1}\,\big|\,X(e^{2\pi is})\,\big|^2\,ds\,.
\end{eqn}
If $\Dl\in\bZ[t^{\pm n}]$, then this integral for 
\[
X=\frac{(1-t^{2n})\Dl(t)}{1-t}=\det(1-\psi_{2n-1})
\]
is controlled from its part over $0<s<\ffrac{1}{2n}$, and
this in turn is controlled by the unitarity of $\psi_{2n-1}$.
To eliminate the denominator, observe that the norm in
\eqref{55'} is just the one in $L^2(S^1)$, and so
$||P\cdot Q||_2\le ||P||_2\cdot ||Q||_2$. \qed

However, in general a bound on the Mahler measure of
arbitrary polynomials of closed $n$-braids does not exist
even for $n=3$.

\begin{prop}
The Mahler measure of $\Dl$ and $V$ is unbounded on closed
3-braids.
\end{prop}

{\def\gm{g_{\min}'}
\proof In the following $x\le O(y)$ means $\limsup |x/y|<\infty$
(the limit following from the context), $x\ge O(y)$ means $\liminf
|x/y|>0$ and $x=O(y)$ means both $x\le O(y)$ and $x\ge O(y)$.
The quantities $C$ and $C'$ will stand for real constants, with
$C$ being positive. They are understood to be independent on the
indices of the terms in the formula they occur in, but
may vary between the various inequalities.

We consider the links $L_n$ given by the closed 3-braids
$(\sg_1\sg_2^{-1})^n$. Using the Burau matrix eigenvalues, we find
\[
V_n\,:=V(L_n)\,=\,t+\frac{1}{t}+e_-^n(t)+e_+^n(t)\,,
\]
where 
\[
e_{\pm}(t)\,=\,-\frac{t+1/t-1}{2}\,\pm
\sqrt{\left( \frac{t+1/t-1}{2} \right)^2-1 }\,.
\]
Similarly
\[
\Dl_n\,:=\,\Dl(L_n)\,=\,\frac{1}{t+1/t+1}\left(2-e_-^n(t)-e_+^n(t)
\right)\,.
\]
Now the Jensen integral gives
\begin{eqn}\label{Ji}
\log\,M(V_n)\,=\,\int_0^1\,\log\,|V_n(e^{2\pi is})|\,ds\,.
\end{eqn}
Since $V_n$ is reciprocal, it suffices to consider the integral
over $s\in[0,1/2]$. Herein, with $t = e^{2\pi is}$, for
$s\in (1/3,1/2]$, we have real $e_{\pm}(t)$, one with norm $>1$. So
\[
\int_{1/3}^{1/2}\,\log\,|V_n(e^{2\pi is})|\,ds\,\ge\,O(n)\,.
\]
We must argue about $s\in[0,1/3]$. In that case, $e_{\pm}(t)$
are conjugate complex numbers of unit norm. This implies that
the integrand is bounded above when $n\to\infty$, but our problem
is to show that it does not decrease too quickly with $n$.
The assignment (with $i=\sqrt{-1}$)
\[
t\,\longmapsto\,f(t)\,=\,-\frac{t+1/t-1}{2}\,+
i\sqrt{ 1-\left( \frac{t+1/t-1}{2} \right)^2 }\,
\]
maps $e^{2\pi i[0,1/3]}\to e^{2\pi i[1/6,1/2]}$, so
\[
\frac{1}{2\pi}\,\arg t\,\stackrel{g}{\longmapsto}\,
\frac{1}{2\pi}\,\arg f(t)
\]
is a map $g\,:\,[0,1/3]\to[1/6,1/2]$, which is
bijective and monotonous, and $C^1$ except in $t=1/3$
where $g'\to\infty$. In particular,
\[
\gm\,=\,\min\,\{\,g'(s)\,:\,s\in[0,1/3]\,\}\,>0\,.
\]
Now consider the functions 
\[
h(s)\,:=\,2\cos 2\pi ng(s)\quad\mbox{and}\quad
m(s)\,:=\,-2\cos 2\pi s
\]
for $s\in[0,1/3]$. The singularities of the Jensen integral \eqref{Ji}
correspond to the intersections of the graphs of $h$ and $m$.

If $n\gg \gm$, then all intersection points\footnote{For the rest
of the proof $i$ is used as an index rather than the complex unit.}
$s_i$ of the graphs of $h$ and $m$ have $|m'(s_i)|>|h'(s_i)|$. (Here
we use also that $m\ne h$ when $s=0$.) In particular, between two
critical points $x_i$ of $h$, it intersects $m$ at most once. So the
number $n'$ of intersections $s_i$ satisfies
\[
\left|\,n'-\frac{2}{3}n\,\right|\,\le\,O(1)\,,\mbox{\quad
and therefore\quad}\,n'=O(n)\,.
\]

Moreover, 
\begin{eqn}\label{20.5}
|m'(s_i)|-|h'(s_i)|\,\ge\,(2-|m(s_i)|)\,\cdot C\,,
\end{eqn}
for a constant $C$ independent on $i$ and $n$. Now when $n$
increases, the $s_i$ will concentrate around $s=\frac{1}{3}$,
but since $\gm>0$, for $s$ small, the $x_i$ will be at distance
$\ge O(\frac{1}{n})$. So
\begin{eqn}\label{20.6}
2-|m(x_i)|\ge C\cdot \frac{i^2}{n^2}
\end{eqn}
for a constant $C$ independent on $i$ and $n$. Then
if $s_1,\dots,s_{n'}$ are the $n'$ intersection points
of $h$ and $m$, we can estimate for values of $s$
between the critical points $x_i$ and $x_{i+1}$ of $h$
\def\tsi{\bigl(2-|m(s_i)|\bigr)}
\begin{eqn}\label{x5}
|m(s)-h(s)|\,\ge\,C\,\cdot \tsi\,\cdot |s-s_i|^2
\end{eqn}
for a number $C$ independent on $i$ and $s$. Inequality
\eqref{x5} follows because we can control the behaviour
of $|m(s)-h(s)|$ from \eqref{20.5} for $s$ around $s_i$,
and though $|m-h|$ decreases slightly around the endpoints
of the interval $[x_i,x_{i+1}]$, this decrease can be
controlled by \eqref{20.6}. So
\begin{eqn}\label{x4}
\int_{x_i}^{x_{i+1}}\,\log\,|\,m(s)-h(s)\,|\,ds\,\ge\,
(x_{i+1}-x_i)\Bigl(\,\log\,\tsi\,+\log C\,\Bigr)\,+\,
\,2\,\int_{x_i}^{x_{i+1}}\,\log |s-s_i|\,ds\,.
\end{eqn}
Since for small $i$ (when $2-|m(s_i)|$ is
small), the $x_i$ will be at distance $\ge O(\frac{1}{n})$, using
\eqref{20.6} and $m(s)\le 1$ when $0\le s\le \myfrac{1}{3}$ we have
\begin{eqn}\label{x2}
\sum_{i=1}^{n'}\,\Bigl|\,\log\,\tsi\,\Bigr|\,\le\,n'\cdot
\left(\,C\cdot\left|\int_{0}^{1}\log |x|\,dx\right|+C'\,\right)\,=
\,O(n)\,.
\end{eqn}

Let $d_i:=x_{i+1}-x_i$. Because of $\gm>0$, we have 
\begin{eqn}\label{63'}
\max_{i}\,d_i\,\le\,O(1/n)\,.
\end{eqn}
So, similarly as in \eqref{x2}
\begin{eqn}\label{x3}
\sum_{i=1}^{n'}\,\Bigl|\,d_i\,\log\,\tsi\,\Bigr|\,
\le O(1)\,.
\end{eqn}
Also, because of \eqref{63'} and $\int 1+\log x\ dx=x\cdot\log x+C$,
for $n$ large
\[
\left|\,\int_{x_i}^{x_{i+1}}\,\bigl(\,1+\log |s-s_i|\,\bigr)\,ds\,
\right|\,\le\,
-C\cdot d_i\cdot \left[\,\log d_i\,+\,C'\,\right]\,\le\,
\sqrt{d_i}\,\le\,O\left(\frac{1}{\sqrt{n}}\right)\,
\]
for some real number $C'$ independent on $i$ and $n$, and so
\[
\sum_{i=1}^{n'}\,\left|\,\int_{x_i}^{x_{i+1}}\,2\cdot
\log |s-s_i|\,ds\,\right|\,\le\,O(\sqrt{n})\,.
\]
This, \eqref{x3} and \eqref{x4} imply that for $s\in[0,1/3]$ the
(part of the) Jensen integral
\eqref{Ji} has a lower bound that behaves like $-O(\sqrt{n})$, so the
dominating part is for $s\in[1/3,1/2]$, and we are done for $V$.

For $\Dl$ the argument is similar. There $m(s)\equiv 2$, so
instead of \eqref{x5} we have for $s\in[x_i,x_{i+1}]$ the estimate
$|m(s)-h(s)|\ge |s-s_i|^2\cdot C\cdot \gm$, with the same further
reasoning.  \qed
}

\begin{rem}
Dan Silver pointed out that for the Alexander polynomial a proof
can be obtained (possibly more elegantly, but not immediately)
from the ideas in \cite{SW}.
\end{rem}

\subsubsection{Skein polynomial and generalized twisting}

By using the representation theory, one can extend the scope of
the results in \cite{SSW,CK} on bounded Mahler measure
to (parallel) multi-strand twisting.
We give a version for the skein polynomial, since its special case,
the Jones polynomial, was discussed in \cite{CK}. One can obtain
this special case by setting $\lm=q$ in the below theorem. (The
reversal of component orientation is not a serious problem for
$V$, unlike for $P$.) We write as before $P_L(v,z)$ for the skein
polynomial and use the skein rule $v^{-1}P_+-vP_-=zP_0$.

A \em{template} $T$ consists, as in \cite{SunThis}, of a number of
strands, that create crossings and \em{slots}, only that for us
strands are oriented and slots have the more general form
(with an arbitrary number of in- and outputs)
\[
\diag{7mm}{3}{2}{
  \picmultigraphics{5}{0.5 0}{\picvecline{0.5 0}{0.5 2}}
  \picfilledbox{1.5 1}{2.8 1}{}
}\,.
\]
A(n oriented) diagram $D$ is \em{associated} to $T$, if $D$ is
obtained from $T$ by inserting into each slot of $T$ a braid of a
certain number of full twists (on the proper number of strands).

The $1$-norm $||P||_1$ of a (Laurent) polynomial $P$ (even in
several variables) is understood to be the sum of absolute values
of all its (non-zero) coefficients (taken in all variables).
We write $M(P)$ for its Mahler measure. It is known that for
polynomials $P$ (of real coefficients), $M(P)\le ||P||_1$
(see \S 2 of \cite{SSW}).

\begin{theorem}
For each $n$ there is a polynomial $D_n(q)$ such that for
each template of $n_1,\dots,n_k$-strand parallel twist
slots and each diagram $D$ associated to $T$ we have
\[
\left|\left|\es\prod_{i=1}^k D_{n_i}(q)\es\cdot\es
P_D(\lm,\sqrt{q}-1/\sqrt{q})\es\right|\right|_1\le C_T\,,
\]
where $C_T$ is a constant that depends only on $T$.
Furthermore $D_n$ are made of a product of terms $1-q^{l}$
for some $1\le l\le n$. In particular, 
\[
\{\,M(\,P_D(\lm,\sqrt{q}-1/\sqrt{q})\,)\,:\,\mbox{\ $D$ is
associated to $T$\,}\,\}
\]
is bounded.
\end{theorem}

\proof Write $D(t_1,\dots,t_k)=D_{t_1,\dots,t_k}$ for the diagram
associated to $T$ by putting into the $i$-th slot $t_i$ full
twists (on $n_i$ strands).

Let first $k=1$. Then the complement of the slot is a room with
$n=n_1$ parallel in- and outputs. The skein module of such
a room is generated by (positive permutation) braids, and
thus it suffices to work with the $n$-strand braid
group. Again, by \cite[Proposition 6.2]{Jones}
\[
P_L(\lm,\sqrt{q}-1/\sqrt{q})=X_L(q,\lm^2/q)\,,
\]
where
% $X_L(q,\lm)$ is the invariant treated in \cite{Jones}.
% By the work of \cite{Jones} we know that
$X_L$ is a writhe-strand normalized (Definition 6.1 in \cite{Jones})
weighted trace sum (equation 5.5), and the full twist has the
effect of multiplying the traces by a power of $q$ (Lemma 9.3).

All terms in the demoninator of the (rational) generating series 
\begin{eqn}\label{GS}
\sum_{t=0}^{\infty}\,X(D_t)(q,\lm^2/q)\,x^t\,
\end{eqn}
are obtained from demoninators in the definition
of $X_L(q,\lm)$ and quantities that differ by the
number $t=t_1$ of twists. To identify these terms we 
will work up to units in $\bZ[q^{\pm 1/2},\lm^{\pm 1/2}]$
that do not depend on the number $t$ of twists but just 
on $T$. The quantities depending on $t$ are only
the power of $q$ in the traces and the term $\sqrt
{\lm}^e$ in Definition 6.1, where $e$ is the exponent sum.
Since the full twist has $e=n(n-1)$, we find the terms
$1-\lm^{n(n-1)/2} q^jx$ for some values of $j\ge 0$. The
exponent $n(n-1)/2$ becomes $n(n-1)$ and $j\ge -n(n-1)/2$
when replacing $\lm^2/q$ for $\lm$. It suffices now
to collect the demoninator terms occurring in the
definition of $X_L$ (and replace $\lm^2/q$ for $\lm$).
{}From Definition 6.1 we have $\lm^{(n-1)/2}$ (which is
a unit, fixed for given $T$, and remains so when
replacing $\lm^2/q$ for $\lm$) and some power of $1-q$.
{}From the text below Figure 5.6 and (5.5) we have
some power of $1-\lm q$ (that becomes $1-\lm^2$ under
substitution) and (in $Q(q)$) products of terms $1-q^l$
for some $1\le l\le n$ (with $l$ being a hook length of a
box in a Young diagram of $n$ boxes). Thus the demoninator
of \eqref{GS} is something we can take to be $D_n(q)$ times 
a product of terms of the form $1-\lm^2$ and $1-\lm^{n
(n-1)}q^{j}x$. 

Now if $k>1$, we first observe that $X(D_{t_1,\dots,t_k})(q,\lm)$
satisfy a linear recurrence (with coefficients in $\bZ[
q^{\pm 1/2},\lm^{\pm 1/2}]$) in $t_i$ for any fixed value
of $t_j, j\ne i$. Moreover that recurrence itself does not
depend on $t_j, j\ne i$ for fixed $i$, only its initial values
do depend. Then one inductively argues over $k$ that the
generating series 
\[
\sum_{t_1,\dots,t_k\ge 0} X(D_{t_1,\dots,t_k})(q,\lm^2/q)\es
x_1^{t_1} \dots x_k^{t_k}
\]
has a demoninator which is the product of the $D_{n_i}$ and terms
$1-\lm^2$ and $1-\lm^{n_i(n_i-1)}q^{j}x_i$ for $1\le i\le k$ and $j\ge 
-n_i(n_i-1)/2$. The cases when $t_i<0$ are analogous. From this
one easily concludes the claim up to the factors $1-\lm^2$.
Since for fixed $T$ the power $p_T$ of $1-\lm^2$ is fixed,
and the $\lm$-span of $P_D$ is bounded by the Morton-Williams-%
Franks inequality (see proposition 15.1 of \cite{Jones}),
one can get disposed of the $1-\lm^2$
factors by linear combinations of the coefficients of
$(1-\lm^2)^{p_T}P_D\cdot \prod D_{n_i}$ in the powers
of $\lm$. \qed

\subsection{Knots with unsharp MWF inequality}

This brief part was motivated by the paper of Kawamuro
\cite{Kawamuro}, which I came across during the continuous work
on my paper. Kawamuro was interested in finding infinitely many
knots for which the MWF inequality is strict, in particular such
satisfying Jones' conjecture in \S\ref{SJS}. The point to make
is that we have at least two ways allowing us to obtain her, and
further such knots much more easily. (In particular it is not
necessary to appeal to the heavy geometric machinery of
Neumann-Rudolph-Giroux.)

We already saw an infinite family of examples with unsharp MWF
(and satisfying Jones' conjecture) in the proof of theorem
\ref{thJC}. This family is, of course, little insightful,
since it is settled already by Birman-Menasco's 3-braid
work. The 4-braid examples they proposed, but could not
decide about, deserve more attention. As given in figure 4
of \cite{Kawamuro}, consider for $(x,y,z,w)\in\bZ^4$ the
4-braids $\bt_{x,y,z,w}=[212^213^x2^y-12^z3^w]$, and let
$K_{x,y,z,w}=\hat\bt_{x,y,z,w}$. Birman-Menasco observe that
$MWF(K_{x,y,z,w})\le 3$ (see lemma 2.9 in \cite{Kawamuro}),
but suspect that many of the $K_{x,y,z,w}$ have braid index 4.
One can exhibit an infinite family of $K_{x,y,z,w}$ of $b=4$,
and thus recover theorem 2.8 of \cite{Kawamuro}, like this.

\begin{prop}
There is a natural number $N$ and 6 tuples $(\ap_1,\dots,
\ap_4)\in \bZ^4$ of different mod-2-reductions in $(\bZ/2)
^4$, such that if $K_\gm:=K_{\gm_1,\dots,\gm_4}$ is a knot,
then $b(K_\gm)=4$ whenever $(\ap_i)\equiv (\gm_i)\bmod
2$ and $\gcd\Bigl(\ffrac{\gm_i-\ap_i}{2}\Bigr)\ge N$.
\end{prop}

\proof Since (complex) roots of unity are dense on the complex
unit circle, choose a root of unity $t$, for which the 3-braid
test in corollary \ref{pp1} applies for $9_{42}$ or $9_{49}$ (see
example \ref{4249}). Then, using the formula for $\Dl$ in lemma 3.1
of \cite{SSW}, we see that for each such $t$ there is an $n\in2\bZ
\sm\{0\}$, such that $\Dl_{\hat\bt}(t)$ is preserved when $\bt$
is replaced by $\sg_i^n\bt$.

Moreover, it is easy to check that for all 6 of the 16 possible
vectors of parities of $(x,y,z,w)$, for which $K_{x,y,z,w}$ is a
knot, there is a representation of (at least one of) $9_{42}$ or
$9_{49}$ of that parity combination. (The representations given
after definition 2.7 of \cite{Kawamuro} show 3 of the parity
types.) Now, take such a representation $\bt_{\ap_1,\dots,\ap_4}$
and vary the parameters $\ap_i$ by multiples of $n$. Then, since
we can choose (for proper $t$) $n$ to be any sufficiently
large even natural number, we obtain in the claim. \qed

Of course, this elementary argument could be further concretified
and strengthened. Presumably, one can do much better using
J.~Murakami's test theorem \reference{TMF} (since it is an equality).
Note that our proof does not confirm Jones' conjecture on any
$K_{x,y,z,w}$, so it may trigger the question: would it help
to find a counterexample? This seems very optimistic, though,
as we will soon see from an alternative approach to Kawamuro's
theorem, using the work in \S\ref{SJS}. First we can settle
(also with regard to Jones' conjecture) the concrete examples
$K_{-1,-2,m,2}$ for $m\ge 2$ even, obtained in her proof.

\begin{prop}\label{pZ}
We have $b(K_{-1,-2,m,2})=4$ for $m\ne -1,0$.
\end{prop}

\proof We give just a brief explanation. We try to calculate the
skein polynomial of the 2-cable $\bt^{[2]}_{-1,-2,m,2}$ of $\bt_
{-1,-2,m,2}$, obtained by replacing $\sg_i$ by $\sg_{2i}\sg_{2i+1}
\sg_{2i-1}\sg_{2i}$, and we may consider just the coefficient
of $v^{4[\bt_{-1,-2,m,2}]+7}$. Again $\pi_Y([2312])$ has at
most 7 different eigenvalues $\dl_i$, for all $Y\vdash 8$ taken
together. But now we have no full twist, and so do not need to
deal with the various $Y$ (and their weights $W_Y$ etc.) one by
one. We can sum $W_Y\cdot c(i,Y)$ over $Y\vdash 8$, and are left
with 7 different $q$-rational expressions $c_i$ to determine. 
Again this can be done from 7 explicit polynomials, which can be
obtained using Morton's program. (If we focus on one parity of $m$
only, we can again reduce the unknowns $c_i$ and test polynomials 
to 5, but need to calculate polynomials from slightly more
complicated braids.)

We calculated in fact 11 polynomials, for $|m|\le 5$, and
used those 9 for $|m|\le 4$ in the determination of $c_i$
to have extra safety. The resulting linear
equation system for $c_i$ (with matrix $\{\dl_i^m\}$) is now
not a serious problem to MATHEMATICA, which gives the solution
(for generic $q$) in just a few seconds. Clearing denominators,
and looking already at the extremal $q$-degrees, we see that for
$|m|\ge 6$ among the 7 terms $c_i\dl_i^m$ the one for $i=2$ has
the lowest minimal or highest maximal degree (as a unique term,
except for $m=6$, where we must use also that the leading
coefficient of $c_6\dl_6^m$ has the same sign). With an
explicit check for the other $m$, we conclude the claim. \qed

This method of proof has also the advantage of easily leading
to a qualitative improvement of Kawamuro's theorem, which is
closer to what one should expect.

\begin{prop}
The tuples $(x,y,z,w)$ for which $K_{x,y,z,w}$ has braid index 4 (and
satisfies Jones' conjecture) are generic in $\bZ^4$, in the sense that
\begin{eqn}\label{bla}
\lim_{n\to\infty}\es
  \frac{
   \big|\,
    \{\,(x,y,z,w)\in \bZ^4\,:\,b(K_{x,y,z,w})=4\,\}\,\cap\, {[-n,n]^4} 
   \,\big|
  \rule[-0.6em]{\z@}{1.3em}}
  {\big|\,[-n,n]^4\,\big|\rule[0.1em]{\z@}{1.0em}}
\es=\,1\,.
\end{eqn}
\end{prop}

\proof Extending the argument for proposition \ref{pZ}, we see that 
\[
\tl P(x,y,z,w):=
  [P(\hat\bt^{[2]}_{x,y,z,w})]_{v^{4[\bt_{x,y,z,w}]+7}}(
  \sqrt{q}-1/\sqrt{q})\cdot \sqrt{q}^{\,4[\bt_{x,y,z,w}]+7}
  \,=\,\sum_{i_x,i_y,i_z,i_w=1}^7\,c_{i_x,i_y,i_z,i_w}\dl_{i_x}^x
  \dl_{i_y}^y\dl_{i_z}^z\dl_{i_w}^w\,,
\]
where $c_{i_x,i_y,i_z,i_w}$ are again some rational expressions
in $q$. They are obviously not all zero, and the polynomials $\tl
P(x,y,z,w)$ have the following property: if one fixes three of the
parameters (say $x,y,z$), and $\tl P(x,y,z,w)=0$ for 7 different
values of the remaining parameter (here $w$), then it vanishes
for all values of that parameter (at the given fixed 3 others).
It is then not hard to deduce \eqref{bla} from this (see e.g.
the proof of lemma 11.1 of \cite{gen2}). \qed

\noindent{\bf Acknowledgement.} I would like to thank to M. Hirasawa,
M. Ishiwata and K.~Murasugi for their contribution, and S.~Kamada,
K.~Kawamuro, T.~Nakamura and D.~Silver for some helpful remarks.
The calculations were performed largely using the programs of
\cite{MorSho2} and \cite{KnotScape}, and MATHEMATI\-CA\TM{}
\cite{Wolfram}. I also wish to thank to my JSPS host Prof. T.~Kohno
at University of Tokyo for his support.

% \vspace{1cm}\mbox{}
% \newpage

\begin{appendix}

\renewcommand{\thesection}{\Alph{section}}
{\def\@seccntformat#1{Appendix \csname the#1\endcsname . \quad}

\section{Postliminaries}
}

\subsection{Fibered Dean knots (Hirasawa-Murasugi)\label{HM}}

Here we present some material due to Hirasawa and Murasugi,
who studied fibering of generalized Dean knots. An overview is given
in \cite{Hirasawa}. As this is work in progress, a longer exposition
may appear subsequently elsewhere.

\begin{defi}
The Dean knot $K(p,q|r,rs)$ is given by the closed $p$-braid
\[
(\sg_{p-1}\sg_{p-2}\cdots \sg_1)^{q}(\sg_1\sg_2\cdots \sg_{r-1})^{rs}\,,
\]
with $p>r>1$ and $q,s$ non-zero integers such that $(q,p)=1$.
\end{defi}

Hirasawa and Murasugi proposed a conjecture on these knots
and obtained so far the following partial results.
 
\begin{conjecture}(Hirasawa-Murasugi)
A Dean's knot $K(p,q|r,rs)$ is a fibred knot if and
only if its Alexander polynomial is monic, that is,
$\Mc\Dl=\pm 1$.
\end{conjecture}

\begin{figure}[th]
\[
\begin{array}{ccc}
\epsfs{6cm}{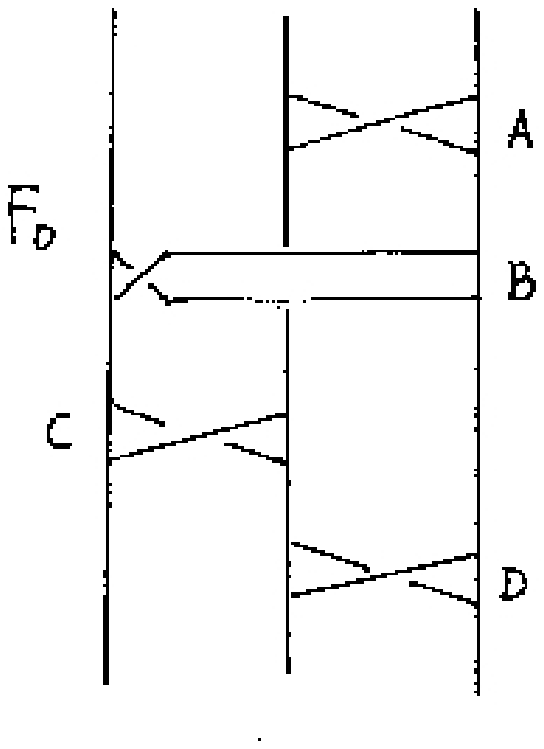} & 
\epsfs{5cm}{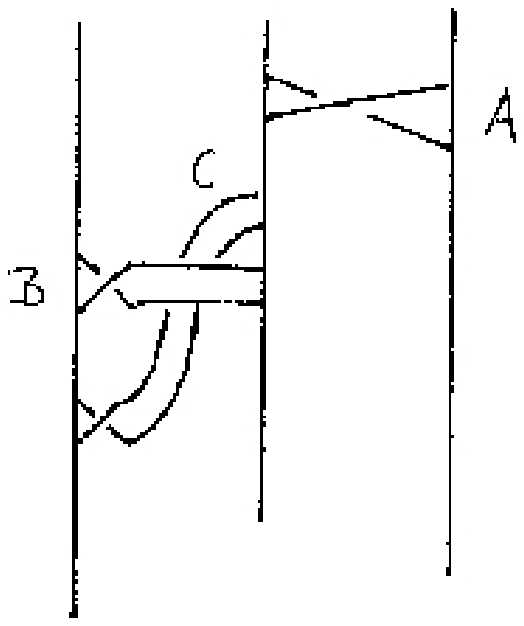} & 
\epsfs{5cm}{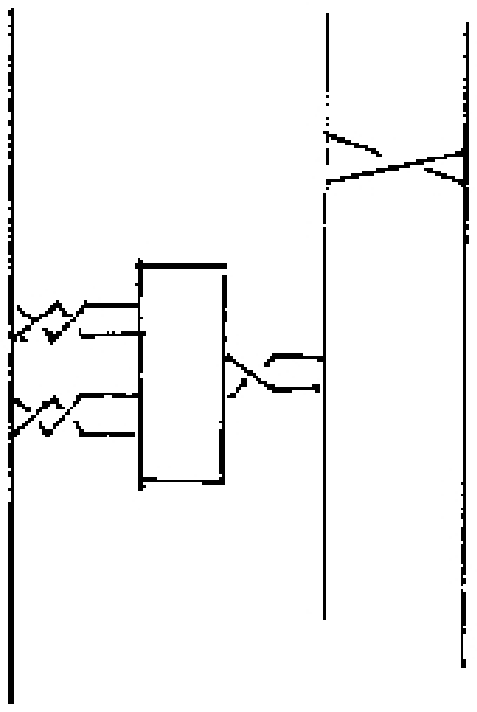}\\
 (a) & (b) & (c) \\
\end{array}
\]
\caption{\label{fig2}}
\end{figure}

\begin{prop}(Hirasawa-Murasugi)
This conjecture has been proven for the following cases.
\def\theenumi{\alph{enumi}}
\def\labelenumi{(\theenumi)}
\begin{enumerate}
\item     $q = kp + 1$, and $r$ and $s$ are arbitrary,
\item     $q =kp - 1$, and  $r$ and $s$ are arbitrary,
\item     $r =  p - 1$, and  $q$ and $s$ are arbitrary.
\end{enumerate}
\end{prop}

The last case implies in particular that the conjecture is
true for $p=3$. Below follows a part of the argument in this
case that settles lemma \reference{lki}. (Other parts of their
proof are very similar to some of our previous arguments.
It seems, for example, that Hirasawa and Murasugi were to
some extent aware of Theorem \reference{Thq}.)

\proof [of lemma \reference{lki}] We consider the braids $[(123)^k-2]$,
with $k>0$ fixed. By isotopies and Hopf plumbings, we modify our
surface 
\begin{eqnarray*}
% & (123)^k-2\to 1223(123)^{k-1}-2\to -21223(123)^{k-1}= & \\
% & 12 -123(123)^{k-1}\to 23(123)^{k-1}12 -1= 2(312)^k-1\to
% (2312)^k-1\,.
&\rx{-1mm} [(123)^k -2] \to [(1223)^k - 2]=[12(2312)^{k-1}23 - 2] \to
[(2312)^{k-1}23 - 212] = [(2312)^{k-1}2312 -1]=[(2312)^k - 1].
\end{eqnarray*}
Let $F$ be the surface of the band representation $\bt=[(2312)^k-1]$.
We will show that $F$ is a fiber surface.

Consider the subsurface $F_0$ of $F$ that spans $[2 3 1 2]$
in the natural manner.
We deform the $k$ copies of $F_0$ to the (isotopic sub)surfaces
$F'_0$ by a series of diagrams, see figure \reference{fig2}.
(Here strands are numbered from right to left and words composed
downward.)

In figure \reference{fig2}, for the move (a) $\to$ (b), we slide
$B$ and $C$, respectively, along $D$ and $A$, then delete $D$
(by deplumbing a Hopf band), and then slide $C$ back along
$A$. For the move (b) $\to$ (c), we slide $C$ along $B$,
then slide $B$ along $C$. The $k$ bands $A$ can be subsequently
removed by Murasugi desumming a $(2,k)$-torus link fiber surface.
Thus the surface $F$, spanned by $\hat\bt$, turns after
(de)summing Hopf bands into a surface $F'$ consisting of $k$
copies of $F'_0$ and one negative band $N$. See figure \ref{fig3} (a).

{
% \tm
\def\@captype{figure}
\long\def\@makecaption#1#2{%
   % \tm
   \vskip 10pt
   {\let\label\@gobble
   \let\ignorespaces\@empty
   \xdef\@tempt{#2}%
   %\typeout{`#2'}%
   }%
   \ea\@ifempty\ea{\@tempt}{%
   \setbox\@tempboxa\hbox{%
      \fignr#1#2}%
      }{%
   \setbox\@tempboxa\hbox{%
      {\fignr#1:}\capt\ #2}%
      }%
   \ifdim \wd\@tempboxa >\captionwidth {%
      \rightskip=\@captionmargin\leftskip=\@captionmargin
      \unhbox\@tempboxa\par}%
   \else
      \hbox to\captionwidth{\hfil\box\@tempboxa\hfil}%
   \fi}%
% %
% \def\fignr{\small\sffamily\bfseries}%
% \def\capt{\small\sffamily}%
% 
% 
% \newdimen\@captionmargin\@captionmargin2cm\relax
% \newdimen\captionwidth\captionwidth\hsize\relax

\
\captionwidth0.3\textwidth\relax
\[
\begin{array}{cc@{\qquad}c}
\epsfs{5cm}{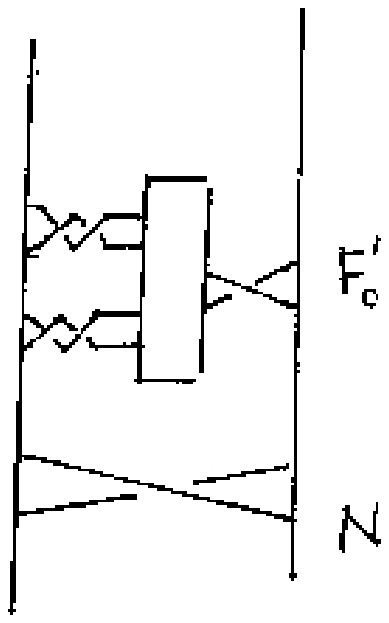} & \epsfs{5cm}{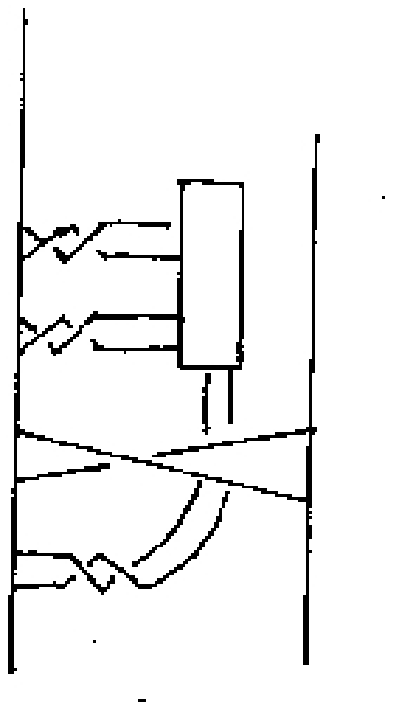} & \epsfs{5cm}{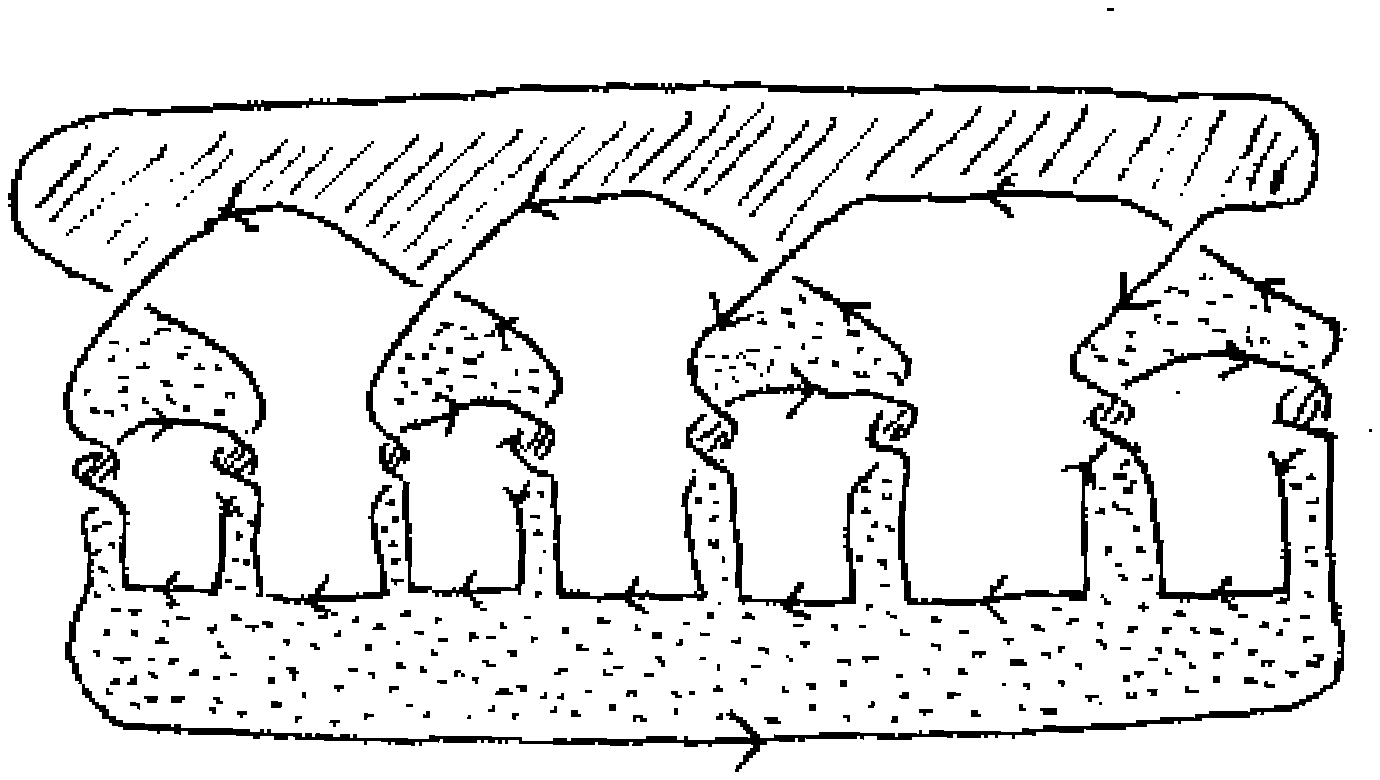} \\
(a) & (b) \\[-2mm]
\multicolumn{2}{c}{ \vbox{\caption{\label{fig3}}} } &
\raisebox{1em}{\vbox{\caption{\label{fig4}}}}
\end{array}
\]

}

In figure \reference{fig3}, we perform the move (a) $\to$ (b), sliding
$B$ along $N$. The last surface $F'$, in figure \reference{fig3} (b),
is Murasugi sum of a fibre surface spanning the $(2,2,-2)$-pretzel
link and $\tl F$, where $\tl F$ consists 
of $k - 1$ copies of $F'_0$ and the band $N$. By induction on
$k$, we see that $F'$ is a fibre surface, and hence $\hat\bt$
is a fibred link. \qed
 
\subsection{$A$-decomposition (joint with Hirasawa-Ishiwata)\label{HI}}

For the proof of theorem \reference{thuq} we introduce the
$A$-decomposition due to T.~Kobayashi \cite{Kobayashi}.

A \em{sutured manifold} in the sense of Gabai \cite{Gabai3} can be
understood as a pair $(L,H)$ consisting of a closed 3-dimensional
submanifold $H$ of $\bR^3$ with boundary $S=\prt H$ a connected
surface, and a set of oriented loops $L\subset H$. We require that
one can orient the connected components of $S\sm L$ so that the
induced orientations on $L$ coincide from both sides of $L$ (in
particular a connected component of $S\sm L$ never bounds to itself
along a loop of $L$), and are given by the orientation of $L$.

Let $F$ be a connected Seifert surface of a(n oriented) link $L=\prt
F$. We embed $F$ as $F\times\{0.5\}$ into the \em{bicolar} $H=
F\times I$ (with $I=[0,1]$). Then $(L,H)$ becomes a sutured
manifold. We call it \em{canonical sutured manifold} $C(F)$ of $F$.

We describe some basic operations on sutured manifolds $(L,H)$.

A \em{decomposition disk} $D$ is a disk with $P=\prt D\subset \prt
H$, properly embedded in the complement of $H$ (i.e. $D\cap H=P$).
We require that $D$ is not parallel to $S=\prt H$, and
satisfies $P\cap F\ne\vn$. We assume also that the intersection of
$P$ and $F$ is transversal, so that it is a collection of points.

Since $L=\prt F$ is separating on $S$, the intersection $D\cap L= P
\cap F$ is an even number of points, and the orientation of $L$ at
the intersection points is alternating (with respect to the
orientation of the loop $P$). See figure \reference{fig5} (a).

\begin{figure}[th]
\[
\begin{array}{c@{\kern1.4cm}c}
\diag{1cm}{3}{3}{
  \piccircle{1.5 d}{0.6}{$D$}
  \pictranslate{1.5 1.5}{
    { \piclinedash{0.05 0.05}{0.02}
      \picmultigraphics[rt]{4}{90}{
        \picline{0.0 0.6}{-0.3 1.4}
      }
    }
    \picputtext{0.4 0.9}{$L$}
    \picmultigraphics[rt]{2}{180}{
      \picvecline{0.3 1.4}{0.0 0.6 0.3 1.4 0.5 conv}
      \picline{0.0 0.6}{0.0 0.6 0.3 1.4 0.5 conv}
      \picrotate{90}{
        \picvecline{0.0 0.6}{0.0 0.6 0.3 1.4 0.5 conv}
        \picline{0.3 1.4}{0.0 0.6 0.3 1.4 0.5 conv}
      }
    }
  }
}
&
\diag{1cm}{2.5}{2.0}{
  \picline{1 0}{1 2}
  \picline{0 0.5}{1 0.5}
  \picvecline{2 0.5}{1 0.5}
  \picvecline{0 1.5}{2 1.5}
  \picputtext{1 2.3}{$D$}
  \picputtext{2.3 0.5}{$L$}
  \picputtext{2.3 1.5}{$L$}
  \picputtext{0.8 1.0}{$a$}
}\quad\lra\quad
\diag{1cm}{2}{2.0}{
  \picline{1 0}{1 2}
  \picputtext{1 2.3}{$D$}
  \pictranslate{1 1}{
    \picmultigraphics[rt]{2}{180}{
      \picline{1 -0.5}{0.2 -0.5}
      \piclineto{0.2 0.5}
      \picveclineto{1 0.5}
    }
  }
}   
\\
\ry{1.7em}(a) & \raisebox{0.8em}{$(b)$} \\
\end{array}
\]
\caption{\label{fig5}}
\end{figure}

Then $L\cap P$ separates $P$ into a collection of intervals or
\em{arcs}. Let $a$ be such an arc. An \em{$A$-operation} on $D$ along
$a$ is a transformation of $(L,H)$ into a sutured manifold $(L',H)$,
where $L'$ is obtained by splicing $L$ along $a$. See figure
\ref{fig5} (b).

A \em{product decomposition} along $D$ is a similar operation, due
to Gabai \cite{Gabai3}, and can be described as an $A$-operation
if $|L\cap D|=2$ (in which case which of the two arcs is chosen
is irrelevant), followed by a subsequent gluing of a $D^2\times I$ into
$H$ along a neighborhood $N(P)\simeq S^1\times I$ of $P$ on $S$.

\begin{defi}
We define a sutured manifold $(L,H)$ to be \em{$A$-decomposable} as
follows:
\def\theenumi{\arabic{enumi}}
\def\labelenumi{\theenumi)}
\begin{enumerate}
\item Assume $H$ is a standardly embedded handlebody (i.e. so that
$\ol{S^3\sm H}$ is also one). If $L$ is a collection of trivial loops
on $\prt H$, and all loops bound \em{disjoint} disks in $\prt H$,
then $(L,H)$ is $A$-decomposable.
\item If $(L',H')$ is obtained from $(L,H)$ by a product decomposition
(along some decomposition disk $D$),
and $(L',H')$ is $A$-decomposable, then so is $(L,H)$.
\item Let $D$ be a decomposition disk of $H$ with $|L\cap D|=2n$ and
choose among the $2n$ arcs on $P=\prt D$ a collection of $n$ cyclically
consecutive arcs $a_1,\dots,a_n$. (Consecutive is to mean that, taken
with their boundary in $L\cap P$, their union is a single interval
in $P$, and not several such intervals.) Let $(L_i,H)$ be obtained
from $(L,H)$ by $A$-decomposition on $D$ along $a_i$ for $i=1,\dots,n$.
Then if all $(L_i,H)$ are $A$-decomposable, so is $(L,H)$.
\end{enumerate}
\end{defi}

\begin{theorem}(T.~Kobayashi \cite{Kobayashi}, O.~Kakimizu
\cite{Kakimizu}) 
\def\theenumi{\arabic{enumi}}
\def\labelenumi{\theenumi)}
\begin{enumerate}
\item A fiber surface is a unique incompressible surface.
\item The property a surface to be a unique incompressible surface 
  is invariant under Hopf (de)plumbing.
\item If $C(F)$ is $A$-decomposable, then $F$ is a
  unique incompressible surface for $L=\prt F$.
\end{enumerate}
\end{theorem}

\begin{figure}%[th]
\[
\epsfsv{4.5cm}{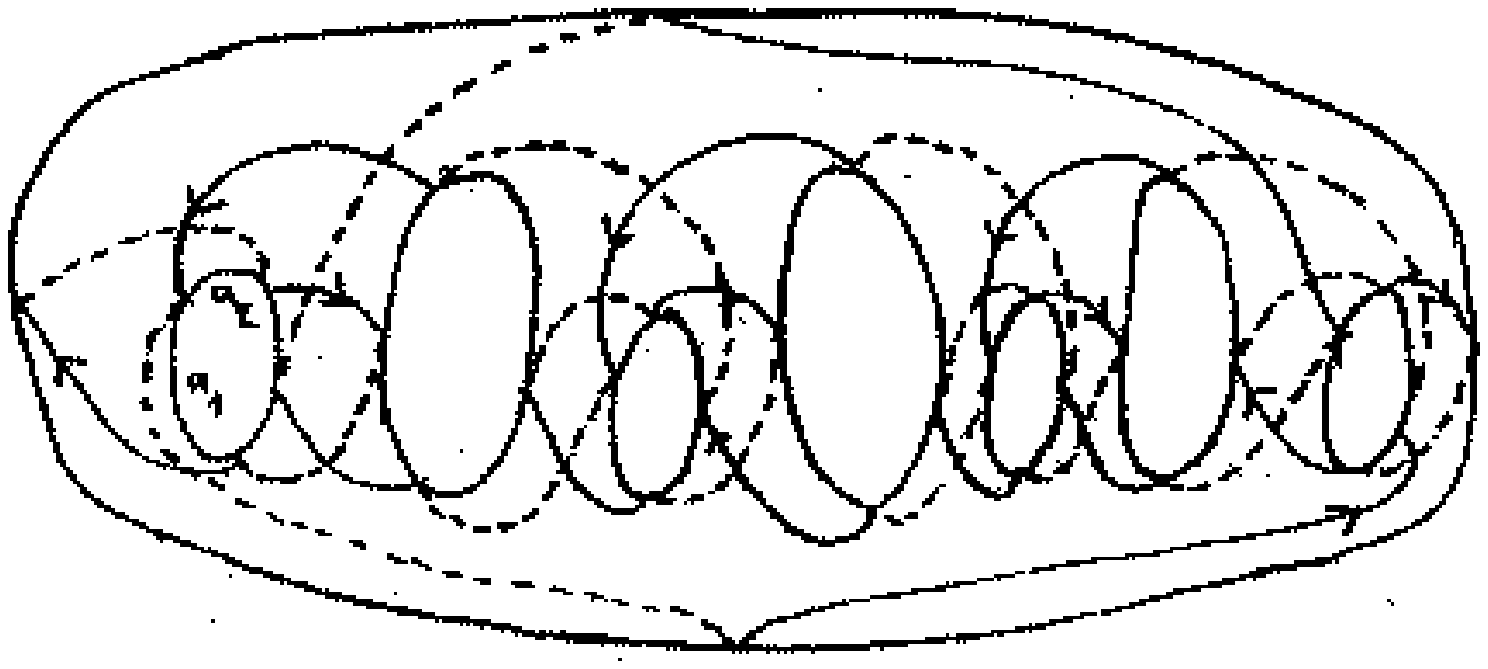} \stackrel{a_1}{\lra} \epsfsv{4.7cm}{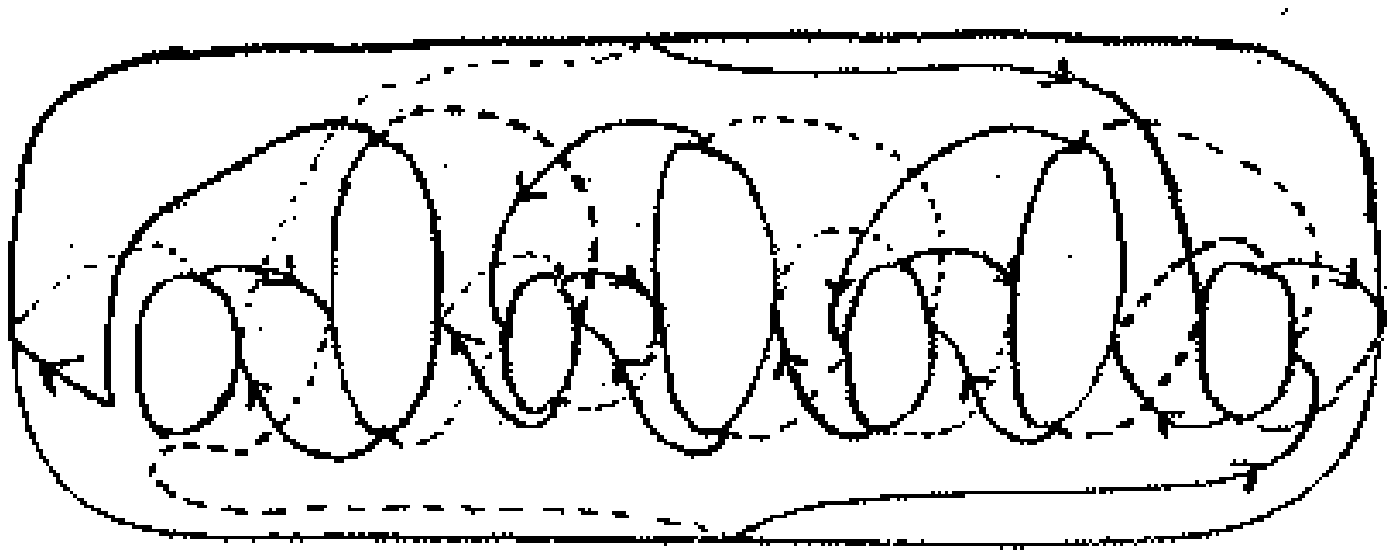}
\]
\[
\epsfsv{4.5cm}{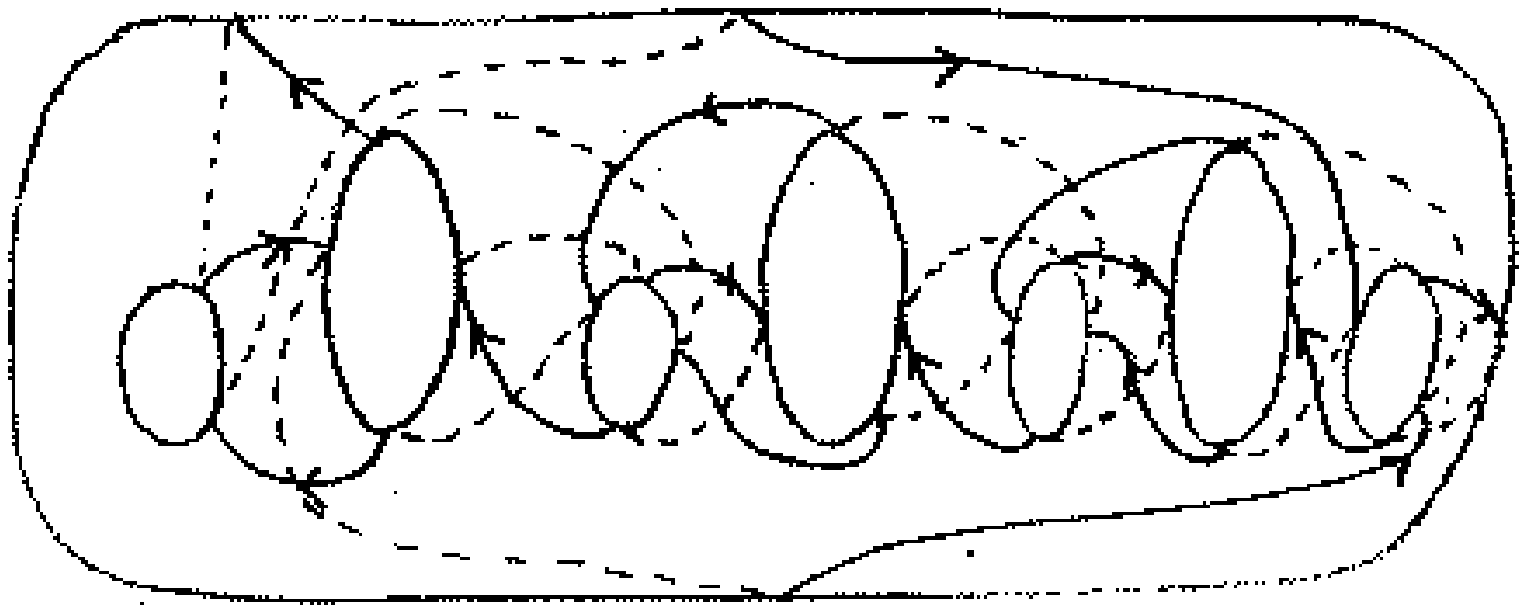} \qquad \epsfsv{4.5cm}{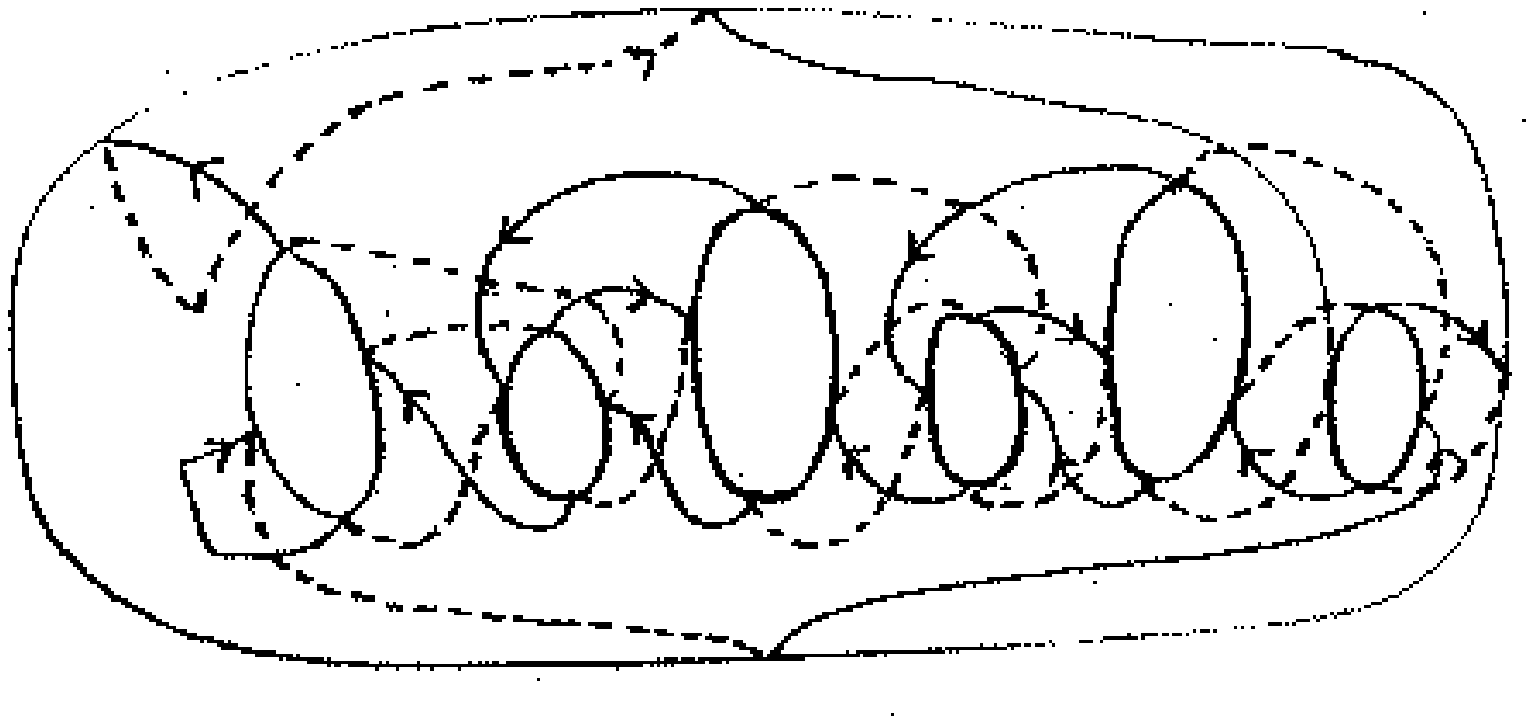} 
\]
\[
\epsfsv{4.5cm}{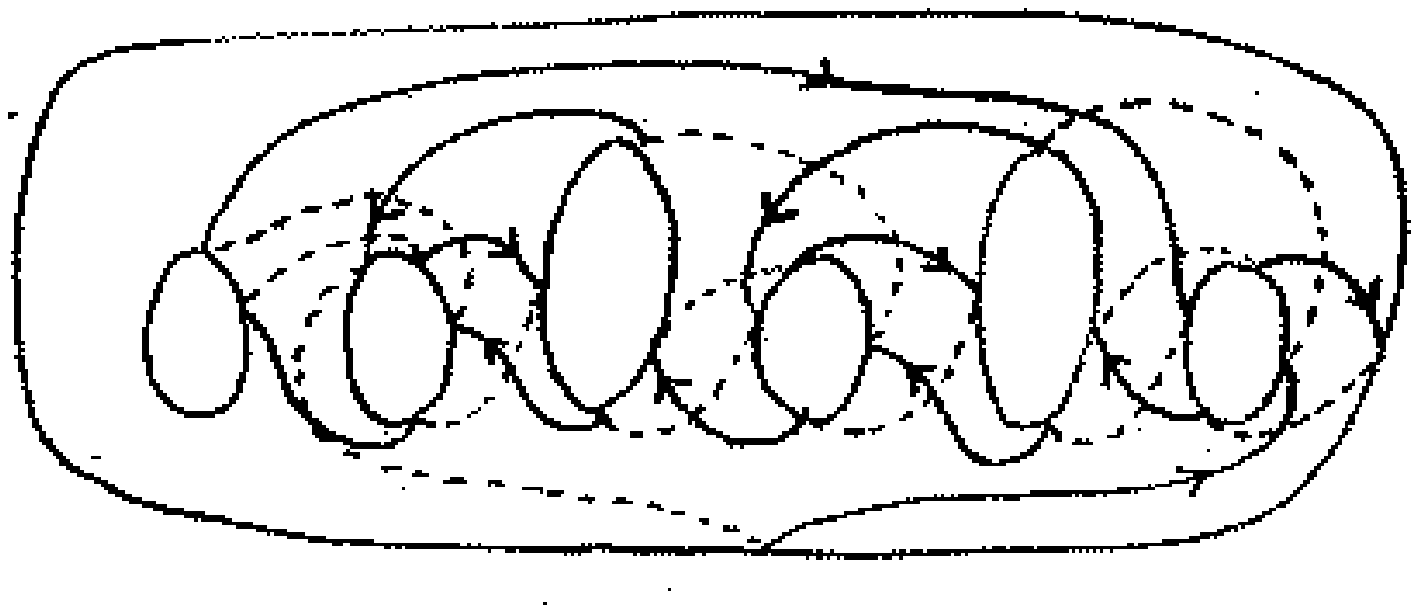} \qquad \epsfsv{4.5cm}{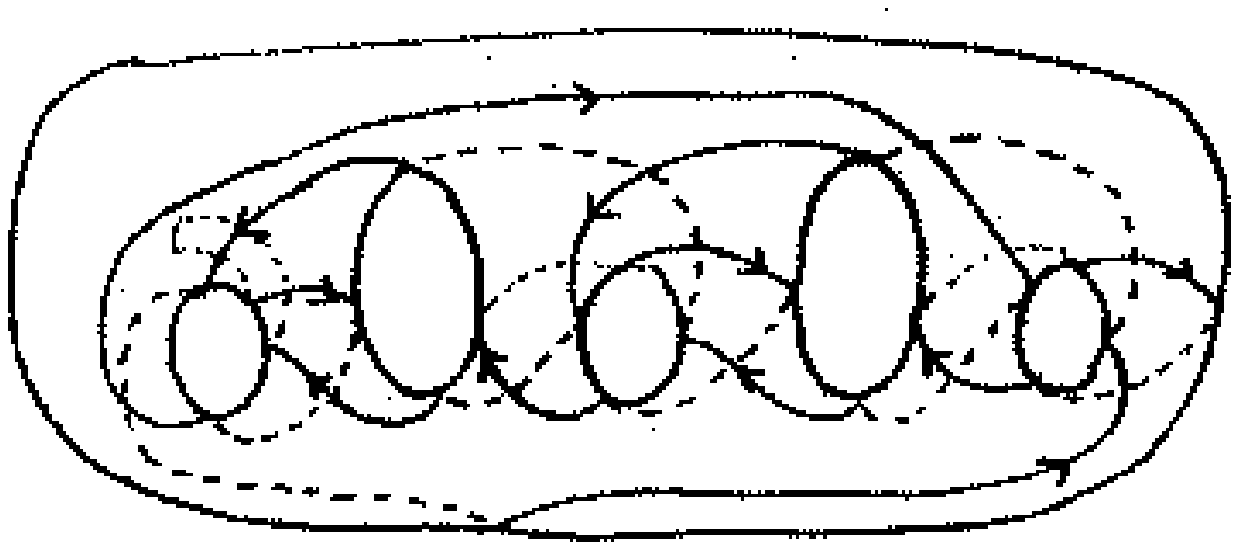}
\]
\[
\epsfsv{4.2cm}{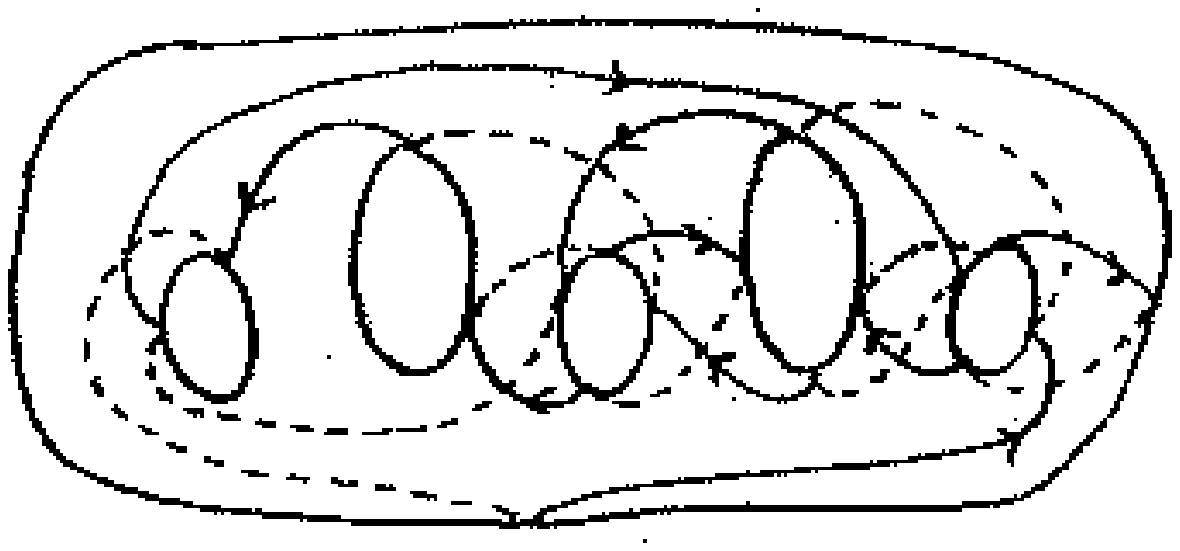} \qquad \epsfsv{4.2cm}{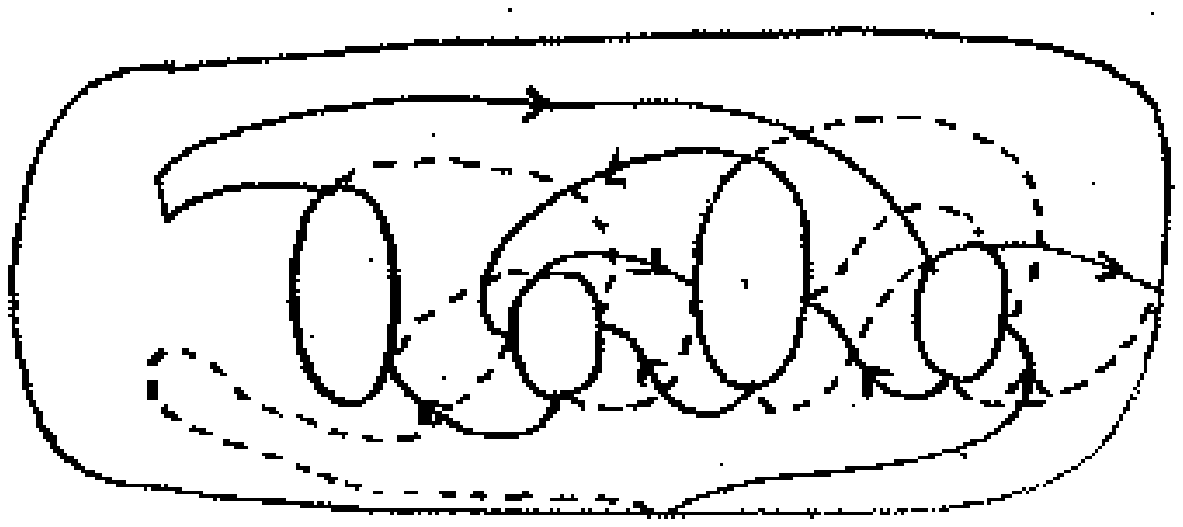}
\]
\[
\epsfsv{3.5cm}{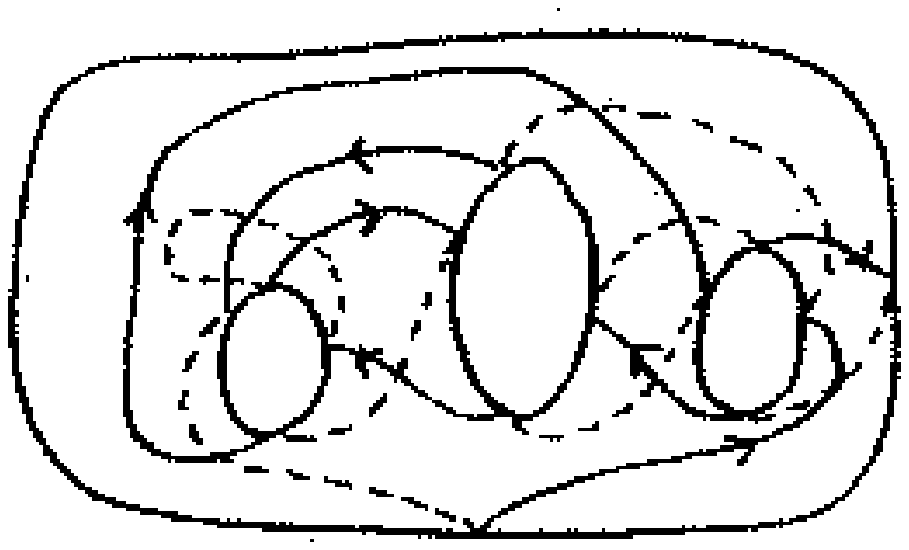} \qquad \epsfsv{3.7cm}{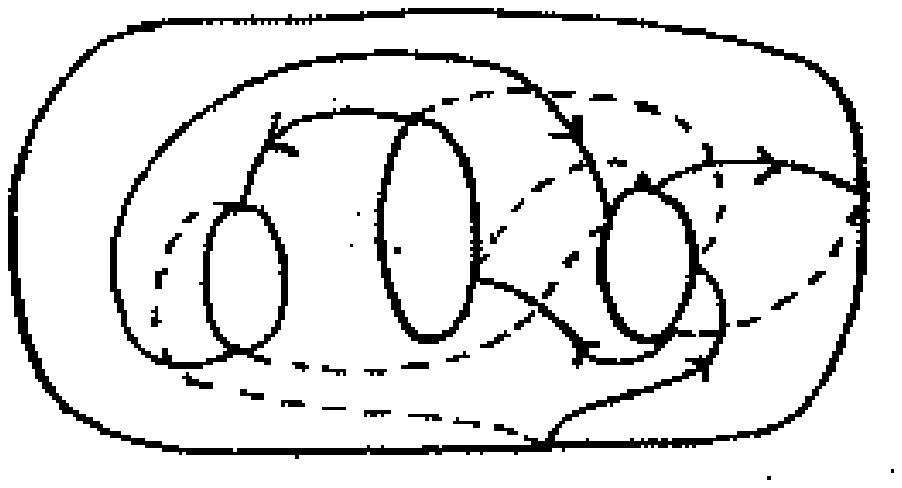}
\]
\[
\epsfsv{3.5cm}{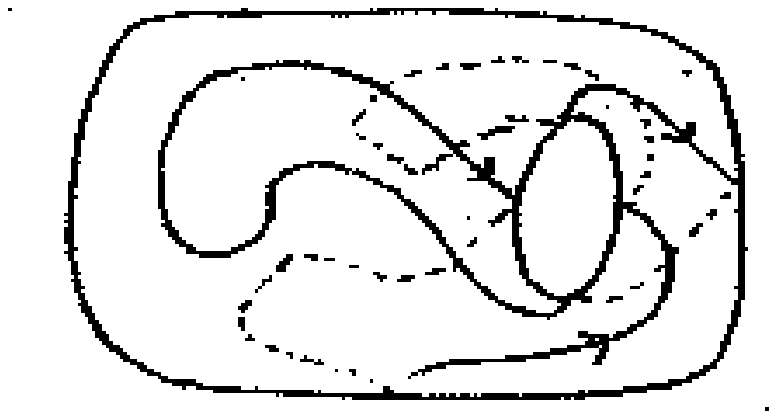} \qquad \epsfsv{3.2cm}{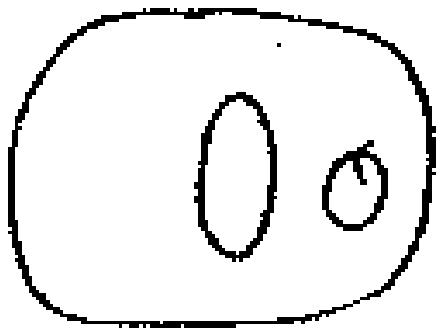}
\]
\caption{\label{fig6} This sequence of diagrams describes the
$A$-decomposition along arc $a_1$ of the canonical sutured
manifold corresponding to the surface in figure \reference{fig4}.
The $A$-operation is the change between the first and the second
diagram. The following diagrams (left to right in each row)
display isotopies and product decompositions, leading finally to
a trivial curve on a torus. The $A$-operation along $a_2$
must be performed too, but the decomposition is similar.}
\end{figure}

Now we complete the proof of theorem \reference{thuq}.

\proof[of theorem \reference{thuq}] 
Let $L$ be a 3-braid link. If $L$ is split, the splitting sphere
also splits any incompressible surface for $L$. Since
$L$ is a split union of a 2-braid link and an unknot
(incl.~2 and 3-component unlinks), the claim is easy.
Excluding split links, corollary \reference{cbd} shows we
need to consider only connected Seifert surfaces. It also
suffices to deal with the non-fibered links only. These are
equivalent under Hopf (de)plumbing to some of $(123)^k$ ($k\ne
0$). For such links, we can modify the transformation of
Hirasawa and Murasugi from \S\reference{HM}. Then we can turn
by Hopf (de)plumbing the surfaces into those like in figure
\reference{fig3} (a), consisting of $k$ copies of $F_0'$, but
now \em{without} the lower band $N$. (Figure \reference{fig4}
shows the case $k=4$.) Then we use $A$-decomposition, as
shown in figure \reference{fig6}. One applies $A$-operations
along the arcs $a_{1,2}$. We show only the result for $a_1$,
the case of $a_2$ and other $k$ is analogous. \qed

\end{appendix}

% \vspace{1cm}\mbox{}

% \tm


{\small
\begin{thebibliography}{BLM}
\bibitem[Al]{Alexander} J. W.~Alexander, \em{A lemma on systems of
        knotted curves}, Proc. Nat. Acad. Sci. U.S.A.
        {\bf 9} (1923), 93--95.
\bibitem[Al2]{Alexander2} \bysame, \em{Topological invariants
        of knots and links}, Trans.\ Amer.\ Math.\ Soc. {\bf 30} (1928),
        275--306.
\bibitem[Ar]{Artin} E.~Artin, \em{Theorie der Z\"opfe}, 
        Abh. Math. Sem. Hamburgischen Univ. {\bf 4} (1926), 47--72.
\bibitem[Ar2]{Artin2} \bysame, \em{Theory of braids}, 
        Ann. of Math. {\bf 48(2)} (1947), 101--126.
\bibitem[Be]{Bennequin} D.~Bennequin, \em{Entrelacements et \'equations
	de Pfaff}, Soc.\ Math.\ de France, Ast\'erisque {\bf 107-108}
	(1983), 87--161.
\bibitem[Bi]{Bigelow} S.~Bigelow, \em{Representations of braid
	groups}, Proceedings of the International Congress of
	Mathematicians, Vol. {\bf II} (Beijing, 2002),  37--45.
\bibitem[B]{Birman} J.\ S.\ Birman, \em{On the Jones polynomial of
	closed $3$-braids}, Invent. Math. {\bf 81(2)} (1985), 287--294.
\bibitem[BKL]{BKL} \bysame, K.\ Ko and S.\ J.\ Lee, \em{A new
	approach to the word and conjugacy problems in the braid
	groups}, Adv. Math. {\bf 139(2)} (1998), 322--353.
\bibitem[BM]{BirMen} \bysame\ and W.~W.~Menasco, \em{Studying knots
	via braids III: Classifying knots which are closed 3 braids},
	Pacific J.~Math. {\bf 161}(1993), 25--113.
\bibitem[BM2]{BirMen2} \bysame\ and \bysame\,, \em{Studying links
	via closed braids II: On a theorem of Bennequin},
	Topology Appl. {\bf 40(1)} (1991), 71--82.
\bibitem[BM3]{BirMen3} \bysame\, and \bysame\,, \em{Studying knots
	via braids VI: A non-finiteness theorem},
	Pacific J.~Math. {\bf 156} (1992), 265--285.
\bibitem[BW]{BW} \bysame\ and H.\ Wenzl, \em{Braids, link polynomials
	and a new algebra}, Trans. Amer. Math. Soc. {\bf 313(1)}
	(1989), 249--273.
\bibitem[Bl]{Bleiler} S.~A.~Bleiler, \em{Realizing concordant
	polynomials with prime knots}, Pacific J. Math. {\bf 100(2)}
	(1982), 249--257.
\bibitem[BLM]{BLM} R.~D.~Brandt, W. B. R. Lickorish and K. Millett,
	\em{A polynomial invariant for unoriented knots and links},
	Inv. Math. {\bf 74} (1986), 563--573.
\bibitem[CK]{CK} A.~Champanerkar and I.~Kofman, \em{On the Mahler
	measure of Jones polynomials under twisting},
	Algebr. Geom. Topol. {\bf 5} (2005), 1--22.
\bibitem[Cr]{Cromwell} P.~R.~Cromwell, {\em Homogeneous links},
	J. London Math. Soc. (series 2) {\bf 39} (1989), 535--552.
\bibitem[DL]{DasLin} O.~Dasbach and X.-S.~Lin, {\em On the Head and
	the Tail of the Colored Jones Polynomial}, Compositio Math.
	{\bf 142(5)} (2006), 1332--1342.
\bibitem[EKT]{EKT} S.\ Eliahou, L.\ H. Kauffman and M.\ Thistlethwaite,
	{\em Infinite families of links with trivial Jones polynomial},
	Topology {\bf 42(1)} (2003), 155--169.
\bibitem[Fi]{Fied} T.~Fiedler, {\em A small state sum for knots},
	Topology {\bf 32~(2)} (1993), 281--294.
\bibitem[Fi2]{Fied2} \bysame, {\em Gauss sum invariants for knots
	and links}, Kluwer Academic Publishers, Mathematics and Its
	Applications Vol {\bf 532} (2001).
\bibitem[FK]{FK} \bysame\, and V.\ Kurlin, \em{A one-parameter
	approach to knot theory}, preprint \web|math.GT/0606381|.
\bibitem[FW]{WilFr} J.\ Franks and R.\ F.\ Williams, \em{Braids and the
	Jones-Conway polynomial}, Trans.\ Amer.\ Math.\ Soc. {\bf 303}
	(1987), 97--108.
\bibitem[F\&]{HOMFLY} P. Freyd, J. Hoste, W. B. R. Lickorish,
	K. Millett, A. Ocneanu and D. Yetter, {\it A new polynomial
	invariant of knots and links}, Bull. Amer. Math. Soc.
	{\bf 12} (1985), 239--246.
\bibitem[Fu]{Fujii} H.~Fujii, \em{Geometric indices and the Alexander
	polynomial of a knot},
	Proc. Amer. Math. Soc. {\bf 124(9)} (1996), 2923--2933.
\bibitem[Ga]{Gabai} D.~Gabai, \em{The Murasugi sum is a natural
	geometric operation}, Low-dimensional topology (San Francisco,
	Calif., 1981), Contemp. Math. {\bf 20}, 131--143, Amer. Math.
	Soc., Providence, RI, 1983.
\bibitem[Ga2]{Gabai2} \bysame, \em{The Murasugi sum is a natural
	geometric operation II}, Combinatorial methods in
	topology and algebraic geometry (Rochester, N.Y., 1982),
	93--100, Contemp. Math. {\bf 44}, Amer. Math. Soc.,
	Providence, RI, 1985.
\bibitem[Ga3]{Gabai3} \bysame, \em{Detecting fibred links in $S^3$},
 	Comment. Math. Helv.  {\bf 61(4)} (1986), 519--555.
\bibitem[Gr]{Garside} F.~Garside, \em{The braid group and other groups},
	Q.~J.~Math.~Oxford {\bf 20} (1969), 235--264.
\bibitem[HM]{Hirasawa} M.~Hirasawa and K.~Murasugi, \em{Double-torus
	fibered knots and pre-fiber surfaces}, Musubime to Teijigen
	Topology (Dec. 1999), 43--49.
\bibitem[HT]{KnotScape} J.~Hoste and M.~Thistlethwaite, {\em
	KnotScape}, a knot polynomial calculation and table access
	program, available at \webb:|http://www.math.utk.edu/\~morwen|.
\bibitem[HTW]{HTW} \bysame, \bysame\, and J.~Weeks, \em{The first
	1,701,936 knots}, Math. Intell. {\bf 20 (4)} (1998), 33--48.
\bibitem[I]{Ishikawa} M.~Ishikawa, {\em On the Thurston-Bennequin
	invariant of graph divide links}, Math. Proc. Cambridge
	Philos. Soc. {\bf 139(3)} (2005), 487--495.
\bibitem[J]{Jones} V.~F.~R.~Jones, {\em Hecke algebra representations 
	of braid groups and link polynomials}, Ann. of Math.
	{\bf 126} (1987), 335--388.
\bibitem[J2]{Jones2} \bysame, {\em A polynomial
	invariant of knots and links via von Neumann algebras},
	Bull. Amer. Math. Soc. {\bf 12} (1985), 103--111.
\bibitem[Kk]{Kakimizu} O.~Kakimizu, \em{Classification of the
	incompressible spanning surfaces for prime knots of $\leq 10$
	crossings}, Hiroshima Math. J. {\bf 35} (2005), 47--92.
\bibitem[K]{Kanenobu} T.~Kanenobu, \em{Relations between the
	Jones and Q polynomials of 2-bridge and 3-braid links},
	Math. Ann. {\bf 285} (1989), 115--124.
\bibitem[K2]{Kanenobu2} \bysame, \em{Examples on polynomial invariants
	of knots and links II}, Osaka J. Math. {\bf 26(3)} (1989),
	465--482.
\bibitem[K3]{Kanenobu3} \bysame, \em{Examples on polynomial invariants
	of knots and links}, Math. Ann. {\bf 275} (1986), 555--572.
\bibitem[K4]{Kanenobu4} \bysame, \em{An evaluation of the first
	derivative of the $Q$ polynomial of a link},
	Kobe J. Math. {\bf 5(2)} (1988), 179--184.
\bibitem[Ka]{Kauffman} L.\ H.\ Kauffman, {\em An invariant of regular
	isotopy}, Trans. Amer. Math. Soc. {\bf 318} (1990), 417--471.
\bibitem[Ka2]{Kauffman2} \bysame, {\em State models and
	the Jones polynomial}, Topology {\bf 26} (1987), 395--407.
\bibitem[Kw]{Kawamuro} K.~Kawamuro, \em{The algebraic crossing number
	and the braid index of knots and links},
	Algebr. Geom. Topol. {\bf 6} (2006), 2313--2350. 
\bibitem[Ki]{Kidwell} M.~E.~Kidwell, \em{On the degree of the
	Brandt-Lickorish-Millett-Ho polynomial of a link},
	Proc. Amer. Math. Soc. {\bf 100(4)} (1987), 755--762.
\bibitem[Kn]{Kneissler} J. A. Kneissler, \em{Woven braids and their
	closures}, J. Knot Theory Ramifications {\bf 8(2)} (1999),
	201--214.
\bibitem[Ko]{Kobayashi} T.~Kobayashi, \em{Uniqueness of minimal
	genus Seifert surfaces for links},
	Topology Appl. {\bf 33(3)} (1989), 265--279.
\bibitem[Kr]{Kreimer} D.\ Kreimer, \em{Knots and Feynman diagrams},
	Cambridge Lecture Notes in Physics {\bf 13},
	Cambridge University Press, Cambridge, 2000.
\bibitem[LM]{LickMil} W. B. R. Lickorish and K.~C.~Millett, {\em A
	polynomial invariant for oriented links}, Topology
	{\bf 26 (1)} (1987), 107--141.
\bibitem[LT]{LickThis} \bysame\, and M.~B.~Thistlethwaite,
	\em{Some links with non-trivial polynomials and their crossing
	numbers}, Comment. Math. Helv. {\bf 63} (1988), 527--539.
\bibitem[MT]{MenThis} W.~W.~Menasco and M.~B.~Thistlethwaite, \em{%
  	The Tait flyping conjecture}, Bull. Amer. Math. Soc. {\bf
    	25 (2)} (1991), 403--412.
\bibitem[Mo]{Morton} H.~R.~Morton, \em{Seifert circles and knot
	polynomials}, Proc. Camb. Phil. Soc. {\bf 99} (1986), 107--109.
\bibitem[Mo2]{MortonPb} \bysame\ (ed.), {\em Problems},  in
	``Braids'', Santa Cruz, 1986 (J.~S.~Birman and A.~L.~Libgober,
	eds.), Contemp. Math. {\bf 78}, 557--574.
\bibitem[MS]{MorSho}  \bysame\ and H. B. Short, \em{The
	$2$-variable polynomial of cable knots}, Math. Proc.
	Cambridge Philos. Soc. {\bf 101(2)} (1987), 267--278.
\bibitem[MS2]{MorSho2} \bysame\ and \bysame, {\tt br9z.p}, a Pascal
	program for calculation of the skein polynomial from braids,
	\web|http://www.liv.ac.uk/~su14/knotprogs.html|.
\bibitem[Mr]{Murakami} J.\ Murakami, \em{The Kauffman polynomial of
	links and representation theory},
	Osaka J. Math. {\bf 24(4)} (1987), 745--758.
\bibitem[Mu]{Murasugi} K.~Murasugi, \em{On the braid index of
	alternating links}, Trans. Amer. Math. Soc. {\bf 326 (1)}
	(1991), 237--260.
\bibitem[Mu2]{Murasugi2} \bysame, \em{On closed 3-braids},
  	Memoirs AMS {\bf 151} (1974), AMS, Providence.
\bibitem[Mu3]{Murasugi3} \bysame\,, \em{Jones polynomial and classical
	conjectures in knot theory}, Topology {\bf 26} (1987), 187--194.
\bibitem[MP]{MurPrz} \bysame\ and J.~Przytycki, \em{The skein
	polynomial of a planar star product of two
	links}, Math. Proc. Cambridge Philos. Soc. {\bf 106(2)}
	(1989), 273--276.
\bibitem[MP2]{MurPrz2} \bysame\ and \bysame, \em{An index of a graph
	with applications to knot theory}, Mem. Amer. Math. Soc.
	{\bf 106 (508)} (1993).
\bibitem[Na]{Nakamura} T.~Nakamura, \em{Notes on the braid index of
	closed positive braids}, Topology Appl. {\bf 135(1-3)}
	(2004), 13--31.
\bibitem[Na2]{Nakamura2} \bysame, \em{Braidzel surfaces and the
	Alexander polynomial}, Proceedings of the Workshop
	``Intelligence of Low Dimensional Topology'', Osaka
	City University (2004), 25--34.
\bibitem[Ni]{Ni} Y.~Ni, \em{Closed 3-braids are nearly fibred},
	preprint \web|math.GT/0510243|.
\bibitem[Oh]{Ohyama} Y.~Ohyama, \em{On the minimal crossing number
	and the braid index of links}, Canad. J. Math. {\bf 45(1)}
	(1993), 117--131.
\bibitem[Or]{Orevkov} S.~Orevkov, \em{Quasipositivity problem for
	3-braids}, Turkish Journal of Math. {\bf 28} (2004),
	89--93; also available at
	\web|http://picard.ups-tlse.fr/homepage/orevkov.html|.
\bibitem[PV]{VirPol} M.~Polyak and O.~Viro,
	{\em Gauss diagram formulas for Vassiliev invariants},
	Int.\ Math.\ Res.\  Notes {\bf 11} (1994), 445--454.
\bibitem[PV2]{VirPol2} \bysame\, and \bysame, {\em On the Casson
	knot invariant}, Knots in Hellas '98, Vol. {\bf 3} (Delphi),
	J. Knot Theory Ramifications {\bf 10(5)} (2001), 711--738.
\bibitem[Ro]{Rolfsen} D.~Rolfsen, {\em Knots and links}, Publish
	or Perish, 1976.
\bibitem[Ru]{Rudolph} L.\ Rudolph, \em{Braided surfaces and
	Seifert ribbons for closed braids}, Comment. Math. Helv.
	{\bf 58} (1983), 1--37.
\bibitem[Ru2]{Rudolph2} \bysame, \em{Quasipositivity as an
	obstruction to sliceness}, Bull. Amer. Math. Soc. (N.S.)
	{\bf 29(1)} (1993), 51--59.
\bibitem[Sc]{Sc} O.~Schreier, \em{\"Uber die Gruppen $A^aB^b=1$},
	Abh. Math. Sem. Univ. Hamburg {\bf 3} (1924), 167--169.
\bibitem[SSW]{SSW} D.~Silver, A.~Stoimenow and S.~G.~Williams, \em{%
	Euclidean Mahler measure and Twisted Links},
	\web|math.GT/0412513|,
	Algebr. Geom. Topol. {\bf 6} (2006), 581--602.
\bibitem[SW]{SW} \bysame\, and S.~Williams, \em{Coloring link diagrams
	with a continuous palette}, Topology {\bf 39} (2000),
	1225--1237.
\bibitem[Sq]{Squier} C.\ Squier, \em{The Burau representation is
  	unitary}, Proc.\ Amer.\ Math.~Soc. {\bf 90} (1984), 199--202.
\bibitem[St]{posqp} A.~Stoimenow, {\em On polynomials and surfaces of
	variously positive links}, \web|math.GT/0202226|,
	Jour. Europ. Math. Soc. {\bf 7(4)} (2005), 477--509.
\bibitem[St2]{3br} \bysame, {\em The skein polynomial of closed
	3-braids}, J. Reine Angew. Math. {\bf 564} (2003), 167--180.
\bibitem[St3]{posex_bcr} \bysame, {\em On the crossing number of
   	positive knots and braids and braid index criteria of Jones
	and Morton-Williams-Franks},
        Trans.\ Amer.\ Math.\ Soc. {\bf 354(10)} (2002), 3927--3954.
\bibitem[St4]{ntriv} \bysame, {\em Coefficients and non-triviality
	of the Jones polynomial}, preprint {\tt math.GT/0606255}.
\bibitem[St5]{gener} \bysame, {\em Diagram genus, generators and
	applications}, preprint.
\bibitem[St6]{bind} \bysame, {\em The braid index and the growth of
	Vassiliev invariants}, J.~Of Knot Theory and Its Ram.
	{\bf 8(6)} (1999), 799--813.
\bibitem[St7]{pos} \bysame, {\em Positive knots, closed braids,
	and the Jones polynomial}, {\tt math.GT/9805078}, Ann. Scuola
	Norm. Sup. Pisa Cl. Sci. {\bf 2(2)} (2003), 237--285.
\bibitem[St8]{gen2}\bysame, {\em Knots of genus two}, preprint
	\web|math.GT/0303012|.
\bibitem[SV]{SV} \bysame\ and A.\ Vdovina, {\em Counting alternating
  	knots by genus}, Math.\ Ann. {\bf 333} (2005), 1--27.
\bibitem[ST]{SunThis} C.~Sundberg and M.~B.~Thistlethwaite, {\em The
	rate of growth of the number of prime alternating links and
	tangles}, Pacific Journal of Math. {\bf 182 (2)} (1998),
	329--358.
\bibitem[Th]{Thistle} M.~B.~Thistlethwaite, {\em On the Kauffman
	polynomial of an adequate link}, Invent.\ Math. {\bf 93(2)}
	(1988), 285--296.
\bibitem[Th2]{Thistle2} \bysame, {\em A spanning tree expansion for the
	Jones polynomial}, Topology {\bf 26} (1987), 297--309.
\bibitem[Tr]{Traczyk} P.~Traczyk, \em{$3$-braids with proportional
  	Jones polynomials}, Kobe J. Math. {\bf 15(2)} (1998), 187--190.
\bibitem[Vo]{Vogel} P.~Vogel, \em{Representation of links by braids:
	A new algorithm}, Comment. Math. Helv. {\bf 65} (1990),
	104--113.
\bibitem[Wi]{Williams} R. F. Williams, \em{Lorenz knots are prime},
	Ergodic Theory Dynam. Systems  {\bf 4(1)} (1984), 147--163. 
\bibitem[Wo]{Wolfram} S. Wolfram, {\em Mathematica --- a system
	for doing mathematics by computer}, Addison-Wesley, 1989.
\bibitem[Xu]{Xu} Peijun Xu, {\em The genus of closed $3$-braids},
	J. Knot Theory Ramifications {\bf 1(3)} (1992), 303--326.
\bibitem[Y]{Yamada} S.~Yamada, \em{The minimal number of Seifert
	circles equals the braid index}, Invent. Math. {\bf 88}
	(1987), 347--356.
\bibitem[Yo]{Yokota} Y.\ Yokota, {\em Polynomial invariants of
	positive links}, Topology {\bf 31(4)} (1992), 805--811.
\end{thebibliography}
}
\end{document}